\documentclass[12pt]{report}
\usepackage{iftex}
\ifpTeX
  \PassOptionsToPackage{dvipdfmx}{graphicx}
  \PassOptionsToPackage{dvipdfmx}{xcolor}
\fi
\usepackage{graphicx}
\usepackage{float}
\usepackage{amsmath}
\usepackage{amssymb}
\usepackage{amsthm}
\usepackage{pb-diagram}
\usepackage{mathtools}
\usepackage{bm}
\usepackage{mathrsfs}
\usepackage{enumerate}
\usepackage[inline]{enumitem}
\usepackage{appendix}
\usepackage{titlesec}

\titleformat{\section}
  {\large\bfseries}
  {\thesection}
  {0.8em}
  {}

\titleformat{\subsection}
  {\normalsize\bfseries}
  {\thesubsection}
  {0.8em}
  {}

\titlespacing*{\section}
  {0pt}{2.5ex plus 1ex minus .2ex}{1.2ex}

\titlespacing*{\subsection}
  {0pt}{2ex plus .8ex minus .2ex}{0.8ex}
\usepackage{libertinus}
\usepackage{tikz-cd}
\usepackage{fancyhdr}
\usepackage[top=20mm,bottom=25mm,left=25mm,right=25mm]{geometry}
\newtheoremstyle{underlined}
  {6pt}
  {6pt}
  {\normalfont}
  {}
  {\bfseries}
  {.}
  {0.5em}
  {\thmname{\underline{#1\ #2}}\thmnote{ (#3)}}

\theoremstyle{underlined}
\newtheorem{dfn}{Definition}[section]
\newtheorem{thm}[dfn]{Theorem}
\newtheorem{prp}[dfn]{Proposition}
\newtheorem{lem}[dfn]{Lemma}
\newtheorem{cor}[dfn]{Corollary}
\newtheorem{ex}[dfn]{Example}
\newtheorem*{lem*}{Lemma}
\newtheorem*{rem}{Remark}
\def\<{\langle}
\def\>{\rangle}

\def\L{\text{L}}

\def\D{\text{D}}

\def\S{\text{S}}

\def\c{\text{c}}
\begin{document}
\title{
{\fontsize{30}{36}\selectfont Stochastic Parallel Transport along Markov Processes\par}

\vspace{2em}

{\LARGE
I. General Theory\par
}
}

\author{
{\LARGE Kenta Misumi}}
\date{}
\setlength{\abovedisplayskip}{4pt}
\setlength{\belowdisplayskip}{4pt}
\makeatletter
\begin{titlepage}
\centering
\vspace*{\fill}
\@title\par
\vfill
{\@author\par}
\vfill
\begin{minipage}{0.9\textwidth}
\fontsize{16}{23}\selectfont
{\centering\fontsize{16}{22}\selectfont\bfseries Abstract\par}
\vspace{0.8em}
\noindent
We give a framework for transporting vectors along Markov processes that may have jumps, without assuming a manifold or other specific geometric structure.
We prove existence and uniqueness of the resulting transport. In doing so, we develop a theory of vector bundles defined quasi-everywhere, in which the vectors live, and an axiomatic framework for energy forms that determine the bundle and the transport.
\end{minipage}\par
\vspace*{\fill}
\end{titlepage}
\makeatother
\chapter*{Introduction}
\renewcommand{\thefigure}{\arabic{figure}}
\markboth{Introduction}{Introduction}
The aim of this paper is to establish a framework for transporting vectors along stochastic processes that may have jumps, without assuming a manifold or other specific geometric structure, and to prove existence and uniqueness of the resulting transport.

\section*{Why is transporting vectors important?}
The reason is that transporting vectors can detect the curvature of a space. That is, after a vector is transported around a closed loop, it need not return to its original direction. Two elementary examples illustrate the idea; no stochastic process is needed here.

\begin{figure}[htbp]
    \centering
    \includegraphics[width=0.85\textwidth]{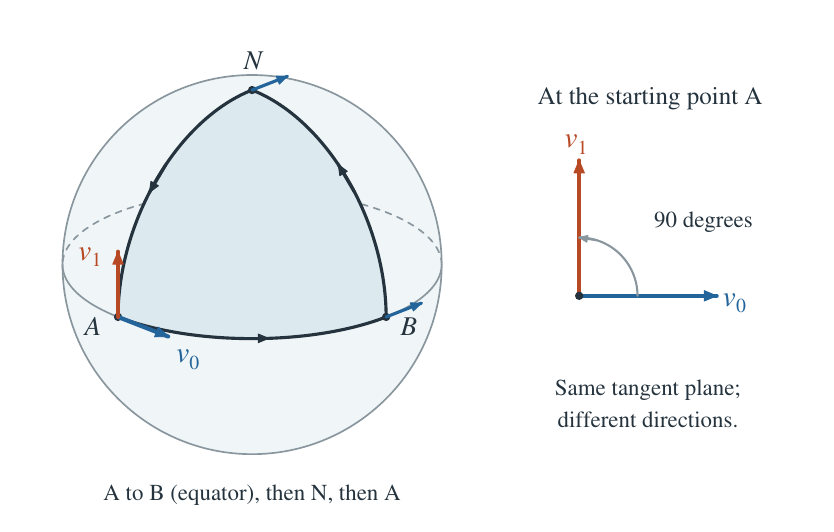}
    \caption{Parallel transport along $A\to B\to N\to A$ on the sphere.
    The initial vector $v_0$ and the returned vector $v_1$ differ by $\pi/2$.
    The right-hand diagram compares them in the same tangent plane at $A$.}
    \label{fig:intro-sphere-transport}
\end{figure}

\medskip
\noindent\textbf{The sphere.}
On the unit sphere, consider the loop that follows the equator eastward for a quarter turn from $A$ to $B$, continues along a meridian to the north pole $N$, and returns to $A$ along another meridian (Figure~\ref{fig:intro-sphere-transport}).

\begin{figure}[htbp]
    \centering
    \includegraphics[width=\textwidth]{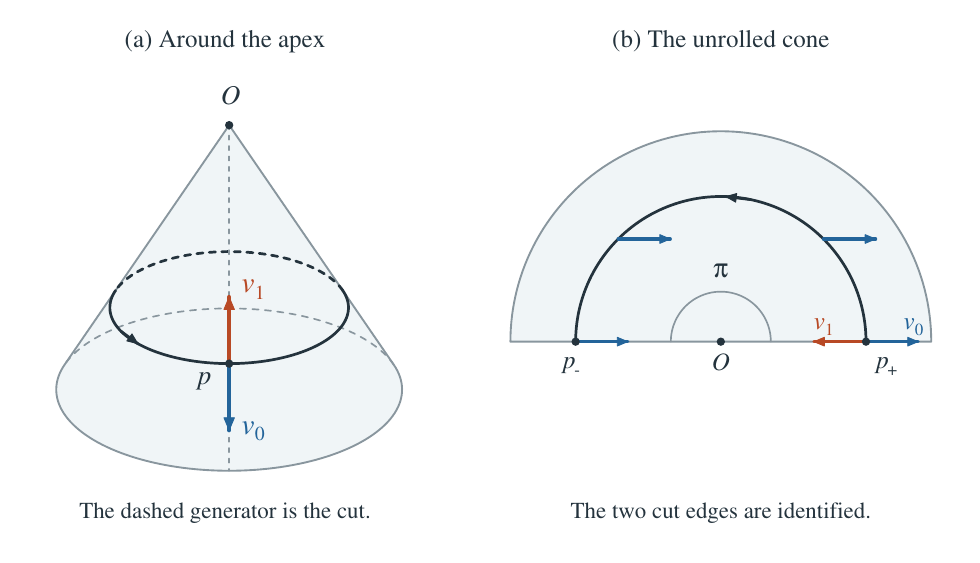}
    \caption{A loop around the apex of a cone of total angle $\pi$.
    The transported vector points to the right throughout the unfolded semicircle.
    Identifying the cut edges rotates the vector by $\pi$,
    so the returned vector is $v_1=-v_0$.}
    \label{fig:intro-cone-transport}
\end{figure}

\medskip
\noindent\textbf{The cone.}
Consider a circular cone of total angle $\pi$ and move once around its apex.
Cutting along a generator unfolds the cone into a semicircle
(Figure~\ref{fig:intro-cone-transport}).
In the unfolded picture, a vector pointing to the right at $p_+$
keeps pointing to the right all the way to $p_-$.
It points away from the apex at $p_+$ and toward it at $p_-$.
The two endpoints represent the same point on the cone, so after one
circuit the vector returns pointing in the opposite direction, giving $v_1=-v_0$.
The cone is intrinsically flat away from its apex. This rotation by $\pi$
detects the curvature concentrated at the apex.

\section*{Why stochastic?}
Both examples above involve transport along deterministic curves, so they do not yet explain why one should consider transport along stochastic paths. A key motivation for doing so, particularly along Markov processes,
is the following idea, which captures the spirit of this paper:

\begingroup
\setlength{\fboxsep}{10pt}
\setlength{\fboxrule}{0.6pt}
\begin{equation}
\tag{\(\star\)}\label{principle:stochastic-transport-heat}
\fbox{%
    \begin{minipage}{0.84\linewidth}
    \centering\large\bfseries
    Stochastic parallel transport is the probabilistic counterpart
    of the heat flow on one-forms.
    \par
    \end{minipage}%
}
\end{equation}
\endgroup

\begin{figure}[H]
    \centering
    \includegraphics[width=0.95\textwidth]{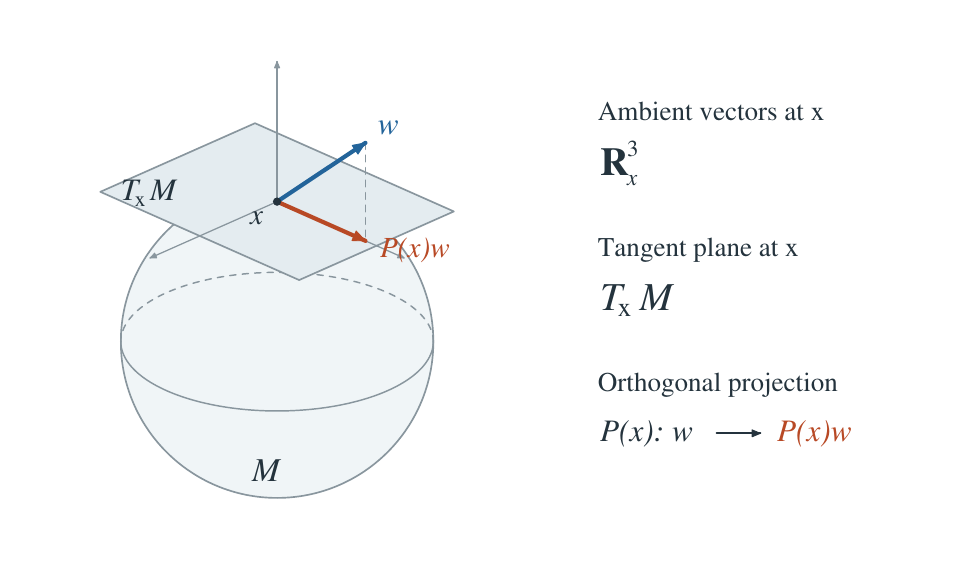}
    \caption{Orthogonal projection onto the tangent plane of a sphere in
    $\mathbb R^3$. An ambient vector $w\in\mathbb R_x^3$ is projected onto
    $P(x)w\in T_xM$.}
    \label{fig:intro-tangent-projection}
\end{figure}

This is precisely analogous to the role of Brownian motion as the probabilistic counterpart of the heat flow on functions.
We now describe this explicitly.
The following discussion may be skipped on a first reading.
Let $M$ be a smooth, compact embedded submanifold of $\mathbb R^N$ without
boundary, equipped with the induced Riemannian metric $g$.
For $x\in M$, write $\mathbb R_x^N$ for the space of ambient vectors based
at $x$, and $T_xM$ for the tangent space to $M$ at $x$.
Let
\[
    P(x)\colon \mathbb R_x^N\longrightarrow T_xM
\]
denote the orthogonal projection
(Figure~\ref{fig:intro-tangent-projection}).

\begin{figure}[H]
    \centering
    \includegraphics[width=0.95\textwidth]{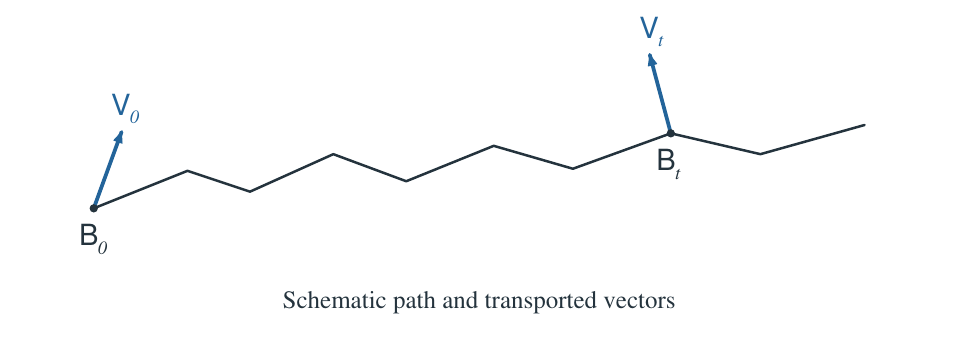}
    \caption{A schematic Brownian path $\mathsf B$ and tangent vectors
    $\mathsf V_t\in T_{\mathsf B_t}M$ based along it.}
    \label{fig:intro-brownian-vector-path}
\end{figure}

Let $(\mathsf B_t)$ be Brownian motion on $M$, started at $x$.
Let $(\mathsf V_t)$ be an $\mathbb R^N$-valued continuous semimartingale
satisfying $\mathsf V_t\in T_{\mathsf B_t}M$.
Figure~\ref{fig:intro-brownian-vector-path} illustrates these processes.

We define the covariant derivative of $\mathsf V$ along
$\mathsf B$ by
\[
    D_{\circ\,d\mathsf B_t}\mathsf V_t
    :=P(\mathsf B_t)\circ d\mathsf V_t.
\]
We say that $\mathsf V$ is parallel along $\mathsf B$ if
\[
    D_{\circ\,d\mathsf B_t}\mathsf V_t=0. \tag{$*$}
\]
Together with the tangency contition
$P(\mathsf B_t)\mathsf V_t=\mathsf V_t$,
$(*)$ gives a linear
Stratonovich SDE with smooth coefficients.
Thus, for each prescribed initial vector $v\in T_xM$, there is a unique $(\mathsf V_t)_{t\geq0}$ with $\mathsf V_0=v$.
This process is called the stochastic parallel transport of $v$ along $(\mathsf B_t)_{t\geq 0}$.
Now let $\alpha$ be any smooth one-form on $M$.
The Stratonovich product rule, parallelism of $\mathsf V$, and It\^o-Stratonovich conversion rule give (the reader can verify the following by a formal calculation)
\[
\begin{aligned}
    d\langle\alpha(\mathsf B_t),\mathsf V_t\rangle
    &=\big\langle D_{\circ\,d\mathsf B_t}(\alpha(\mathsf B_t)),
        \mathsf V_t\big\rangle
      +\big\langle\alpha(\mathsf B_t),
        D_{\circ\,d\mathsf B_t}\mathsf V_t\big\rangle\\
    &=\langle(D\alpha)(\mathsf B_t),\mathsf V_t\rangle\,d\mathsf B_t
      +\frac12\big\langle(D^2\alpha)(\mathsf B_t)
        (d\mathsf B_t\otimes d\mathsf B_t),\mathsf V_t\big\rangle\\
    &=\langle(D\alpha)(\mathsf B_t),\mathsf V_t\rangle\,d\mathsf B_t
      -\frac12\langle(\Delta^C\alpha)(\mathsf B_t),\mathsf V_t\rangle\,dt,
\end{aligned}
\]
where the formal rule
$d\mathsf B_t\otimes d\mathsf B_t=g^{-1}(\mathsf B_t)\,dt$
contracts $D^2\alpha$ to its trace $\Delta^C\alpha:=-\operatorname{tr}_gD^2\alpha$. ($\Delta^C$ is called the connection Laplacian.)
Define the parallel transport operator by
\[
    \mathsf U_t(\omega)\colon
    T_{\mathsf B_0(\omega)}M\longrightarrow T_{\mathsf B_t(\omega)}M,
    \qquad
    v\longmapsto\mathsf V_t(\omega).
\]
Then we can write the equation by
\[d(\mathsf U_t^* \alpha(B_t))=\mathsf U_t^*((D\alpha )(\mathsf B_t))\,d\mathsf B_t-\frac{1}{2}\mathsf U_t^*((\Delta^C\alpha)(\mathsf B_t))\,dt.\]
Taking expectation,
\[
    (T_t^C\alpha)(x)
    =
    \mathbb E_x\!\left[\mathsf U_t^*\alpha(\mathsf B_t)\right],
    \qquad
    T_t^C=e^{-t\Delta^C/2}.
\]
Thus the connection heat flow is obtained by averaging one-forms transported back to the starting point.

\section*{How can we formulate stochastic parallel transport without assuming a manifold structure?}
Our setting consists of a quasi-regular Dirichlet space $(X,m,\mathcal E)$ and an associated right process $(\mathsf X_t,\mathbb P_x)$;
see \cite{MaRockner1992,FukushimaOshimaTakeda2011} for standard references.
This framework includes Dirichlet forms on smooth Riemannian and
sub-Riemannian manifolds
\cite{Hsu2002,BaudoinBonnefontGarofaloMunive2014}, as well as on
infinite-dimensional spaces \cite{MaRockner1992},
fractals \cite{Kigami2001}, and RCD spaces
\cite{AmbrosioGigliSavare2014}.
In this general setting, neither a manifold structure nor a notion of
tangent space is given a priori.
In the spirit of principle~\eqref{principle:stochastic-transport-heat},
our approach is to axiomatize the energy of one-forms on $(X,m,\mathcal E)$
and to construct parallel transport as the corresponding stochastic object.
We must address three problems:
\begin{enumerate}
    \item\label{prob:intro-vector-spaces}
    Define a structure, analogous to the tangent bundle,
    in which vectors can live.
    \item\label{prob:intro-one-form-energy}
    Axiomatize the notion of energy for one-forms.
    \item\label{prob:intro-transport-existence-uniqueness}
    Formulate the existence and uniqueness problem for stochastic
    parallel transport using these structures.
\end{enumerate}

Problem~\ref{prob:intro-vector-spaces} is addressed in
Chapter~\ref{chap:quasi-hilbert-bundles} by introducing the notion of a
\emph{quasi-Hilbert bundle}.
Informally, this is a framework for vector bundles that disregards
sets of zero capacity, which are invisible to the associated Markov process.
For smooth real metric bundles of the same rank,
Theorem~\ref{thm:smooth-bundle-quasi-classification} gives the following
characterization:

\begin{center}
\begingroup
\setlength{\fboxsep}{10pt}
\setlength{\fboxrule}{0.6pt}
\fbox{%
    \begin{minipage}{\dimexpr\linewidth-2\fboxsep-2\fboxrule\relax}
    \centering\large\bfseries
    Two smooth vector bundles are isomorphic as quasi-Hilbert bundles precisely when
    they are smoothly isomorphic outside a set of zero capacity.
    \par
    \end{minipage}%
}
\endgroup
\end{center}

At the level of almost Hilbert bundles, which disregard sets of measure zero,
the only invariant of a smooth vector bundle of constant rank is its rank
(Proposition~\ref{prp:almost-triviality}).
By contrast, the quasi-Hilbert bundle framework distinguishes the trivial real
line bundle over $S^1$ from the M\"obius line bundle.
However, it does not distinguish all smooth vector bundles: for example,
it identifies the tangent bundle $TS^2$ with the trivial rank-two bundle
$S^2\times\mathbb R^2$.
Indeed, the north pole $p$ has zero capacity, and $TS^2$ is trivial on
$S^2\setminus\{p\}\cong\mathbb R^2$.
Theorem~\ref{thm:smooth-bundle-quasi-classification} also shows that the quasi-isomorphism class of a smooth real metric bundle of fixed rank is determined by its determinant line, or equivalently by its first Stiefel--Whitney class in $H^1(M;\mathbb Z/2\mathbb Z)$.
The underlying Dirichlet space alone does not specify which
quasi-Hilbert bundle should carry the transport.
This brings us to Problem~\ref{prob:intro-one-form-energy}.
In Chapter~\ref{chap:bundle-forms-quasi-realization}, we axiomatize energies
of one-forms through the notion of a \emph{Kato form}
(Definition~\ref{dfn:kato-form}).
Its defining property is equivalent to domination of the associated
semigroup by the scalar heat semigroup
(Theorem~\ref{thm:kato-domination}).
This connects the definition with classical Kato inequalities and
semigroup domination; see
\cite{HessSchraderUhlenbrock1977,Ouhabaz1996,LenzSchmidtWirth2020}.
The main result of the chapter is the \emph{quasi-realization theorem}
(Theorems~\ref{thm:quasi-realization-existence}
and~\ref{thm:quasi-realization-uniqueness}).
Informally, it states:

\begin{center}
\begingroup
\setlength{\fboxsep}{10pt}
\setlength{\fboxrule}{0.6pt}
\fbox{%
    \begin{minipage}{\dimexpr\linewidth-2\fboxsep-2\fboxrule\relax}
    \centering\large\bfseries
    Every Kato form admits a natural realization on a quasi-Hilbert
    bundle, unique up to isomorphism.
    \par
    \end{minipage}%
}
\endgroup
\end{center}

On this bundle, sections of finite energy have quasi-continuous
representatives.
In this sense, the result may be viewed as a quasi-analogue of the
Serre--Swan theorem \cite{Swan1962}: the energy determines the bundle
on which its sections live.
For example, let $M$ be a smooth Riemannian manifold with smooth boundary,
equipped with the scalar Neumann Dirichlet form.
We use the metric to identify one-forms with tangent vectors.
The connection energy with Neumann boundary conditions is a Kato form,
and its quasi-realization is the full tangent bundle $TM$.
In particular, at a boundary point $x$, the fiber is the whole space
$T_xM$, including the normal direction.
Now consider the connection energy with mixed boundary conditions:
the normal component vanishes at the boundary, while the tangential
components satisfy the natural Neumann condition.
This is also a Kato form, but its quasi-realization has fibers
\[
    (T_{\mathrm{tan}}M)_x=
    \begin{cases}
        T_xM,&x\in M^\circ,\\
        T_x\partial M,&x\in\partial M.
    \end{cases}
\]
Thus the normal direction disappears from the boundary fibers.
These two examples are treated in
Sections~\ref{subsec:neumann-connection-energy}
and~\ref{subsec:mixed-connection-energy}.
Since the boundary has positive capacity, our framework can distinguish
between these two cases.
Problem~\ref{prob:intro-transport-existence-uniqueness} is addressed in
Chapter~\ref{chap:multiplicative-stochastic-transports}.
The main result of this chapter, and of the present paper, is the
following existence and uniqueness result
(Theorems~\ref{thm:multiplicative-transport-existence}
and~\ref{thm:multiplicative-transport-uniqueness}).
Informally, it states:
\begin{center}
\begingroup
\setlength{\fboxsep}{10pt}
\setlength{\fboxrule}{0.6pt}
\fbox{%
    \begin{minipage}{\dimexpr\linewidth-2\fboxsep-2\fboxrule\relax}
    \centering\large\bfseries
    Every Kato form admits a unique associated stochastic transport
    on its quasi-realization.
    \par
    \end{minipage}%
}
\endgroup
\end{center}

\begin{figure}[H]
    \centering
    \includegraphics[width=\linewidth]{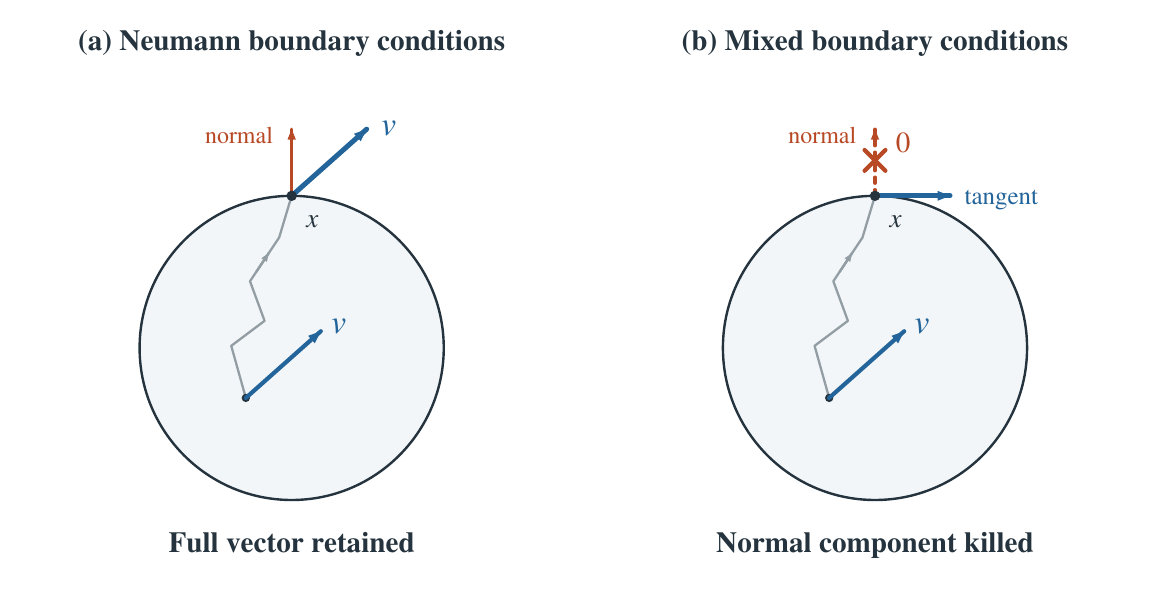}
    \caption{Transport at the first boundary hit on a Euclidean disk.
    The gray curves indicate schematic paths to a boundary point $x$.
    Left: transport with Neumann condition retains the full vector.
    Right: transport with mixed condition kills the normal component (crossed out)
    and retains the tangential component.}
    \label{fig:boundary-transport-disk}
\end{figure}

Returning to the preceding examples, the Neumann connection energy
gives a transport that retains both tangential and normal components
at the boundary.
The mixed connection energy gives a transport that kills the normal
component there, leaving only the tangential component.
On a Euclidean disk, at the first boundary hit, the former leaves the
vector unchanged, while the latter projects it onto the boundary
tangent line (Figure~\ref{fig:boundary-transport-disk}).
In particular, a purely normal vector is sent to zero in the mixed case.
For parallel transport along reflecting Brownian motion, see
\cite[Section~2]{Hsu2002Boundary}; for the removal of the normal
component at the boundary, see \cite[Section~3]{Hsu2002Boundary}
and \cite{ArnaudonLi2017}.

\begin{figure}[H]
    \centering
    \includegraphics[width=\linewidth]{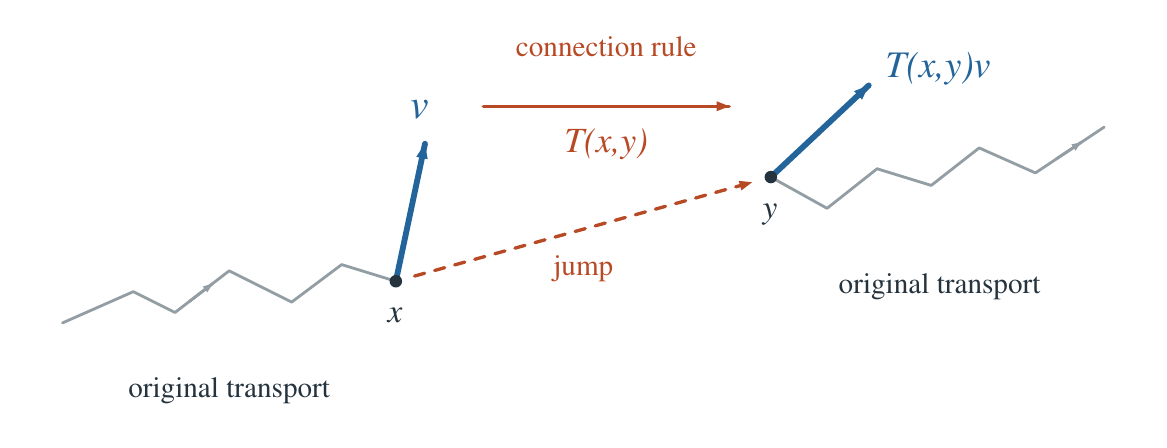}
    \caption{A connection rule at a jump.
    The dashed arrow indicates a jump from $x$ to $y$,
    at which the vector $v$ is replaced by $T(x,y)v$.
    Between jumps, vectors follow the original transport.}
    \label{fig:jump-connection-rule}
\end{figure}

Another important energy on one-forms is the Hodge energy.
By the Weitzenb\"ock formula, it is a potential perturbation of the
connection energy, involving Ricci curvature and, when a boundary is
present, the second fundamental form.
We treat this relation in Section~\ref{subsec:hodge-energy-mixed-boundary},
within the potential-perturbation framework of
Section~\ref{subsec:potential-perturbation}.
In Section~\ref{subsec:damped-transports}, we also establish existence
and uniqueness of the corresponding transport, called
\emph{damped transport}
(Theorems~\ref{thm:damped-transport-existence}
and~\ref{thm:damped-transport-association}).
To identify the associated form, we use the resolvent formula for vector measures of finite energy integral (Theorem~\ref{thm:vector-revuz-resolvent}).
In the smooth setting, this recovers the classical damped parallel
transport for the Hodge heat flow
\cite[Section~7.2]{Hsu2002}, including the boundary constructions of
Hsu~\cite{Hsu2002Boundary} and Arnaudon and Li~\cite{ArnaudonLi2017}.
Jump perturbations provide another important class of energies.
Here we specify a \emph{connection rule} $T(x,y)$: a linear map
describing how a vector is carried across a jump from $x$ to $y$
(Definition~\ref{dfn:connection-rule}).
Together with a jump measure $j$, this rule defines a jump energy which is added to the original energy
(Section~\ref{subsec:jump-perturbation}).
For continuous original paths and bounded jump rates, the corresponding
\emph{jump transport} follows the original transport between jumps
and applies $T(x,y)$ at each jump
(Example~\ref{ex:finite-activity-jump-transport}
in Section~\ref{subsec:jump-transports}).
Thus a jump from $x$ to $y$ sends
\[
    v\longmapsto T(x,y)v,
\]
as illustrated in Figure~\ref{fig:jump-connection-rule}.
The analogous construction on weighted graphs, using unitary connection
maps, and its Feynman--Kac representation are studied in
\cite{GuneysuMilatovicTruc2014}.

\section*{Why do we need such a general theory?}
The purpose of this general framework is to develop stochastic parallel
transport and Malliavin calculus on RCD spaces (metric measure spaces
with a synthetic lower Ricci curvature bound and a quadratic Sobolev
energy; see \cite{AmbrosioGigliSavare2014}).
The details will be developed in Part~II of this series.
These developments will build on Gigli's differential calculus on
RCD spaces \cite{Gigli2018}, Braun's vector calculus and semigroup
domination results for tamed Dirichlet spaces \cite{Braun2024},
and the theory of quasi-continuous vector fields of Debin, Gigli,
and Pasqualetto \cite{DebinGigliPasqualetto2021}.
There we plan to study stochastic parallel transport along Brownian 
motion and Brownian bridge on RCD spaces and its relation to parallel 
transport along families of curves
as developed by Gigli, Pasqualetto, and collaborators
\cite{GigliPasqualetto2020,CaputoGigliPasqualetto2025}.
We also plan to discuss Malliavin calculus on RCD path spaces, logarithmic Sobolev inequalities, and characterizations of two-sided Ricci curvature bounds in terms of energy forms on path space.
A guiding precedent is Naber's characterization of two-sided Ricci
curvature bounds in the smooth setting through estimates on path space,
including suitable logarithmic Sobolev inequalities
\cite{Naber2015}; see also \cite{HaslhoferNaber2018}.

\chapter*{Prerequisites}
\markboth{Prerequisites}{Prerequisites}
Fukushima, Oshima and Takeda (FOT) \cite{FukushimaOshimaTakeda2011} is the main reference for the Dirichlet form theory used throughout. On a first reading for this paper, one may assume that the form is regular and take its associated Hunt process (FOT, \S7.2). For Section~1.1, the bundle construction uses Hilbert spaces and $L^2$ integration, and the cotangent bundle uses carré du champ (FOT, \S3.2). For Section~1.2, the scalar notions of capacity, nests and quasi-continuity in FOT, \S2.1, suffice to define quasi-Hilbert bundles. The geometric argument in Section~1.3 uses transversality to the zero section on a manifold with boundary, for which Guillemin and Pollack \cite[Chapter~2, Sections~1 and~3]{GuilleminPollack1974} give a short account. Section~2.1 uses closed forms and their semigroups (FOT, \S1.3) and measures of finite energy integral (FOT, \S2.2). In Section~\ref{sec:bundle-form-operations}, potential perturbations and restrictions use smooth measures (FOT, \S2.2); jump perturbations use jumping measures (FOT, \S3.2). Section~2.3 assumes familiarity with metric connections, differential forms and integration by parts on manifolds with boundary. Section~3.1 uses the Markov property and martingale regularization (FOT, Appendix~A.3). We also use the optional section theorem to deduce indistinguishability from almost sure equality at every bounded stopping time. For the statement, see Cohen and Elliott \cite[Chapter~7, p.~159]{CohenElliott2015}, and for a proof, see Dellacherie and Meyer \cite[Theorem~IV.84]{DellacherieMeyer1978}. For Section~\ref{sec:transport-operations}, the Revuz correspondence is treated in FOT, \S5.1, scalar potential perturbations in FOT, \S6.1, monotone convergence of closed forms in \cite[Section~4]{Simon1978}, and jumps in FOT, \S4.5. For Section~3.3, Hsu \cite[Section~2]{Hsu2002Boundary} treats stochastic parallel transport along reflecting Brownian motion and the role of boundary local time. Readers who wish to work in the full quasi-regular setting can consult Ma and R\"ockner \cite[Chapter~III, Sections~2--3, and Chapter~IV, Sections~2--3]{MaRockner1992} for nests, quasi-continuity and the associated right process.

\chapter*{Basic setting}
\markboth{Basic setting}{Basic setting}
We first fix the setting for the paper. Let $X$ be a Polish space, $m$ a locally finite Borel measure on $X$, and $(\mathcal E,\D(\mathcal E))$ a symmetric Dirichlet form on $\L^2(m)$. We assume that $\mathcal E$ is quasi-regular in the sense of \cite[Chapter~IV, Section~3]{MaRockner1992} and call $(X,m,\mathcal E)$ a quasi-regular Dirichlet space. A nest is an increasing sequence $(F_n)$ of closed subsets of $X$ such that the functions in $\D(\mathcal E)$ that vanish $m$-a.e.\ outside some $F_n$ are dense in $\D(\mathcal E)$ for the form norm $\|u\|_{\mathcal E_1}=(\mathcal E(u,u)+\|u\|_2^2)^{1/2}$. A set is exceptional if it is disjoint from the union of some nest. Every exceptional set is $m$-null, and a countable union of exceptional sets is exceptional. We say that a property holds quasi-everywhere (q.e.) if it holds outside an exceptional set. A scalar function is quasi-continuous if its restriction to each $F_n$ of some nest is continuous. A countable family of quasi-continuous functions admits a common such nest. Quasi-continuous functions that agree $m$-a.e.\ agree q.e., and every $u\in\D(\mathcal E)$ has a quasi-continuous representative $\tilde u$. If $u_n\to u$ in the form norm, then some subsequence of $(\tilde u_n)$ converges to $\tilde u$ q.e. \cite[Chapter~III, Sections~2--3]{MaRockner1992}.
On a measurable space $(\Omega,\mathcal F)$ we fix an $m$-symmetric right process $(\mathsf X_t,\mathbb P_x)$ properly associated with $\mathcal E$ \cite[Chapter~IV, Sections~2--3]{MaRockner1992}. It takes values in $X_\Delta=X\cup\{\Delta\}$, where $\Delta$ is the absorbing cemetery state. Under $\mathbb P_x$ the process starts at $x$. Let $\mathcal F_t^{\mathsf X}=\sigma(\mathsf X_s\mid 0\leq s\leq t)$ be the natural filtration, $\mathcal F^{\mathsf X}=\sigma(\mathsf X_s\mid s\geq0)$ the $\sigma$-field generated by the entire path, and $\mathcal F_{t+}^{\mathsf X}=\bigcap_{u>t}\mathcal F_u^{\mathsf X}$ the right limit at $t$. We write $(\mathcal F_t^{\mathsf X,x})_{t\geq0}$ for the $\mathbb P_x$-augmentation of $(\mathcal F_{t}^{\mathsf X})_{t\geq0}$ obtained by adjoining all subsets of $\mathbb P_x$-null sets in $\mathcal F^{\mathsf X}$. This filtration satisfies the usual conditions: it is complete and right-continuous. The process is strong Markov with respect to $(\mathcal F_{t+}^{\mathsf X})_{t\geq0}$. Its paths are right-continuous. We write $\zeta=\inf\{t\geq0\mid \mathsf X_t=\Delta\}$, with $\inf\varnothing=\infty$, and $\mathsf X_t=\Delta$ for $t\geq\zeta$. The time shift $\theta_s$ satisfies $\mathsf X_t\circ\theta_s=\mathsf X_{s+t}$. Proper association means that for each bounded Borel $f\in\L^2(m)$ and $t>0$, the function $x\mapsto\mathbb E_x[f(\mathsf X_t)]$, with $f(\Delta)=0$, is a quasi-continuous $m$-version of $T_tf$. For an increasing sequence $(F_n)$ of closed subsets of $X$, let $\tau_n=\inf\{t\geq0\mid \mathsf X_t\notin F_n\}$. Then $(F_n)$ is a nest if and only if $\mathbb P_x(\lim_n\tau_n=\zeta)=1$ for q.e.\ $x$ \cite[Chapter~IV, Section~5]{MaRockner1992}. Thus along almost every path from q.e.\ $x$, each interval $[0,t]$ with $t<\zeta$ lies before some $\tau_n$. In particular, for each exceptional set $N$, the path avoids $N$ throughout $[0,\zeta)$ almost surely for q.e.\ starting point, and a quasi-continuous function evaluated along the path is almost surely right-continuous on $[0,\zeta)$ for q.e.\ starting point.

\chapter*{Notation}
\markboth{Notation}{Notation}
\begingroup
\small
\renewcommand{\arraystretch}{1.15}
\noindent\begin{tabular}{@{}p{0.30\textwidth}@{\hspace{0.03\textwidth}}p{0.67\textwidth}@{}}
\textbf{Notation} & \textbf{Meaning} \\
\hline
$m$-a.e., q.e. & Outside an $m$-null set or an exceptional set, respectively. \\
$\mathcal H$, $\widetilde{\mathcal H}$ & An almost Hilbert bundle over $(X,m)$ and a quasi-Hilbert bundle over $(X,m,\mathcal E)$, respectively. \\
$\L^p(m)$ & Almost measurable real functions modulo $m$-a.e.\ equality, with $\int_X |f|^p\,dm<\infty$ for $1\le p<\infty$. For $p=\infty$ we require essential boundedness, and for $p=0$ we impose no integrability condition. \\
$\L^p(m;\mathcal H)$ & Almost measurable sections $\xi$ of $\mathcal H$ with $|\xi|\in\L^p(m)$, modulo $m$-a.e.\ equality. \\
$\tilde\alpha$, $\Gamma^{\mathrm{qc}}(\widetilde{\mathcal H})$ & A quasi-continuous representative of $\alpha$, and the space of quasi-continuous sections of $\widetilde{\mathcal H}$. \\
$\mathcal E$, $L$, $T_t$, $G_\lambda$ & The scalar Dirichlet form on $\L^2(m)$, its generator $L$, its semigroup $T_t=e^{tL}$ and its resolvent $G_\lambda=(\lambda-L)^{-1}$ for $\lambda>0$. \\
$\breve{\mathcal E}$, $\breve L$, $\breve T_t$, $\breve G_\lambda$ & A bundle form on $\L^2(m;\mathcal H)$ and its associated generator, semigroup and resolvent. \\
$(\mathcal H^{\breve{\mathcal E}},J)$ & A quasi-realization of $\breve{\mathcal E}$. The map $J$ identifies its fibers with those of the original bundle $m$-a.e. \\
$\mathsf X$, $\mathsf B$, $\mathsf U$, $\mathsf V$ & Sans-serif letters denote stochastic processes, including the underlying Markov process, Brownian motion, stochastic transport and a transported vector, respectively. \\
$\mathbb P_x$, $\mathbb E_x$, $\zeta$ & Law and expectation for the process started at $x$, and its lifetime. \\
$\mathcal F_t^{\mathsf X}$, $\mathcal F_{t+}^{\mathsf X}$, $\mathcal F^{\mathsf X}$ & The natural filtration, its right limit, and the $\sigma$-field generated by the entire path, respectively. \\
$\mathcal F_t^{\mathsf X,x}$ & The usual $\mathbb P_x$-augmentation of $(\mathcal F_{t}^{\mathsf X})_{t\geq0}$. \\
$A^*$, $A^\dagger$ & The dual operator and the Hilbert space adjoint of $A$, respectively. \\
$\mathsf U_t$, $\mathsf U_t^*$ & Transport from the dual fiber at $\mathsf X_0$ to the dual fiber at $\mathsf X_t$, and its dual operator acting in the reverse direction on the primal fibers. \\
$\theta_s$, $\mathsf U_{s,t}$ & The path shift by $s$ and transport from time $s$ to time $t$, given by $\mathsf U_{s,t}=\mathsf U_{t-s}\circ\theta_s$. \\
$j$, $T(x,y)$ & A jump measure and the contraction from the dual fiber at $x$ to that at $y$ assigned by the connection rule. \\
$R\mu$, $A^\mu$ & A smooth potential with operator density $R$ and the positive continuous additive functional with Revuz measure $\mu$. \\
\end{tabular}
\endgroup

\clearpage
\renewcommand{\thefigure}{\thechapter.\arabic{figure}}
\tableofcontents
\clearpage
\chapter{Quasi-Hilbert Bundles}\label{chap:quasi-hilbert-bundles}
\section*{Overview}
This chapter introduces two notions of a Hilbert bundle over a quasi-regular Dirichlet space $(X,m,\mathcal E)$. An almost Hilbert bundle is described by measurable sections, with fibers identified up to sets of $m$-measure zero. A quasi-Hilbert bundle is described by quasi-continuous sections, with fibers identified up to sets of zero capacity. In particular, when $\mathcal E$ admits a carr\'e du champ, we show that it canonically determines an almost cotangent bundle (Theorem~\ref{thm:almost-cotangent-bundle}). Section~\ref{sec:riemannian-bundles} shows what this distinction means on a Riemannian manifold. Every smooth metric vector bundle of constant rank is trivial as an almost Hilbert bundle (Proposition~\ref{prp:almost-triviality}). Two smooth real metric bundles of the same rank are isomorphic as quasi-Hilbert bundles precisely when they are smoothly isomorphic outside a closed set of zero capacity (Theorem~\ref{thm:smooth-bundle-quasi-classification}).

\section*{Prior work}
Cipriani and Sauvageot \cite{CiprianiSauvageot2003} represented Dirichlet forms by Hilbert-bimodule-valued derivations. For regular symmetric Dirichlet forms, Hinz, R\"ockner and Teplyaev \cite[Section~2]{HinzRocknerTeplyaev2013} represented the space of one-forms as a direct integral of Hilbert spaces over an energy-dominant measure. The bundle formulation in Theorem~\ref{thm:almost-cotangent-bundle} is equivalent to Gigli's Hilbert-module formulation \cite[Theorem~1.4.11 and Remark~1.4.12]{Gigli2018}. On the capacity side, Debin, Gigli and Pasqualetto \cite{DebinGigliPasqualetto2021} introduced a tangent module defined up to sets of zero capacity on RCD spaces.
\section{Almost Hilbert bundles}\label{sec:almost-bundles}
\subsection{Almost Hilbert bundles}
A scalar-valued function on $X$ is said to be almost measurable
if it coincides $m$-a.e.\ with a measurable function.
\begin{dfn}
Let $\{\mathcal H_x\}_{x\in X}$ be a family of Hilbert spaces.
A countable family of sections $(\xi_k)_{k\geq1}$, with
$\xi_k(x)\in \mathcal H_x$ for every $x\in X$, is called an
almost fundamental family if:
\begin{enumerate}
    \item For every $k,l\geq1$, the function
    $
        x\mapsto\langle\xi_k(x),\xi_l(x)\rangle
    $
    is almost measurable.

    \item For $m$-a.e.\ $x\in X$,
    $
        \mathcal H_x
        =
        \overline{\operatorname{Span}
        \{\xi_k(x)\mid k\geq1\}}.
    $
\end{enumerate}
\end{dfn}

\begin{dfn}
Let $(\xi_k)_{k\geq1}$ and $(\xi'_{k'})_{k'\geq1}$
be almost fundamental families of
$\{\mathcal H_x\}_{x\in X}$.
They are said to be equivalent if, for every $k,k'\geq1$,
the function
$
    x\mapsto\langle\xi_k(x),\xi'_{k'}(x)\rangle
$
is almost measurable.
\end{dfn}

\begin{rem}
The above relation on almost fundamental families
is an equivalence relation.
\end{rem}

\begin{proof}
Reflexivity and symmetry are immediate. To prove transitivity,
suppose that $(\xi_k)_{k\geq1}$ and $(\xi'_{k'})_{k'\geq1}$
are equivalent, and that $(\xi'_{k'})_{k'\geq1}$ and
$(\xi''_{k''})_{k''\geq1}$ are equivalent.

Apply the Gram--Schmidt process pointwise to
$(\xi'_{k'})_{k'\geq1}$, setting the resulting vector equal
to zero whenever the residual vector vanishes.
This gives sections $(e_j)_{j\geq1}$ such that each $e_j$
is a finite linear combination of $\xi'_1,\ldots,\xi'_j$
with almost measurable coefficients, and, for $m$-a.e.\
$x\in X$, the nonzero members of $(e_j(x))_{j\geq1}$
form an orthonormal basis of $\mathcal H_x$.

For every $k,k''\geq1$, we therefore have
\[
    \langle\xi_k(x),\xi''_{k''}(x)\rangle
    =
    \sum_{j=1}^{\infty}
    \langle\xi_k(x),e_j(x)\rangle
    \langle e_j(x),\xi''_{k''}(x)\rangle
    \quad\text{for }m\text{-a.e. }x\in X.
\]
Each summand is almost measurable by the assumed equivalences.
Hence the left-hand side is almost measurable, which proves
transitivity.
\end{proof}

\begin{dfn}
An almost Hilbert bundle $\mathcal H$ over $(X,m)$
is a family of Hilbert spaces $\{\mathcal H_x\}_{x\in X}$
together with an equivalence class of almost fundamental families.
\end{dfn}

\begin{dfn}
Let $\mathcal H$ be an almost Hilbert bundle over $(X,m)$.
A section $X\ni x\mapsto\xi(x)\in \mathcal H_x$ is said to be
almost measurable if, for some (equivalently, every)
representative $(\xi_k)_{k\geq1}$ of the equivalence class
of almost fundamental families associated with $\mathcal H$,
the function
$
    x\mapsto\langle\xi(x),\xi_k(x)\rangle
$
is almost measurable for every $k\geq1$.
\end{dfn}
\begin{dfn}
Let $\mathcal H^{(1)}$ and $\mathcal H^{(2)}$ be almost Hilbert
bundles over $(X,m)$.
Let $\Phi=(\Phi(x))_{x\in X}$ be a family of maps
$\Phi(x)\colon \mathcal H_x^{(1)}\to \mathcal H_x^{(2)}$.
\begin{itemize}
    \item $\Phi$ is called an almost Hilbert
    bundle morphism if $\Phi(x)$ is bounded linear
    for $m$-a.e.\ $x\in X$ and, for every almost measurable
    section $\xi$ of $\mathcal H^{(1)}$, the section
    $
        \Phi\xi\colon x\mapsto\Phi(x)\xi(x)
    $
    of $\mathcal H^{(2)}$ is almost measurable.

    \item $\Phi$ is called an almost Hilbert
    bundle isomorphism if it is an almost Hilbert bundle morphism and there exists an almost Hilbert bundle morphism
    $\Psi\colon \mathcal H^{(2)}\to \mathcal H^{(1)}$ such that
    \[
        \Psi(x)\Phi(x)=\operatorname{id}_{\mathcal H_x^{(1)}},
        \qquad
        \Phi(x)\Psi(x)=\operatorname{id}_{\mathcal H_x^{(2)}}
        \quad\text{for }m\text{-a.e. }x\in X.
    \]
\end{itemize}
An almost Hilbert bundle morphism $\Phi$
is said to be isometric if, for $m$-a.e.\ $x\in X$,
\[
    |\Phi(x)v|=|v|
    \qquad\text{for every }v\in \mathcal H_x^{(1)}.
\]
\end{dfn}
Throughout, we identify two almost Hilbert bundle morphisms
$\Phi,\Psi\colon \mathcal H^{(1)}\to \mathcal H^{(2)}$ whenever
\[
    \Phi(x)=\Psi(x)
    \quad\text{for }m\text{-a.e. }x\in X.
\]
We use the same notation for a morphism and its
equivalence class under this identification.
Composition is defined using representatives and is
independent of their choice.
\begin{dfn}
Two almost Hilbert bundles $\mathcal H^{(1)}$ and $\mathcal H^{(2)}$
are said to be isomorphic if there exists an almost Hilbert bundle isomorphism between them, and isometrically
isomorphic if such an isomorphism can be chosen to be isometric.
\end{dfn}
Any almost Hilbert bundle isomorphism $\Phi$ can be replaced
by the isometric almost Hilbert bundle isomorphism
\[
    \Phi_{\mathrm{iso}}
    :=
    \Phi|\Phi|^{-1}
    \qquad\text{for }m\text{-a.e.\ }x.
\]
Here $|\Phi|(x):=(\Phi(x)^\dagger\Phi(x))^{1/2}$ is the
fiberwise absolute value, and $\dagger$ denotes the Hilbert
space adjoint. Almost measurability is preserved by the
fiberwise functional calculus.

\subsection{Basic operations on almost Hilbert bundles}
\begin{dfn}
Let $\mathcal H$ be an almost Hilbert bundle over $(X,m)$.
An almost Hilbert subbundle of $\mathcal H$ is a family
$\mathcal K=(\mathcal K_x)_{x\in X}$ of closed linear subspaces $\mathcal K_x\subset \mathcal H_x$
for which there exists a countable family $(\eta_j)_{j\ge1}$
of almost measurable sections of $\mathcal H$ such that
\[
    \mathcal K_x
    =
    \overline{\operatorname{Span}\{\eta_j(x)\mid j\ge1\}}
    \qquad m\text{-a.e.}
\]
We equip $\mathcal K$ with the inherited fiberwise inner products
and the almost Hilbert bundle structure induced by $(\eta_j)$.
\end{dfn}
The natural inclusion $\mathcal K\hookrightarrow \mathcal H$ is an isometric
almost Hilbert bundle morphism.

\begin{dfn}
Let $\mathcal H$ be an almost Hilbert bundle over $(X,m)$,
and let $(\xi_k)_{k\geq1}$ be a representative of its
equivalence class of almost fundamental families.
For each $x\in X$, let $\mathcal H_x^*$ be the continuous dual of $\mathcal H_x$,
endowed with its canonical Hilbert space structure.
Define sections $(\xi_k^*)_{k\geq1}$ by
\[
    \<\xi_k^*(x),v\>=\langle v,\xi_k(x)\rangle,
    \qquad x\in X,\quad v\in \mathcal H_x.
\]
The dual bundle $\mathcal H^*$ is the almost Hilbert bundle
with fibers $\{\mathcal H_x^*\}_{x\in X}$ and the equivalence class
of almost fundamental families represented by
$(\xi_k^*)_{k\geq1}$.
This equivalence class is independent of the choice of
$(\xi_k)_{k\geq1}$.
\end{dfn}

\begin{dfn}
Let $\mathcal H^{(1)}$ and $\mathcal H^{(2)}$ be almost Hilbert
bundles over $(X,m)$.
Let $(\xi_k^{(i)})_{k\geq1}$ be a representative
of the equivalence class of almost fundamental
families associated with $\mathcal H^{(i)}$.
The Hilbert tensor product $\mathcal H^{(1)}\widehat{\otimes}\mathcal H^{(2)}$
is the almost Hilbert bundle with fibers
$\{\mathcal H_x^{(1)}\widehat{\otimes}\mathcal H_x^{(2)}\}_{x\in X}$
and the equivalence class of almost fundamental
families represented by
$(\xi_k^{(1)}\otimes\xi_l^{(2)})_{k,l\geq1}$, where
\[
    (\xi_k^{(1)}\otimes\xi_l^{(2)})(x)
    =
    \xi_k^{(1)}(x)\otimes\xi_l^{(2)}(x).
\]
This equivalence class is independent of the chosen
representatives.
\end{dfn}

\begin{dfn}
Let $\mathcal H$ be an almost Hilbert bundle over $(X,m)$,
and let $(\xi_j)_{j\geq1}$ be an almost fundamental family
representing its almost Hilbert bundle structure.
For every integer $k\geq1$, the exterior power
$\Lambda^k\mathcal H$ is the almost Hilbert bundle with
fibers $\Lambda^k\mathcal H_x$ and almost fundamental family
\[
    \bigl(\xi_{i_1}\wedge\cdots\wedge\xi_{i_k}\bigr)
    _{i_1,\ldots,i_k\geq1}.
\]
Here $\Lambda^k\mathcal H_x$ denotes the completion of
the algebraic exterior power with respect to the inner product
\[
    \left\langle
    u_1\wedge\cdots\wedge u_k,\,
    v_1\wedge\cdots\wedge v_k
    \right\rangle
    =
    \det\bigl(\langle u_i,v_j\rangle\bigr)_{i,j=1}^k,
\]
and exterior products of sections are defined pointwise:
\[
    (\alpha_1\wedge\cdots\wedge\alpha_k)(x)
    =
    \alpha_1(x)\wedge\cdots\wedge\alpha_k(x).
\]
\end{dfn}
The determinant identity shows that the displayed family
is almost fundamental.
Its equivalence class is independent of the chosen
representative, and exterior products of almost measurable
sections are almost measurable.
The pointwise norm and the pointwise inner product of
almost measurable sections are almost measurable.
\begin{dfn}
Let $\mathcal H$ be an almost Hilbert bundle over $(X,m)$.
We denote by $\L^0(m;\mathcal H)$ the space of almost measurable
sections of $\mathcal H$ modulo $m$-a.e.\ equality.
For $p\in[1,\infty]$ and $\xi\in\L^0(m;\mathcal H)$, set
\[
    \|\xi\|_p
    \coloneqq
    \begin{cases}
        \displaystyle
        \left(\int_X|\xi(x)|^p\,dm(x)\right)^{1/p},
        &1\leq p<\infty,\\[6pt]
        \displaystyle
        \operatorname*{ess\,sup}_{x\in X}|\xi(x)|,
        &p=\infty.
    \end{cases}
\]
The space of $\L^p$ sections of $\mathcal H$ is defined by
\[
    \L^p(m;\mathcal H)
    \coloneqq
    \{\xi\in\L^0(m;\mathcal H)\mid \|\xi\|_p<\infty\}.
\]
For $p=2$, we equip $\L^2(m;\mathcal H)$ with the inner product
\[
    \langle\xi,\eta\rangle_2
    \coloneqq
    \int_X\langle\xi(x),\eta(x)\rangle\,dm(x).
\]
\end{dfn}
\begin{prp}\label{prp:lp-sections-complete}
Let $\mathcal H$ be an almost Hilbert bundle over $(X,m)$.
For every $p\in[1,\infty]$, the space $\L^p(m;\mathcal H)$
with norm $\|\cdot\|_p$ is a Banach space.
In particular, $\L^2(m;\mathcal H)$ with the inner product
$\langle\cdot,\cdot\rangle_2$ is a Hilbert space.
\end{prp}
\begin{proof}
It suffices to show that every series
$\sum_{n=1}^\infty\alpha_n$ with
\[
    \alpha_n\in\L^p(m;\mathcal H),
    \qquad
    \sum_{n=1}^\infty\|\alpha_n\|_p<\infty
\]
converges in $\L^p(m;\mathcal H)$.
Choose almost measurable representatives of the $\alpha_n$.
To control their pointwise sum, set
\[
    g(x):=\sum_{n=1}^\infty|\alpha_n(x)|.
\]
For $p<\infty$, Minkowski's inequality for finite sums followed
by monotone convergence gives
\[
    \|g\|_{p}
    \leq\sum_{n=1}^\infty\|\alpha_n\|_p<\infty.
\]
For $p=\infty$, the same bound follows by taking a common null
set outside which
$|\alpha_n(x)|\leq\|\alpha_n\|_\infty$ for every $n$.
Thus $g$ is finite $m$-a.e., so completeness of each fiber
allows us to define
\[
    \alpha(x):=\sum_{n=1}^\infty\alpha_n(x)
    \quad\text{for }m\text{-a.e.\ }x.
\]
Set $\alpha=0$ on the remaining null set.
For every section $\eta$ in an almost fundamental family,
\[
    \langle\alpha(x),\eta(x)\rangle
    =
    \lim_{N\to\infty}
    \sum_{n=1}^N\langle\alpha_n(x),\eta(x)\rangle
    \quad\text{for }m\text{-a.e.\ }x.
\]
The right-hand side is almost measurable. Hence $\alpha$ is
an almost measurable section, and $|\alpha|\leq g$ shows that $\alpha\in\L^p(m;\mathcal H)$.
The preceding scalar norm
bound give
\[
    \left\|\alpha-\sum_{n=1}^N\alpha_n\right\|_p
    \leq
    \left\|\sum_{n>N}|\alpha_n|\right\|_{\L^p(m)}
    \leq
    \sum_{n>N}\|\alpha_n\|_p
    \longrightarrow0.
\]
Thus every absolutely convergent series converges, proving completeness.
\end{proof}

\subsection{Almost tangent bundles}
Assume that $\mathcal E$ admits a carr\'e du champ $\Gamma$.

\begin{thm}[The almost cotangent bundle]\label{thm:almost-cotangent-bundle}
There exist an almost Hilbert bundle $\mathcal H$ over $(X,m)$
and a linear map
\[
    d\colon\D(\mathcal E)\longrightarrow\L^2(m;\mathcal H)
\]
satisfying the following conditions:
\begin{enumerate}
    \item For every $f,g\in\D(\mathcal E)$,
    \[
        \langle df,dg\rangle=\Gamma(f,g)
        \qquad\text{in }\L^1(m).
    \]

    \item There exists a sequence $(f_n)_{n\geq1}\subset\D(\mathcal E)$
    such that representatives of $(df_n)_{n\geq1}$ form an
    almost fundamental family representing the almost Hilbert
    bundle structure of $\mathcal H$.
\end{enumerate}
Moreover, if $(\mathcal H^{(1)},d^{(1)})$ and
$(\mathcal H^{(2)},d^{(2)})$ satisfy these conditions,
then there exists a unique isometric almost Hilbert bundle
isomorphism
$
    \Phi\colon\mathcal H^{(1)}\to\mathcal H^{(2)}
$
such that the following diagram commutes:
\[
\begin{tikzcd}
    & \D(\mathcal E)
    \arrow[dl,"d^{(1)}"']
    \arrow[dr,"d^{(2)}"]
    &
    \\
    \L^2(m;\mathcal H^{(1)})
    \arrow[rr,"\Phi"']
    &&
    \L^2(m;\mathcal H^{(2)}).
\end{tikzcd}
\]
Here the bottom arrow denotes the induced map
$\xi\mapsto\Phi\xi$.
\end{thm}

\begin{proof}\;
\begin{itemize}
    \item We first construct the almost Hilbert bundle $\mathcal H$.
    Choose a sequence $(f_i)_{i\geq1}$ dense in $\D(\mathcal E)$
    for the form norm, and choose symmetric measurable representatives
    $G_{ij}$ of $\Gamma(f_i,f_j)$.
    For each $n$ and each $a\in\mathbb Q^n$, we have
    \[
        \sum_{i,j=1}^n a_i a_j G_{ij}
        =\Gamma\Bigl(\sum_{i=1}^n a_i f_i,\sum_{i=1}^n a_i f_i\Bigr)
        \geq0 \qquad m\text{-a.e.}
    \]
    Since there are only countably many such pairs $(n,a)$, we may
    remove a single null set $N$. By continuity in $a$, every finite
    matrix $(G_{ij}(x))_{i,j=1}^n$ is positive semidefinite for $x\notin N$.
    For $x\notin N$, equip the space of finitely supported real sequences
    with the positive semidefinite form
    \[
        (a,b)_x:=\sum_{i,j}a_i b_j G_{ij}(x).
    \]
    Let $\mathcal H_x$ be the completion of the quotient by the null
    space of this form, and denote the image of the $i$-th coordinate
    vector by $e_i(x)$. On $N$, set $\mathcal H_x=\{0\}$ and $e_i(x)=0$.
    The sections $e_i$ have measurable pairwise inner products, equal
    to $G_{ij}$ off $N$, and their span is dense in each fiber.
    Thus $(e_i)_{i\geq1}$ is an almost fundamental family and defines
    the almost Hilbert bundle $\mathcal H$.

    \item We next construct the differential $d$.
    For a finite sum $u=\sum_i a_i f_i$, set $d_0u:=\sum_i a_i e_i$.
    The Gram matrix identity gives
    \[
        \int_X\Bigl|\sum_i a_i e_i\Bigr|^2\,dm
        =\int_X\Gamma(u,u)\,dm.
    \]
    In particular, two representations of the same $u$ give the same
    element of $\L^2(m;\mathcal H)$. The integrated carr\'e du champ
    is bounded by a constant times $\mathcal E(u,u)$, so $d_0$ extends
    by continuity in the form norm to a linear map
    $d\colon\D(\mathcal E)\to\L^2(m;\mathcal H)$.
    On the span of the $f_i$, the same Gram matrix identity gives
    $\langle d_0u,d_0v\rangle=\Gamma(u,v)$ in $\L^1(m)$.
    Both sides are continuous in the form norm with values in $\L^1(m)$:
    for the left-hand side this follows from Cauchy--Schwarz in
    $\L^2(m;\mathcal H)$, and for the right-hand side from
    \[
        \|\Gamma(u,v)\|_{\L^1(m)}
        \leq
        \Bigl(\int_X\Gamma(u,u)\,dm\Bigr)^{1/2}
        \Bigl(\int_X\Gamma(v,v)\,dm\Bigr)^{1/2}.
    \]
    Density now gives $\langle df,dg\rangle=\Gamma(f,g)$ for all
    $f,g\in\D(\mathcal E)$. Since $df_i$ is represented by $e_i$,
    the second condition of the theorem also holds.

    \item We finally prove that the pair $(\mathcal H,d)$ is unique up to a unique isometric almost Hilbert bundle isomorphism intertwining the differentials.
    Suppose $(\mathcal H^{(1)},d^{(1)})$ and
    $(\mathcal H^{(2)},d^{(2)})$ both satisfy the two conditions.
    Form a single sequence $(h_i)_{i\geq1}$ from the two sequences in
    the second condition and a form-norm dense sequence in $\D(\mathcal E)$.
    Choose representatives of the sections $d^{(k)}h_i$ for $k=1,2$.
    Outside a single null set, these sections span a dense subspace
    in each fiber of their respective bundles and, for all $i,j$,
    \[
        \langle d^{(1)}h_i(x),d^{(1)}h_j(x)\rangle
        =\langle d^{(2)}h_i(x),d^{(2)}h_j(x)\rangle,
    \]
    since both inner products represent $\Gamma(h_i,h_j)$.
    Consequently the assignment $d^{(1)}h_i(x)\mapsto d^{(2)}h_i(x)$
    preserves the norm of every finite linear combination. It therefore
    extends uniquely to a surjective linear isometry
    $\Phi(x)\colon\mathcal H^{(1)}_x\to\mathcal H^{(2)}_x$.
    To check almost measurability, let $\xi$ be an almost measurable
    section of $\mathcal H^{(1)}$. For every $i$,
    \[
        \langle\Phi\xi,d^{(2)}h_i\rangle
        =\langle\xi,d^{(1)}h_i\rangle \qquad m\text{-a.e.}
    \]
    The right-hand side is almost measurable, and the sections on the
    left include an almost fundamental family for $\mathcal H^{(2)}$.
    Hence $\Phi\xi$ is almost measurable. The same argument applies to
    the inverse, so $\Phi$ is an isometric almost Hilbert bundle isomorphism.
    By construction, $\Phi d^{(1)}h_i=d^{(2)}h_i$ for every $i$.
    Both differentials are continuous in the form norm by the first
    condition, so density gives $\Phi d^{(1)}f=d^{(2)}f$ for every
    $f\in\D(\mathcal E)$. Any other such isomorphism agrees with $\Phi$
    on all the $d^{(1)}h_i$ outside a single null set, hence on their
    closed span in each fiber there. This proves uniqueness.
\end{itemize}
\end{proof}

\begin{dfn}
The almost Hilbert bundle $\mathcal H$ constructed in Theorem~\ref{thm:almost-cotangent-bundle}, equipped with $d$, is called the almost cotangent bundle
associated with $(X,m,\mathcal E,\Gamma)$ and is denoted by $T^*X$.
Its dual bundle
\[
    TX\coloneqq(T^*X)^*
\]
is called the almost tangent bundle.
\end{dfn}

\clearpage
\section{Quasi-Hilbert bundles}\label{sec:quasi-bundles}
\subsection{Quasi-Hilbert bundles}
\begin{dfn}
Let $\{\widetilde{\mathcal H}_x\}_{x\in X}$ be a family of Hilbert spaces.
A countable family $(\xi_k)_{k\geq1}$ of sections of $\widetilde{\mathcal H}$
is called a quasi-fundamental family if:
\begin{enumerate}
    \item For every $k,l\geq1$, the function
    $
        x\mapsto\langle\xi_k(x),\xi_l(x)\rangle
    $
    is quasi-continuous.

    \item For q.e.\ $x\in X$,
    $
        \widetilde{\mathcal H}_x
        =
        \overline{\operatorname{Span}
        \{\xi_k(x)\mid k\geq1\}}.
    $
\end{enumerate}
\end{dfn}
In particular, every quasi-fundamental family is an almost fundamental family, since quasi-continuous functions are almost measurable and every exceptional set is $m$-null.
\begin{dfn}
Let $\{\widetilde{\mathcal H}_x\}_{x\in X}$ be a family of Hilbert spaces.
Let $(\xi_k)_{k\geq1}$ and $(\xi'_{k'})_{k'\geq1}$
be quasi-fundamental families of $\widetilde{\mathcal H}$.
They are said to be equivalent if, for every $k,k'\geq1$,
the function
$
    x\mapsto\langle\xi_k(x),\xi'_{k'}(x)\rangle
$
is quasi-continuous.
\end{dfn}

\begin{rem}
The above relation on quasi-fundamental families
is an equivalence relation.
\end{rem}

\begin{proof}
Reflexivity and symmetry are immediate.
To prove transitivity, suppose that
$(\xi_k)_{k\geq1}$ and $(\xi'_{k'})_{k'\geq1}$ are equivalent,
and that $(\xi'_{k'})_{k'\geq1}$ and
$(\xi''_{k''})_{k''\geq1}$ are equivalent.
Choose a common nest $(F_n)_{n\geq1}$ such that all the
inner-product functions appearing in these assumptions
and in the definitions of the three quasi-fundamental
families are continuous on each $F_n$, and all three
families have dense linear span in $\widetilde{\mathcal H}_x$ for every
$x\in\bigcup_{n\geq1}F_n$.
Fix $n,k,k''\geq1$, $x\in F_n$, and $\epsilon>0$.
Choose $N\geq1$ and $a_1,\ldots,a_N\in\mathbb R$ such that
\[
    \left|\xi_k(x)-\sum_{j=1}^{N}a_j\xi'_j(x)\right|
    <\epsilon.
\]
For $y\in F_n$, set
\[
    f(y)\coloneqq\langle\xi_k(y),\xi''_{k''}(y)\rangle,
    \qquad
    g(y)\coloneqq
    \sum_{j=1}^{N}a_j
    \langle\xi'_j(y),\xi''_{k''}(y)\rangle.
\]
The functions $g$,
$y\mapsto|\xi''_{k''}(y)|$, and
$
    y\mapsto
    \left|\xi_k(y)-\sum_{j=1}^{N}a_j\xi'_j(y)\right|
$
are continuous on $F_n$.
By the Cauchy--Schwarz inequality,
\[
    |f(y)-g(y)|
    \leq
    \left|\xi_k(y)-\sum_{j=1}^{N}a_j\xi'_j(y)\right|
    |\xi''_{k''}(y)|.
\]
Consequently,
\[
    \varlimsup_{\substack{y\to x\\y\in F_n}}
    |f(y)-f(x)|
    \leq 2\epsilon|\xi''_{k''}(x)|.
\]
Since $\epsilon>0$ is arbitrary, $f$ is continuous on $F_n$.
Thus $f$ is quasi-continuous, proving transitivity.
\end{proof}
\begin{dfn}
A quasi-Hilbert bundle $\widetilde{\mathcal H}$ is a family of Hilbert spaces
$\{\widetilde{\mathcal H}_x\}_{x\in X}$ together with an equivalence class of
quasi-fundamental families.
This equivalence class is called the quasi-structure of $\widetilde{\mathcal H}$.
\end{dfn}
A quasi-structure on $\{\widetilde{\mathcal H}_x\}_{x\in X}$ canonically induces
an almost Hilbert bundle structure on the same family,
since equivalent quasi-fundamental families are equivalent
as almost fundamental families.
\begin{dfn}
Let $\widetilde{\mathcal H}$ be a quasi-Hilbert bundle.
A section $X\ni x\mapsto\xi(x)\in \widetilde{\mathcal H}_x$ is said to be
quasi-continuous if the functions
\[
    x\mapsto|\xi(x)|,
    \qquad
    x\mapsto\langle\xi(x),\xi_k(x)\rangle,
    \quad k\geq1,
\]
are quasi-continuous, where $(\xi_k)_{k\geq1}$
is a representative of the quasi-structure of $\widetilde{\mathcal H}$.
We denote the set of quasi-continuous sections by
$\Gamma^{\text{qc}}(\widetilde{\mathcal H})$.
\end{dfn}

\begin{rem}
The definition of a quasi-continuous section is independent
of the choice of the representative of the quasi-structure.
Moreover, the pointwise inner product of two quasi-continuous
sections is quasi-continuous.
\end{rem}

\begin{proof}
Both assertions follow from the same local approximation
argument used in the preceding proof, using the
quasi-continuity of the norms of the sections.
\end{proof}
Quasi-continuous sections are closed under addition and
multiplication by quasi-continuous scalar-valued functions.
Every member of a quasi-fundamental family representing
the quasi-structure of $\widetilde{\mathcal H}$ is quasi-continuous.
\begin{dfn}
Let $\widetilde{\mathcal H}^{(1)}$ and $\widetilde{\mathcal H}^{(2)}$ be quasi-Hilbert bundles.
Let $\Phi=(\Phi(x))_{x\in X}$ be a family of maps
$\Phi(x)\colon \widetilde{\mathcal H}_x^{(1)}\to \widetilde{\mathcal H}_x^{(2)}$.
\begin{itemize}
    \item $\Phi$ is called a quasi-Hilbert bundle morphism
    if $\Phi(x)$ is bounded linear for q.e.\ $x\in X$
    and, for every $\xi\in\Gamma^{\text{qc}}(\widetilde{\mathcal H}^{(1)})$, the section
    \[
        \Phi\xi\colon x\mapsto\Phi(x)\xi(x)
    \]
    belongs to $\Gamma^{\text{qc}}(\widetilde{\mathcal H}^{(2)})$.

    \item $\Phi$ is called a quasi-Hilbert bundle
    isomorphism if it is a quasi-Hilbert bundle morphism
    and there exists a quasi-Hilbert bundle morphism
    $\Psi\colon \widetilde{\mathcal H}^{(2)}\to \widetilde{\mathcal H}^{(1)}$ such that
    \[
        \Psi(x)\Phi(x)=\operatorname{id}_{\widetilde{\mathcal H}_x^{(1)}},
        \qquad
        \Phi(x)\Psi(x)=\operatorname{id}_{\widetilde{\mathcal H}_x^{(2)}}
        \quad\text{for q.e. }x\in X.
    \]
\end{itemize}
A quasi-Hilbert bundle morphism $\Phi$ is said to be
isometric if, for q.e.\ $x\in X$,
\[
    |\Phi(x)v|=|v|
    \qquad\text{for every }v\in \widetilde{\mathcal H}_x^{(1)}.
\]
\end{dfn}

We identify two quasi-Hilbert bundle morphisms
$\Phi,\Psi\colon \widetilde{\mathcal H}^{(1)}\to \widetilde{\mathcal H}^{(2)}$ whenever
\[
    \Phi(x)=\Psi(x)
    \quad\text{for q.e. }x\in X.
\]
We use the same notation for a morphism and its
equivalence class under this identification.
Composition is defined using representatives and is
independent of their choice.
Every quasi-Hilbert bundle morphism induces an almost
Hilbert bundle morphism between the underlying almost
Hilbert bundles. This follows from the same pointwise
approximation argument as in the almost setting, since
the images of the quasi-fundamental sections are almost
measurable.

\begin{dfn}
Two quasi-Hilbert bundles $\widetilde{\mathcal H}^{(1)}$ and $\widetilde{\mathcal H}^{(2)}$
are said to be quasi-isomorphic if there exists a
quasi-Hilbert bundle isomorphism between them, and
quasi-isometrically isomorphic if such an isomorphism
can be chosen to be isometric.
\end{dfn}
If a quasi-Hilbert bundle isomorphism $\Phi$ has a fiberwise
adjoint $\Phi^\dagger$ that is also a quasi-Hilbert bundle
morphism, then
\[
    \Phi_{\mathrm{iso}}(x)
    :=
    \Phi(x)|\Phi|(x)^{-1}
    \qquad\text{for q.e.\ }x
\]
is an isometric quasi-Hilbert bundle isomorphism.
This condition is not automatic
for general quasi-Hilbert bundle isomorphisms under the
present definition.

\subsection{Basic operations on quasi-Hilbert bundles}
\begin{dfn}
Let $\widetilde{\mathcal H}$ be a quasi-Hilbert bundle over $X$.
A quasi-Hilbert subbundle of $\widetilde{\mathcal H}$ is a family
$\widetilde{\mathcal K}=(\widetilde{\mathcal K}_x)_{x\in X}$ of closed linear subspaces $\widetilde{\mathcal K}_x\subset \widetilde{\mathcal H}_x$
for which there exists a countable family $(\eta_j)_{j\ge1}$
of quasi-continuous sections of $\widetilde{\mathcal H}$ such that
\[
    \widetilde{\mathcal K}_x
    =
    \overline{\operatorname{Span}\{\eta_j(x)\mid j\ge1\}}
    \qquad\text{q.e.}
\]
We equip $\widetilde{\mathcal K}$ with the inherited fiberwise inner products
and the quasi-Hilbert bundle structure induced by $(\eta_j)$.
\end{dfn}
The natural inclusion $\widetilde{\mathcal K}\hookrightarrow \widetilde{\mathcal H}$ is an isometric
quasi-Hilbert bundle morphism.
\begin{dfn}
Let $\widetilde{\mathcal H}$ be a quasi-Hilbert bundle, and let
$(\xi_k)_{k\geq1}$ be a representative of its
quasi-structure.
For each $k\geq1$, define a section $\xi_k^*$ of
the underlying almost dual bundle $\widetilde{\mathcal H}^*$ by
\[
    \<\xi_k^*(x),v\>=\langle v,\xi_k(x)\rangle,
    \qquad x\in X,\quad v\in \widetilde{\mathcal H}_x.
\]
The dual quasi-Hilbert bundle $\widetilde{\mathcal H}^*$ is the almost dual
bundle equipped with the quasi-structure represented
by $(\xi_k^*)_{k\geq1}$.
\end{dfn}

This quasi-structure is independent of the chosen
representative, by the same inner-product identities
as in the almost setting.
The fiberwise Riesz maps
\[
    \widetilde{\mathcal H}_x\ni v\longmapsto\langle\cdot,v\rangle\in \widetilde{\mathcal H}_x^*
\]
define an isometric quasi-Hilbert bundle isomorphism
from $\widetilde{\mathcal H}$ to $\widetilde{\mathcal H}^*$.
In particular, the evaluation of a quasi-continuous
dual section on a quasi-continuous section is
quasi-continuous.

\begin{dfn}
Let $\widetilde{\mathcal H}^{(1)}$ and $\widetilde{\mathcal H}^{(2)}$ be quasi-Hilbert bundles.
For $i=1,2$, let $(\xi_k^{(i)})_{k\geq1}$
be a representative of the quasi-structure of $\widetilde{\mathcal H}^{(i)}$.
The Hilbert tensor product
$\widetilde{\mathcal H}^{(1)}\widehat{\otimes}\widetilde{\mathcal H}^{(2)}$
is the almost Hilbert tensor product equipped with
the quasi-structure represented by
$(\xi_k^{(1)}\otimes\xi_l^{(2)})_{k,l\geq1}$, where
\[
    (\xi_k^{(1)}\otimes\xi_l^{(2)})(x)
    =
    \xi_k^{(1)}(x)\otimes\xi_l^{(2)}(x).
\]
\end{dfn}

The tensor-product family is quasi-fundamental, and its
equivalence class is independent of the chosen
representatives.
The proof is the same as in the almost setting, using
the identity
\[
    \langle u_1\otimes u_2,v_1\otimes v_2\rangle
    =
    \langle u_1,v_1\rangle\langle u_2,v_2\rangle
\]
and the fact that products of quasi-continuous scalar
functions are quasi-continuous.
In particular, if $\xi^{(i)}\in\Gamma^{\text{qc}}(\widetilde{\mathcal H}^{(i)})$,
$i=1,2$, then
\[
    \xi^{(1)}\otimes\xi^{(2)}
    \in
    \Gamma^{\text{qc}}(\widetilde{\mathcal H}^{(1)}\widehat{\otimes}\widetilde{\mathcal H}^{(2)}),
\]
since its pairings with the tensor-product fundamental
family are quasi-continuous and
\[
    |\xi^{(1)}(x)\otimes\xi^{(2)}(x)|
    =
    |\xi^{(1)}(x)|\,|\xi^{(2)}(x)|.
\]
\begin{dfn}
Let $\widetilde{\mathcal H}$ be a quasi-Hilbert bundle,
and let $(\xi_j)_{j\geq1}$ represent its quasi-structure.
For every integer $k\geq1$, the exterior power
$\Lambda^k\widetilde{\mathcal H}$ is the almost exterior
power equipped with the quasi-structure represented by
\[
    \bigl(\xi_{i_1}\wedge\cdots\wedge\xi_{i_k}\bigr)
    _{i_1,\ldots,i_k\geq1}.
\]
\end{dfn}
The determinant formula for inner products shows that this
family is quasi-fundamental and that its equivalence class
is independent of the chosen representative.
Exterior products are taken pointwise, and the exterior
product of quasi-continuous sections is quasi-continuous.

\subsection{Quasi-topology on quasi-Hilbert bundles}
Let $\widetilde{\mathcal H}$ be a quasi-Hilbert bundle over $(X,m,\mathcal{E})$.
\begin{dfn}
Let $(\xi_k)_{k\geq1}$ be a quasi-fundamental family of $\widetilde{\mathcal H}$.
A nest $(F_n)_{n\geq1}$ in $X$ is called adapted to
$(\xi_k)_{k\geq1}$ if, for every $n\geq1$:
\begin{enumerate}
    \item
    $\langle\xi_k,\xi_l\rangle|_{F_n}$ is continuous
    for every $k,l\geq1$.

    \item
    For every $x\in F_n,\;\widetilde{\mathcal H}_x=
        \overline{\operatorname{Span}
        \{\xi_k(x)\mid k\geq1\}}.$
\end{enumerate}
\end{dfn}
\begin{dfn}
Let $(\xi_k)_{k\geq1}$ be a quasi-fundamental family of $\widetilde{\mathcal H}$,
and let $(F_n)_{n\geq1}$ be a nest adapted to $(\xi_k)_{k\geq1}$.
For every $n\geq1$, we equip $\widetilde{\mathcal H}|_{F_n}$ with the weakest
topology for which
\begin{itemize}
    \item
    $\widetilde{\mathcal H}|_{F_n}\ni(x,a)\mapsto x\in F_n$,
    \item
    $\widetilde{\mathcal H}|_{F_n}\ni(x,a)\mapsto
    a_k:=\langle a,\xi_k(x)\rangle\in\mathbb R$,
    $k\geq1$,
    \item
    $\widetilde{\mathcal H}|_{F_n}\ni(x,a)\mapsto |a|\in\mathbb R$
\end{itemize}
are continuous.
\end{dfn}
\begin{rem}
Even if a nest $(F_n)_{n\geq1}$ is adapted to two equivalent
quasi-fundamental families $(\xi_k)_{k\geq1}$ and
$(\xi'_k)_{k\geq1}$, the induced topologies on $\widetilde{\mathcal H}|_{F_n}$
need not coincide.
Thus, in general, these topologies depend on both the nest
and the choice of the quasi-fundamental family.
\end{rem}
\begin{prp}\label{prp:compatible-nest-topologies}
Let $(\xi_k)_{k\geq1}$ and $(\xi'_k)_{k\geq1}$ be two equivalent
quasi-fundamental families of $\widetilde{\mathcal H}$, and let $(F_n)_{n\geq1}$
and $(F'_n)_{n\geq1}$ be nests adapted to these families,
respectively.
Then there exists a nest $(G_n)_{n\geq1}$ such that
\[
    G_n\subset F_n\cap F'_n
    \qquad\text{for every }n\geq1,
\]
and, for every $n\geq1$, the topologies on $\widetilde{\mathcal H}|_{G_n}$
induced by $(\xi_k)_{k\geq1}$ and $(\xi'_k)_{k\geq1}$ coincide.
\end{prp}
\begin{proof}
By the equivalence of the two quasi-fundamental families,
all the functions $\langle\xi_k,\xi'_l\rangle$ are
quasi-continuous.
Hence we may choose a nest $(G_n)_{n\geq1}$ such that
$G_n\subset F_n\cap F'_n$ and
$\langle\xi_k,\xi'_l\rangle|_{G_n}$ is continuous
for every $n,k,l\geq1$.
Fix $n\geq1$, and suppose that
$(x_\lambda,a_\lambda)\to(x,a)$ in $\widetilde{\mathcal H}|_{G_n}$
with respect to the topology induced by $(\xi_k)_{k\geq1}$.
Then
\[
    x_\lambda\to x,\qquad
    |a_\lambda|\to|a|,\qquad
    \langle a_\lambda,\xi_k(x_\lambda)\rangle
    \to\langle a,\xi_k(x)\rangle
    \quad\text{for every }k\geq1.
\]
Fix $l\geq1$.
For any finite linear combination
$\eta\in\operatorname{Span}\{\xi_k\mid k\geq1\}$,
the function $|\xi'_l-\eta|$ is continuous on $G_n$,
and
\[
    \langle a_\lambda,\eta(x_\lambda)\rangle
    \to\langle a,\eta(x)\rangle.
\]
The triangle inequality therefore gives
\[
    \varlimsup_\lambda
    \left|
        \langle a_\lambda,\xi'_l(x_\lambda)\rangle
        -\langle a,\xi'_l(x)\rangle
    \right|
    \leq
    2|a|\,|\xi'_l(x)-\eta(x)|.
\]
By the generating property at $x$, the right-hand side
can be made arbitrarily small.
Thus
\[
    \langle a_\lambda,\xi'_l(x_\lambda)\rangle
    \to\langle a,\xi'_l(x)\rangle.
\]
Consequently, the identity map from the topology induced by
$(\xi_k)_{k\geq1}$ to that induced by $(\xi'_k)_{k\geq1}$
is continuous.
Interchanging the two families proves the converse,
so the two topologies coincide.
\end{proof}
\clearpage
\section{On Riemannian manifolds}\label{sec:riemannian-bundles}
\subsection{Dirichlet energy with Neumann boundary conditions}
Let $M$ be a smooth Riemannian manifold with smooth boundary
$\partial M$.
We denote by $\mathrm{vol}$ the Riemannian volume measure on $M$,
by $\nu$ the outward unit normal vector field along $\partial M$,
and by $\sigma$ the induced Riemannian surface measure on $\partial M$.
Smoothness is understood up to the boundary.
We regard the differential $d$ as a densely defined operator
from $\L^2(\mathrm{vol})$ to $\L^2(\mathrm{vol};T^*M)$ with domain
\[
\D(d)
:=
\left\{
f\in C^\infty(M)
\,\middle|\,
f\in\L^2(\mathrm{vol}),\
df\in\L^2(\mathrm{vol};T^*M)
\right\}.
\]
Let $d^\dagger$ denote its formal adjoint with respect to
$\mathrm{vol}$ and $g$.
The integration-by-parts formula reads
\[
\int_M\langle df,\alpha\rangle\,d\mathrm{vol}
=
\int_M f\,d^\dagger\alpha\,d\mathrm{vol}
+
\int_{\partial M}f\,\alpha(\nu)\,d\sigma,
\]
for $f\in C^\infty(M)$ and $\alpha\in\Gamma_c(T^*M)$.
By testing this identity with $\alpha$ supported in the interior
of $M$, we see that $d$ is closable.
We denote its closure again by $d$ and set
\[
W^{1,2}:=\D(d),
\qquad
\|f\|_{W^{1,2}}^2
:=\|f\|_2^2+\|df\|_2^2.
\]
We define
\[
\mathcal E(f,g)
:=
\frac12\int_M\langle df,dg\rangle\,d\mathrm{vol}
\quad,\quad
f,g\in\D(\mathcal E):=W^{1,2}.
\]
Then $\mathcal E$ is a densely defined, symmetric, nonnegative
closed form on $\L^2(\mathrm{vol})$.
The chain rule implies the Markov property, so $\mathcal E$
is a Dirichlet form.
We call it the Dirichlet form with Neumann boundary conditions
associated with $M$.
\subsection{Quasi-structure of smooth metric vector bundles}
\begin{dfn}\label{dfn:natural-quasi-structure}
Let $E\to M$ be a smooth metric vector bundle.
Choose a countable family $(\xi_k)_{k\geq1}$ of smooth sections
with compact support in $M$ such that
\[
    E_x=\operatorname{Span}\{\xi_k(x)\mid k\geq1\}
    \qquad\text{for every }x\in M.
\]
Since each pairing $\langle\xi_k,\xi_l\rangle$ is smooth,
this family is quasi-fundamental with respect to $\mathcal E$.
We call the quasi-structure represented by this family
the natural quasi-structure on $E$.
It is independent of the choice of the smooth generating family,
since the pairings between any two such families are smooth.
\end{dfn}
\begin{rem}
    All sections are understood to be smooth up to the boundary,
and their compact supports may meet $\partial M$.
The fibers over the boundary are the full fibers of $E$.
In particular, for the tangent bundle we retain
\[
    T_xM=T_x\partial M\oplus\mathbb R\nu_x
    \qquad\text{for every }x\in\partial M.
\]
\end{rem}
\begin{prp}\label{prp:almost-triviality}
Let $E\to M$ be a smooth metric vector bundle of rank $r$.
Then
\[
    E\cong\underline{\mathbb R^r}
\]
as almost Hilbert bundles.
\end{prp}
\begin{proof}
Choose a countable open cover $(U_n)_{n\geq1}$ of $M$
such that $E$ admits a smooth orthonormal frame on each $U_n$.
The sets
$
    A_n:=U_n\setminus\bigcup_{j<n}U_j
$
form a Borel partition of $M$.
Gluing the local orthonormal frames along this partition yields
a global Borel measurable orthonormal frame of $E$,
which gives the required isometric almost isomorphism.
\end{proof}
\begin{ex}
    By Proposition~\ref{prp:almost-triviality}, the trivial and M\"obius line bundles over $S^1$ are isomorphic as almost Hilbert bundles.
\end{ex}

\begin{thm}[Quasi-isomorphism of smooth metric bundles]\label{thm:smooth-bundle-quasi-classification}
Let $E^{(1)}$ and $E^{(2)}$ be smooth real metric vector bundles
of the same rank $r$ over $M$, equipped with the quasi-structures of Definition~\ref{dfn:natural-quasi-structure}.
Then the following conditions are equivalent:
\begin{enumerate}
    \item
    $E^{(1)}\cong E^{(2)}$ as quasi-Hilbert bundles.

    \item
    $\det E^{(1)}\cong\det E^{(2)}$
    as smooth line bundles over $M$.

    \item
    There exists an exceptional closed set $N\subset M$ such that
    \[
        E^{(1)}|_{M\setminus N}
        \cong
        E^{(2)}|_{M\setminus N}
    \]
    as smooth vector bundles over $M\setminus N$.
\end{enumerate}
Here $\det E^{(i)}:=\bigwedge^r E^{(i)}$.
All the above isomorphisms may be chosen to be isometric.
\end{thm}
\begin{proof}
    For $1\Rightarrow 2$, the induced isometry of determinant lines gives a quasi-trivialization of $(\det E^{(1)})^*\otimes\det E^{(2)}$. Lemma~\ref{lem:quasi-trivial-line} makes this line bundle smoothly trivial, so the determinant lines are smoothly isomorphic.
\begin{lem}\label{lem:quasi-continuous-signs}
Let $U\subset M$ be a connected open set, and let
$f\colon U\to\mathbb R$ be quasi-continuous.
Suppose that
\[
    f(x)\in\{-1,1\}
    \qquad\text{q.e. on }U.
\]
Then there exists $c\in\{-1,1\}$ such that
\[
    f(x)=c
    \qquad\text{q.e. on }U.
\]
\end{lem}
\begin{proof}
Let $V\Subset U$ be an arbitrary connected open set.
By the definition of quasi-continuity and the assumption on $f$,
there exists a decreasing sequence $(G_n)_{n\ge1}$ of open subsets
of $M$ such that
\[
    \operatorname{Cap}_1(G_n)\longrightarrow0
    \qquad\text{as }n\to\infty,
\]
and, for every $n\ge1$, the restriction
$f|_{\overline V\setminus G_n}$ is continuous and
\[
    f(x)\in\{-1,1\}
    \qquad\text{for every }x\in\overline V\setminus G_n.
\] 
Let $e_n$ be the $1$-equilibrium potential of $G_n$,
identified with its quasi-continuous representative.
By definition,
\[
    0\le e_n\le1
    \quad\text{q.e. on }M,
    \qquad
    e_n=1
    \quad\text{q.e. on }G_n,
\]
and
\[
    \|e_n\|_{W^{1,2}}\longrightarrow0
    \qquad\text{as }n\to\infty.
\]
Since
$\{f|_{\overline V\setminus G_n}=-1\}$ and
$\{f|_{\overline V\setminus G_n}=1\}$
are disjoint compact subsets of $M$, there exists
$h_n\in C^\infty(M)$ with $|h_n|\le1$ which equals
$-1$ and $1$ on neighborhoods of these sets, respectively. Define 
\[f_n:=h_n(1-e_n):V\longrightarrow\mathbb{R}.\]
Since $h_n=f$ on $V\setminus G_n$, the uniform
bound $|h_n|\le1$ and the defining properties of $e_n$
imply that
\[
    f_n\longrightarrow f
    \qquad\text{in }\L^2(V).
\]
Since $dh_n=0$ on $V\setminus G_n$ and $e_n=1$ q.e. on $G_n$,
we have
\[
    df_n=(1-e_n)\,dh_n-h_n\,de_n=-h_n\,de_n
    \qquad\text{a.e. on }V.
\]
Thus, the bound $|h_n|\le1$ and the defining properties of $e_n$
yield
\[
    df_n\longrightarrow0
    \qquad\text{in }\L^2(V;T^*V).
\]
Thus, by the closability of $d$, we have $f\in W^{1,2}(V)$
and $df=0$ a.e. on $V$.
Since $V\Subset U$ was arbitrary, a countable cover of $U$
by such sets yields $f\in W_{\mathrm{loc}}^{1,2}(U)$ and
\[
    df=0
    \qquad\text{a.e. on }U.
\]
Since $U$ is connected, $f$ is constant a.e. on $U$.
As $f$ is quasi-continuous, $f$ is constant q.e. on $U$.
This proves the claim.
\end{proof}
\begin{lem}\label{lem:quasi-trivial-line}
Let $L$ be a smooth real metric line bundle over $M$.
If $L$ is quasi-isometric to the trivial line bundle $\underline {\mathbb R}$, then
$L$ is smoothly  isomorphic
to the trivial line bundle $\underline {\mathbb {R}}$.
\end{lem}
\begin{proof}
    Let $\Phi:\underline{\mathbb R}\to L$ be a quasi-isometric isomorphism. Choose a countable cover $(U_i)_{i\ge1}$ of $M$ by connected trivializing neighborhoods of $L$, and, for each $i\ge1$, a smooth unit section $e_i$ of $L|_{U_i}$. By Lemma~\ref{lem:quasi-continuous-signs}, for each $i\ge1$, there exists
a constant $c_i\in\{-1,1\}$ such that
\[
    \Phi(1)|_{U_i}=c_i e_i
    \qquad\text{q.e. on }U_i.
\]
Set $s_i:=c_i e_i$. On each overlap $U_i\cap U_j$, the smooth
sections $s_i$ and $s_j$ agree outside a set of zero capacity.
Since every nonempty open set has positive capacity, continuity
implies that $s_i=s_j$ everywhere on $U_i\cap U_j$.
Hence the sections $(s_i)_{i\ge1}$ glue to a smooth global
unit section of $L$. This proves the claim.
\end{proof}
For $2\Rightarrow 3$, the rank-one case follows from $E=\det E$. In higher rank, apply Lemma~\ref{lem:determinant-splitting} to both bundles. Outside the union of the two closed exceptional sets, the splittings and the isomorphism of determinant lines give the required smooth isomorphism.
\begin{lem}\label{lem:transverse-section}
    Let $E$ be a smooth vector bundle of rank $r\ge2$ over $M$. Then there exist a closed smooth submanifold $N\subset M$
of codimension $r$, possibly empty, and a smooth section $\xi$
of $E$ such that $\xi$ is nowhere vanishing on $M\setminus N$.
\end{lem}
\begin{proof}
    For sufficiently large $L\in\mathbb N$, we may identify $E$
isometrically with a smooth subbundle of $M\times \mathbb{R}^L$. Let $P(x)\colon\mathbb{R}^L\to E_x$
denote the fiberwise orthogonal projection. Define a smooth map $F:M\times \mathbb{R}^L\to E$ by
\[
    F(x,a):=P(x)a.
\]
Identifying $M$ with the image of the zero section in $E$,
we see that $F$ is transverse to $M$. Thus, by the parametric transversality theorem on manifolds with boundary, for almost every
$a\in\mathbb R^L$, the section $F_a\coloneqq F(\cdot,a)$ is transverse to the zero section $M\subset E$. Fix such an $a\in\mathbb R^L$. Then the set $N\coloneqq F_a^{-1}(M)$
is either empty or a closed smooth submanifold
of $M$ of codimension $r$. Thus, $\xi\coloneqq F_a$ and the set $N$ have the required properties.
\end{proof}

\begin{lem}\label{lem:determinant-splitting}
Let $E$ be a smooth real metric vector bundle of rank $r\ge 2$
over $M$. Then there exists a closed exceptional set $N\subset M$
such that
\[
    E|_{M\setminus N}
    \cong
    \bigl(\underline{\mathbb R}^{\,r-1}\oplus\det E\bigr)|_{M\setminus N}
\]
as smooth metric vector bundles.
\end{lem}
\begin{proof}
By Lemma~\ref{lem:transverse-section}, there exist a closed smooth submanifold
$N\subset M$ of codimension $r$, possibly empty, and a smooth
section $\xi$ of $E$ which is nowhere vanishing on $M\setminus N$.
Since $r\ge2$, the set $N$ is exceptional.
The smooth unit section $\xi/|\xi|$ spans a trivial line
subbundle of $E|_{M\setminus N}$. Let $E'$ denote its
orthogonal complement. Then we have a smooth isometric isomorphism
\[
    E|_{M\setminus N}
    \cong \underline{\mathbb R}|_{M\setminus N}\oplus E'.
\]
We repeat this construction while the remaining bundle has
rank at least two. At each step, the new zero set is closed
in the current open manifold and has codimension at least two.
Its union with the previously removed set is therefore closed
in $M$ and exceptional.
Thus, after enlarging $N$, we obtain a smooth metric line
bundle $L$ over $M\setminus N$ and a smooth isometric isomorphism
\[
    E|_{M\setminus N}
    \cong \underline{\mathbb R}^{r-1}|_{M\setminus N}\oplus L.
\]
Taking determinants yields a smooth isometric isomorphism
\[
    L\cong(\det E)|_{M\setminus N},
\]
which proves the claim.
\end{proof}
For $3\Rightarrow 1$, a smooth isometry on $M\setminus N$ preserves the natural quasi-structures of Definition~\ref{dfn:natural-quasi-structure}, since $N$ is exceptional.
\end{proof}
Since $w_1(\det E)=w_1(E)$ \cite{MilnorStasheff1974}, the classical classification of real line bundles \cite[Proposition~3.10]{Hatcher2017} identifies condition~2 of Theorem~\ref{thm:smooth-bundle-quasi-classification} with the equality $w_1(E^{(1)})=w_1(E^{(2)})$.
\begin{cor}\label{cor:orientable-quasi-triviality}
    Let $E$ be a smooth real metric vector bundle of rank $r$
over $M$. If $E$ is orientable, then
$
    E\cong\underline{\mathbb R}^{r}
$
isometrically as quasi-Hilbert bundles.
\end{cor}
\begin{proof}
    Since $E$ is orientable, its determinant line is trivial. Theorem~\ref{thm:smooth-bundle-quasi-classification} applied to $E$ and $\underline{\mathbb R}^{r}$ gives the conclusion.
\end{proof}
\begin{ex}\;
    \begin{enumerate}
    \item
Let $L\to S^1$ be the Möbius line bundle.
Then the trivial bundle $\underline{\mathbb R}\to S^1$ and $L\to S^1$ are not isomorphic as quasi-Hilbert bundles, since $L$ is not smoothly trivial.
    \item
    The tangent bundle $TS^2\to S^2$ is quasi-isometrically trivial,
    by Corollary~\ref{cor:orientable-quasi-triviality}, since $S^2$ is orientable.

    \item
    Let $L\to S^1$ be the Möbius line bundle.
    The product bundle
    \[
        L\times L\longrightarrow T^2=S^1\times S^1
    \]
    has rank two and is not quasi-isometrically trivial.
    Indeed, its determinant restricts to the Möbius line bundle
    on $S^1\times\{p\}$ for any $p\in S^1$, and is therefore
    nontrivial.
\end{enumerate}
\end{ex}

\clearpage
\chapter{Kato Forms and Quasi-realization}\label{chap:bundle-forms-quasi-realization}
\section*{Overview}
We axiomatize the connection energy of one-forms on a Riemannian manifold by requiring that the heat semigroup of a nonnegative bundle form be dominated by the scalar heat semigroup. We call a bundle form satisfying this condition a Kato form (Definition~\ref{dfn:kato-form}). Theorem~\ref{thm:kato-domination} gives equivalent formulations in terms of the form, its generator and its resolvent. Our main result (Theorems~\ref{thm:quasi-realization-existence} and~\ref{thm:quasi-realization-uniqueness}) is that every Kato form determines a quasi-Hilbert bundle, unique up to isomorphism, on which each finite-energy section has a quasi-continuous representative. On this quasi-realization we define the resolvent of a vector measure of finite energy integral with respect to a closed set. Section~\ref{sec:bundle-form-operations} studies smooth potential perturbations, restrictions to subbundles, and jump perturbations. A restriction can be recovered by increasing the potential on the orthogonal complement of the subbundle (Theorem~\ref{thm:restriction-as-supremum}). Section~\ref{sec:riemannian-bundle-forms} illustrates how the boundary condition determines the fibers. The Neumann connection energy retains the full fiber at the boundary, whereas the mixed condition retains only the tangential fiber. The Hodge energy naturally corresponds to mixed boundary conditions (Section~\ref{subsec:hodge-energy-mixed-boundary}).

\section*{Prior work}
The name ``Kato form'' refers to Kato's distributional inequality \cite{Kato1972}. Hess, Schrader and Uhlenbrock \cite{HessSchraderUhlenbrock1977} characterized semigroup domination by a Kato inequality for the generators. Ouhabaz \cite{Ouhabaz1996} gave a form criterion on scalar $L^2$. Shigekawa used the carr\'e du champ to study domination for vector-valued semigroups \cite{Shigekawa1997} and for the Hodge heat semigroup on manifolds with boundary \cite{Shigekawa2000}. Manavi, Vogt and Voigt \cite{ManaviVogtVoigt2005} extended the form criterion to Hilbert-valued $L^2$, and Lenz, Schmidt and Wirth \cite{LenzSchmidtWirth2020} proved it in the setting of ordered Hilbert spaces.

\section{Kato forms and quasi-realization}\label{sec:bundle-forms-realization}
\subsection{Kato forms}
\begin{dfn}
A bundle form on $(X,m)$ is a pair $(\breve{\mathcal E},\mathcal H)$,
where $\mathcal H$ is an almost Hilbert bundle over $(X,m)$ and
$
    \breve{\mathcal E}\colon
    \D(\breve{\mathcal E})\times\D(\breve{\mathcal E})
    \rightarrow\mathbb R
$
is a densely defined, lower-semibounded, closed symmetric
bilinear form on $\L^2(m;\mathcal H)$.
We denote by $\breve L$, $(\breve T_t)_{t\geq0}$,
and $(\breve G_\lambda)$ the generator, semigroup,
and resolvent associated with $\breve{\mathcal E}$,
respectively, with $\lambda$ taken sufficiently large.
\end{dfn}
When $\mathcal H$ is clear from the context, we simply refer to
$\breve{\mathcal E}$ as a bundle form.
\begin{thm}[Kato inequality and semigroup domination]\label{thm:kato-domination}
Let $\breve{\mathcal{E}}$ be a non-negative bundle form on $(X,m)$. Then the following conditions are equivalent. 
\begin{enumerate} 
\item Kato property\;:\;for every $\alpha\in\D( \breve{\mathcal{E}})$ and $f\in\D(\mathcal{E})_+$ with $f\leq|\alpha|$, one has $ |\alpha|\in\D(\mathcal{E})$ and $\frac{f}{|\alpha|}\alpha\in\D(\breve{\mathcal{E}})$. Moreover, 
\[ \mathcal{E}(|\alpha|,f)\leq \breve{\mathcal{E}} \left( \alpha,\frac{f}{|\alpha|}\alpha \right)\] 
\item For every $\alpha\in\D(\breve{L})$ and $f\in\D(L)_+$, \[ \int_X f\,\frac{\<\breve{L}\alpha,\alpha\>}{|\alpha|} \,dm \leq \int_X (L f)|\alpha|\,dm \] 
\item Semigroup domination\;:\;for every $t>0$ and $\alpha\in\L^2(m;\mathcal H)$, \[|\breve{T}_t\alpha|\leq T_t|\alpha|. \] 
\item Resolvent domination\;:\;for every $\lambda>0$ and $\alpha\in\L^2(m;\mathcal H)$, \[|\lambda\breve{G}_\lambda\alpha|\leq\lambda G_{\lambda}|\alpha|. \] \end{enumerate} Here, the fraction is understood to be zero on $\{|\alpha|=0\}$.
\end{thm}
\begin{proof}\;
\begin{itemize}
    \item $1\Rightarrow 4$.
Set
$h\coloneqq |\lambda\breve{G}_\lambda\alpha|-\lambda G_\lambda|\alpha|.$
Since $\lambda\breve{G}_\lambda\alpha\in\D(\breve{\mathcal E})$, $1$ implies that
$|\lambda\breve{G}_\lambda\alpha|\in\D(\mathcal E)$. Hence
$h,h_+\in\D(\mathcal E)$, and
$0\leq h_+\leq|\lambda\breve{G}_\lambda\alpha|.$
By 1 and the correspondence between the resolvent and the form, we obtain
\begin{align*}
    \mathcal{E}_\lambda
    (|\lambda\breve{G}_\lambda\alpha|,h_+)
    &\leq\breve{\mathcal{E}}_\lambda\left(\lambda\breve{G}_\lambda\alpha,
        \frac{h_+}{|\lambda\breve{G}_\lambda\alpha|}
        \lambda\breve{G}_\lambda\alpha\right) \\
    &=\lambda\int_X\left\langle\alpha,\frac{h_+}{|\lambda\breve{G}_\lambda\alpha|}\lambda\breve{G}_\lambda\alpha\right\rangle\,dm \\
    &\leq \lambda\int_X h_+|\alpha|\,dm \\
    &=\mathcal{E}_\lambda(\lambda G_\lambda|\alpha|,h_+).
\end{align*}
Thus,
$\mathcal{E}_\lambda(h,h_+)\leq0.$
On the other hand, by the lattice property of $\mathcal E$,
\[\mathcal{E}_\lambda(h,h_+)\geq\mathcal{E}_\lambda(h_+,h_+)
\geq0.
\]
Hence $h_+=0$. Therefore
$|\lambda\breve{G}_\lambda\alpha|\leq\lambda G_\lambda|\alpha|.
$
\item  $4\Rightarrow 2$. By $4$ and the symmetry of $G_\lambda$,
\[\int_X f\,|\lambda\breve G_\lambda\alpha|\,dm\leq\int_X (\lambda G_\lambda f)|\alpha|\,dm.
\]
Since
$|\lambda\breve G_\lambda\alpha|-|\alpha|\geq\<\lambda\breve G_\lambda\alpha-\alpha,\alpha\>/|\alpha|,
$ we obtain
\[\int_X f\,\frac{\<\lambda\breve G_\lambda\alpha-\alpha,\alpha\>}{|\alpha|}\,dm\leq\int_X(\lambda G_\lambda f-f)|\alpha|\,dm.
\]
Multiplying by $\lambda$ and letting $\lambda\to\infty$, we conclude that
\[\int_X f\,\frac{\<\breve L\alpha,\alpha\>}{|\alpha|}\,dm
\leq\int_X (L f)|\alpha|\,dm.
\]
\item $2\Rightarrow 3$.
It is enough to prove the assertion for $\alpha\in\D(\breve L)$.
Fix $t>0$ and let $f\in\D(L)_+$ be arbitrary. Set
\[
    \Phi(s)\coloneqq
    \int_X T_{t-s}f\,|\breve T_s\alpha|\,dm,
    \qquad s\in[0,t].
\]
We have $\breve T_s\alpha\in\D(\breve L)$ and
$T_{t-s}f\in\D(L)_+$.
Moreover, the maps
$s\mapsto T_{t-s}f$ and $s\mapsto|\breve T_s\alpha|$
are absolutely continuous as $\L^2(m)$-valued maps, and
\[
    \frac{d}{ds}T_{t-s}f=-L T_{t-s}f,
    \qquad
    \frac{d}{ds}|\breve T_s\alpha|
    =
    \frac{\langle\breve L\breve T_s\alpha,
    \breve T_s\alpha\rangle}{|\breve T_s\alpha|}
\]
for a.e.\ $s\in(0,t)$.
Thus $\Phi$ is absolutely continuous and, by $2$,
\begin{align*}
    \Phi'(s)
    &=
    -\int_X (L T_{t-s}f)|\breve T_s\alpha|\,dm
    +\int_X T_{t-s}f\,
    \frac{\langle\breve L\breve T_s\alpha,
    \breve T_s\alpha\rangle}{|\breve T_s\alpha|}\,dm
    \leq0
\end{align*}
for a.e.\ $s\in(0,t)$.
Therefore,
\[
    \int_X f|\breve T_t\alpha|\,dm
    =\Phi(t)\leq\Phi(0)
    =\int_X fT_t|\alpha|\,dm.
\]
Since this holds for every $f\in\D(L)_+$, we obtain
$|\breve T_t\alpha|\leq T_t|\alpha|$.
\item $3\Rightarrow 1$. 
Let $\alpha\in\D(\breve{\mathcal E})$ be arbitrary. By $3$, for every $t>0$,
\begin{align*}
    \mathcal E^{(t)}(|\alpha|)
    &=\frac1t\left(
        \|\alpha\|_2^2
        -\langle T_t|\alpha|,|\alpha|\rangle_2
    \right)\\
    &\leq\frac1t\left(
        \|\alpha\|_2^2
        -\langle\breve T_t\alpha,\alpha\rangle_2
    \right)\\
    &=\breve{\mathcal E}^{(t)}(\alpha)
    \leq\breve{\mathcal E}(\alpha).
\end{align*}
Thus the semigroup characterization of the form domain yields
$|\alpha|\in\D(\mathcal E)$. Next, let $f\in\D(\mathcal E)_+$ satisfy $f\leq|\alpha|$.
Since
\[
    \left|\alpha-\frac f{|\alpha|}\alpha\right|=|\alpha|-f,
    \qquad
    \left|\alpha+\frac f{|\alpha|}\alpha\right|=|\alpha|+f,
\]
condition $3$ yields
\[
    \left\langle
        \breve T_t\left(\alpha-\frac f{|\alpha|}\alpha\right),
        \alpha+\frac f{|\alpha|}\alpha
    \right\rangle_2
    \leq
    \langle T_t(|\alpha|-f),|\alpha|+f\rangle_2
\]
for every $t>0$.
Together with
\[
    \left\langle
        \alpha-\frac f{|\alpha|}\alpha,
        \alpha+\frac f{|\alpha|}\alpha
    \right\rangle_2
    =
    \langle|\alpha|-f,|\alpha|+f\rangle_2,
\]
this implies
\[
    \breve{\mathcal E}^{(t)}
    \left(
        \alpha-\frac f{|\alpha|}\alpha,
        \alpha+\frac f{|\alpha|}\alpha
    \right)
    \geq
    \mathcal E^{(t)}(|\alpha|-f,|\alpha|+f).
\]
By symmetry, we obtain
\[
    \breve{\mathcal E}^{(t)}(\alpha)
    -
    \breve{\mathcal E}^{(t)}
    \left(\frac f{|\alpha|}\alpha\right)
    \geq
    \mathcal E^{(t)}(|\alpha|)
    -
    \mathcal E^{(t)}(f).
\]
Consequently,
\begin{align*}
    \breve{\mathcal E}^{(t)}
    \left(\frac f{|\alpha|}\alpha\right)
    &\leq
    \breve{\mathcal E}^{(t)}(\alpha)
    -\mathcal E^{(t)}(|\alpha|)
    +\mathcal E^{(t)}(f)\\
    &\leq
    \breve{\mathcal E}(\alpha)+\mathcal E(f).
\end{align*}
The right-hand side is independent of $t>0$.
Thus the semigroup characterization of the form domain yields
$
    \frac f{|\alpha|}\alpha\in\D(\breve{\mathcal E}).
$
Finally, by $3$,
\begin{align*}
    \mathcal E^{(t)}(|\alpha|,f)
    &=\frac1t\left(
        \langle|\alpha|,f\rangle_2
        -\langle T_t|\alpha|,f\rangle_2
    \right)\\
    &\leq\frac1t\left(
        \left\langle\alpha,\frac f{|\alpha|}\alpha\right\rangle_2
        -\left\langle
            \breve T_t\alpha,\frac f{|\alpha|}\alpha
        \right\rangle_2
    \right)\\
    &=\breve{\mathcal E}^{(t)}
    \left(\alpha,\frac f{|\alpha|}\alpha\right).
\end{align*}
Since all arguments belong to the corresponding form domains,
letting $t\downarrow0$ gives
\[
    \mathcal E(|\alpha|,f)
    \leq
    \breve{\mathcal E}
    \left(\alpha,\frac f{|\alpha|}\alpha\right).
\]
This proves $1$.
\end{itemize}
\end{proof}
\begin{dfn}\label{dfn:kato-form}
   $\breve{\mathcal E}$ is called a Kato form on $(X,m,\mathcal E)$ if it satisfies the equivalent conditions of Theorem~\ref{thm:kato-domination}.
\end{dfn}
Define 
\begin{align*}
\D(\breve{\mathcal E})^\infty
&\coloneqq
\{\alpha\in\D(\breve{\mathcal E})\mid |\alpha|\in \L^\infty(m)\},\\
\D(\breve L)^\infty
&\coloneqq
\{\alpha\in\D(\breve L)\mid |\alpha|\in \L^\infty(m)\}.
\end{align*}
\begin{prp}
Let $\breve{\mathcal E}$ be a Kato form on $(X,m,\mathcal E)$. For every
$f\in\D(\mathcal E)^\infty$ and
$\alpha\in\D(\breve{\mathcal E})^\infty$, one has
$f\alpha\in\D(\breve{\mathcal E})^\infty$ and
\[
\breve{\mathcal E}(f\alpha)^{1/2}
\leq
\|f\|_{\infty}
\breve{\mathcal E}(\alpha)^{1/2}
+
\||\alpha|\|_{\infty}
\mathcal E(f)^{1/2}.
\]
\end{prp}
\begin{proof}
By homogeneity, we may assume that $\|f\|_{\infty}=1$. Since $|\alpha|\in\D(\mathcal E)^\infty$ and
$\D(\mathcal E)^\infty$ is an algebra, we have
$f_\pm|\alpha|\in\D(\mathcal E)_+$.
Moreover, $f_\pm|\alpha|\leq|\alpha|$.
Hence the Kato property yields
$f_\pm\alpha\in\D(\breve{\mathcal E})$.
Therefore, $f\alpha=f_+\alpha-f_-\alpha\in\D(\breve{\mathcal E})^\infty$. For every $\varepsilon>0$, we have
\[
0\leq
\frac{\varepsilon}{2+\varepsilon}(1+\varepsilon-f)|\alpha|
\leq
(1+\varepsilon+f)|\alpha|
=
|(1+\varepsilon+f)\alpha|.
\]
Applying the Kato property and cancelling the positive factor
$\varepsilon/(2+\varepsilon)$, we obtain
\[
\mathcal E\bigl((1+\varepsilon+f)|\alpha|,
(1+\varepsilon-f)|\alpha|\bigr)
\leq
\breve{\mathcal E}\bigl((1+\varepsilon+f)\alpha,
(1+\varepsilon-f)\alpha\bigr).
\]
By symmetry, expanding both sides and letting $\varepsilon\downarrow0$ yields
\[
\breve{\mathcal E}(f\alpha)
\leq
\mathcal E(f|\alpha|)
+\breve{\mathcal E}(\alpha)
-\mathcal E(|\alpha|).
\]
The scalar Leibniz inequality and the Kato property therefore give
\begin{align*}
\breve{\mathcal E}(f\alpha)
&\leq
\left(
\mathcal E(|\alpha|)^{1/2}
+\||\alpha|\|_{\infty}\mathcal E(f)^{1/2}
\right)^2
+\breve{\mathcal E}(\alpha)
-\mathcal E(|\alpha|)\\
&\leq
\left(
\breve{\mathcal E}(\alpha)^{1/2}
+\||\alpha|\|_{\infty}\mathcal E(f)^{1/2}
\right)^2.
\end{align*}
Taking square roots and rescaling $f$ proves the assertion.
\end{proof}

\subsection{Quasi-realization theorem}
\begin{dfn}\label{dfn:quasi-realization}
Let $\breve{\mathcal E}$ be a bundle form on $(X,m)$.
A quasi-realization of $\breve{\mathcal E}$ is a pair
$(\widetilde{\mathcal H},J)$ consisting of a quasi-Hilbert
bundle $\widetilde{\mathcal H}$ over $(X,m,\mathcal E)$ and an isometric
almost Hilbert bundle isomorphism
$
    J\colon \widetilde{\mathcal H}\rightarrow \mathcal H
$
satisfying the following properties:
\begin{enumerate}
    \item For every $\alpha\in\D(\breve{\mathcal E})$,
    there exists
    $
        \tilde\alpha\in
        \Gamma^{\mathrm{qc}}(\widetilde{\mathcal H}),
    $
    unique up to q.e.\ equality, such that
    \[
        J\tilde\alpha=\alpha
        \quad m\text{-a.e.}
    \]
    \item There exists a countable family
    $(\alpha_k)_{k\geq1}\subset\D(\breve{\mathcal E})$,
    dense with respect to the $\breve{\mathcal E}_1$-norm,
    such that $(\tilde\alpha_k)_{k\geq1}$ is a
    quasi-fundamental family representing the quasi-structure
    of $\widetilde{\mathcal H}$.
\end{enumerate}
\end{dfn}
\begin{thm}[Uniqueness of quasi-realization]\label{thm:quasi-realization-uniqueness}
If $(\widetilde{\mathcal H},J)$ and $(\widetilde{\mathcal H}',J')$ are
quasi-realizations of $\breve{\mathcal E}$, then there exists
an isometric quasi-Hilbert bundle isomorphism
$
    \Phi\colon \widetilde{\mathcal H}\rightarrow \widetilde{\mathcal H}',
$
unique up to q.e.\ equality, such that the following
diagram commutes $m$-a.e.:
\[
\begin{tikzcd}
    \widetilde{\mathcal H}
    \arrow[rr,"\Phi"]
    \arrow[dr,"J"']
    &&
    \widetilde{\mathcal H}'
    \arrow[dl,"J'"]
    \\
    & \mathcal H &
\end{tikzcd}
\]
\end{thm}
\begin{proof}\;
\begin{itemize}
    \item We first compare the quasi-continuous representatives in the two quasi-realizations.
    Let $(\widetilde{\mathcal H},J)$ and $(\widetilde{\mathcal H}',J')$ be two
    quasi-realizations of $\breve{\mathcal E}$.
        For each $\alpha\in\D(\breve{\mathcal E})$, denote its quasi-continuous representatives in $\widetilde{\mathcal H}$ and $\widetilde{\mathcal H}'$ by $\tilde\alpha$ and $\tilde\alpha'$, respectively.
    Choose countable families $(\alpha_k)_{k\geq1}$ and $(\beta_k)_{k\geq1}$ as in Property~2 for the respective pairs.
    We first observe that, for every $\alpha,\beta\in\D(\breve{\mathcal E})$ and $c\in\mathbb R$,
    \[
        \widetilde{\alpha+c\beta}
        =
        \tilde\alpha+c\tilde\beta,
        \qquad
        \widetilde{\alpha+c\beta}^{\,\prime}
        =
        \tilde\alpha'+c\tilde\beta'
        \quad\text{q.e.}
    \]
    Indeed, for every $k\geq1$, the identity
    \[
        \langle\widetilde{\alpha+c\beta},\tilde\alpha_k\rangle
        =
        \langle\tilde\alpha,\tilde\alpha_k\rangle
        +
        c\langle\tilde\beta,\tilde\alpha_k\rangle
    \]
    holds $m$-a.e.\ by Property~1 and the linearity and isometry of $J$.
    Both sides are quasi-continuous, so the identity holds q.e.
    Countability and the generating property of $(\tilde\alpha_k)_{k\geq1}$ give the first assertion.
    The second follows in the same way using $(\tilde\beta_k')_{k\geq1}$ and $J'$.
    Moreover, since $J$ and $J'$ are isometric $m$-a.e., scalar quasi-continuous uniqueness gives
    \[
        |\tilde\alpha|
        =
        \widetilde{|\alpha|}
        =
        |\tilde\alpha'|
        \quad\text{q.e.}
    \]
    Consequently, polarization yields
    \[
        \langle\tilde\alpha,\tilde\beta\rangle
        =
        \frac14\bigl(
            \widetilde{|\alpha+\beta|}^{\,2}
            -
            \widetilde{|\alpha-\beta|}^{\,2}
        \bigr)
        =
        \langle\tilde\alpha',\tilde\beta'\rangle
        \quad\text{q.e.}
    \]
    In particular, these inner products are quasi-continuous.

    \item We use their inner products to construct an isometry between the fibers for q.e.\ $x$.
    Let $(\gamma_n)_{n\geq1}$ enumerate
    $\{\alpha_k\mid k\geq1\}\cup\{\beta_k\mid k\geq1\}$.
    The families $(\tilde\gamma_n)_{n\geq1}$ and $(\tilde\gamma_n')_{n\geq1}$ are quasi-fundamental families representing the quasi-structures of $\widetilde{\mathcal H}$ and $\widetilde{\mathcal H}'$, respectively.
    Indeed, each contains the corresponding original quasi-fundamental family, and all the required mixed inner products are quasi-continuous by the preceding identity.
    By countability, there exists an exceptional set $N$ such that, for every $x\in X\setminus N$,
    \[
        \langle\tilde\gamma_n(x),\tilde\gamma_l(x)\rangle
        =
        \langle\tilde\gamma_n'(x),\tilde\gamma_l'(x)\rangle
        \qquad\text{for every }n,l\geq1,
    \]
    and both families have dense linear span in their respective fibers.
    For such $x$, the assignment
    \[
        \sum_{n=1}^{r}c_n\tilde\gamma_n(x)
        \longmapsto
        \sum_{n=1}^{r}c_n\tilde\gamma_n'(x),
        \qquad r\geq1,\quad c_1,\ldots,c_r\in\mathbb R,
    \]
    is therefore well-defined and isometric.
    Since its domain and range are dense, it extends uniquely to an isometric isomorphism
    \[
        \Phi_x\colon \widetilde{\mathcal H}_x\longrightarrow \widetilde{\mathcal H}_x'.
    \]
    Set $\Phi_x=0$ for $x\in N$.

    \item We now show that $\Phi$ is a quasi-Hilbert bundle isomorphism and that $J'\Phi=J$ $m$-a.e.
    If $\xi\in\Gamma^{\mathrm{qc}}(\widetilde{\mathcal H})$, then
    \[
        |\Phi\xi|=|\xi|,
        \qquad
        \langle\Phi\xi,\tilde\gamma_n'\rangle
        =
        \langle\xi,\tilde\gamma_n\rangle
        \quad\text{q.e.}
    \]
    Since the two families represent the respective quasi-structures, these identities show that $\Phi\xi$ is quasi-continuous.
    The same argument applies to $\Phi^{-1}$.
    Hence $\Phi$ is an isometric quasi-Hilbert bundle isomorphism.
    For every $n\geq1$, we have
    \[
        J'\Phi\tilde\gamma_n
        =
        J'\tilde\gamma_n'
        =
        \gamma_n
        =
        J\tilde\gamma_n
        \quad m\text{-a.e.}
    \]
    Outside a common $m$-null set, these identities hold simultaneously, the family $(\tilde\gamma_n(x))_{n\geq1}$ has dense linear span, and $J_x$ and $J_x'$ are isometric.
    By continuity, it follows that
    \[
        J_x'\circ\Phi_x=J_x
        \quad\text{for }m\text{-a.e.\ }x\in X.
    \]
    Thus the diagram commutes $m$-a.e.

    \item We show that the commutation relation determines $\Phi$ up to q.e.\ equality.
    Finally, let $\Psi$ be another isometric quasi-Hilbert bundle isomorphism satisfying $J'\Psi=J$ $m$-a.e.
    For every $n\geq1$, the section $\Psi\tilde\gamma_n$ is quasi-continuous and satisfies
    \[
        J'\Psi\tilde\gamma_n
        =
        J\tilde\gamma_n
        =
        \gamma_n
        \quad m\text{-a.e.}
    \]
    Property~1 for $(\widetilde{\mathcal H}',J')$ therefore gives
    \[
        \Psi\tilde\gamma_n
        =
        \tilde\gamma_n'
        =
        \Phi\tilde\gamma_n
        \quad\text{q.e.}
    \]
    Taking a common exceptional set and using the dense linear span of $(\tilde\gamma_n(x))_{n\geq1}$, we conclude that $\Psi_x=\Phi_x$ for q.e.\ $x\in X$.
\end{itemize}
\end{proof}

\begin{thm}[Existence of quasi-realization]\label{thm:quasi-realization-existence}
Let $\breve{\mathcal E}$ be a Kato form on $(X,m,\mathcal E)$.
Then $\breve{\mathcal E}$ admits a quasi-realization
$(\widetilde{\mathcal H},J)$.
\end{thm}
\begin{proof}\;
\begin{itemize}
    \item We first construct the quasi-Hilbert bundle $\widetilde{\mathcal H}$.
    The Kato property implies that, for every $\alpha\in\D(\breve{\mathcal E})$,
    \[
        |\alpha|\in\D(\mathcal E),
        \qquad
        \mathcal E_1(|\alpha|)\leq\breve{\mathcal E}_1(\alpha).
    \]
    For each $\alpha\in\D(\breve{\mathcal E})$, choose a quasi-continuous version
    $\widetilde{|\alpha|}$ of $|\alpha|$.
    We shall repeatedly use the fact that an equality or inequality between
    quasi-continuous scalar functions which holds $m$-a.e.\ also holds q.e.
    Since $\L^2(m;\mathcal H)$ is separable and $\breve{\mathcal E}$ is closed,
    $\D(\breve{\mathcal E})$ is separable with respect to the
    $\breve{\mathcal E}_1$-norm.
    Choose a countable $\mathbb Q$-linear subspace
    $C\subset\D(\breve{\mathcal E})$ which is dense in this norm.
    By countability, there exists an exceptional set $N$ such that,
    for every $x\in X\setminus N$, all $\widetilde{|\alpha|}(x)$,
    $\alpha\in C$, are finite and nonnegative, $\widetilde{|0|}(x)=0$, and
    \begin{align*}
        \widetilde{|q\alpha|}(x)
        &=|q|\widetilde{|\alpha|}(x),\\
        \widetilde{|\alpha+\beta|}(x)
        &\leq\widetilde{|\alpha|}(x)+\widetilde{|\beta|}(x),\\
        \widetilde{|\alpha+\beta|}(x)^2+\widetilde{|\alpha-\beta|}(x)^2
        &=2\widetilde{|\alpha|}(x)^2+2\widetilde{|\beta|}(x)^2
    \end{align*}
    hold simultaneously for every $\alpha,\beta\in C$ and $q\in\mathbb Q$.
    For $x\in X\setminus N$, set
    \[
        N_x:=\{\alpha\in C\mid\widetilde{|\alpha|}(x)=0\}.
    \]
    The quotient $C/N_x$ is a pre-Hilbert space over $\mathbb Q$
    with inner product
    \[
        \langle[\alpha],[\beta]\rangle_x
        :=\frac14\bigl(
            \widetilde{|\alpha+\beta|}(x)^2
            -
            \widetilde{|\alpha-\beta|}(x)^2
        \bigr).
    \]
    Let $\widetilde{\mathcal H}_x$ be its completion, with scalar multiplication
    extended from $\mathbb Q$ to $\mathbb R$.
    For $x\in N$, set $\widetilde{\mathcal H}_x:=\{0\}$.
    For $\alpha\in C$, define $\tilde\alpha(x)$ to be the image of $\alpha$
    in $\widetilde{\mathcal H}_x$ when $x\notin N$, and zero otherwise.
    Enumerate $C$ as $(\alpha_k)_{k\geq1}$.
    By construction,
    \[
        \langle\tilde\alpha_k,\tilde\alpha_l\rangle
        =
        \frac14\bigl(
            \widetilde{|\alpha_k+\alpha_l|}^{\,2}
            -
            \widetilde{|\alpha_k-\alpha_l|}^{\,2}
        \bigr)
        \quad\text{q.e.},
    \]
    so these inner products are quasi-continuous.
    Moreover,
    \[
        \widetilde{\mathcal H}_x
        =
        \overline{\operatorname{Span}
        \{\tilde\alpha_k(x)\mid k\geq1\}}
        \quad\text{for q.e.\ }x\in X.
    \]
    Thus $(\tilde\alpha_k)_{k\geq1}$ is a quasi-fundamental family.
    We equip $\{\widetilde{\mathcal H}_x\}_{x\in X}$ with the quasi-structure
    represented by this family and denote the resulting quasi-Hilbert bundle
    by $\widetilde{\mathcal H}$.

    \item We next construct the realization map $J$.
    Choose almost measurable representatives in $\mathcal H$ of the elements of $C$.
    Outside a single $m$-null set, these representatives respect all
    $\mathbb Q$-linear relations in $C$ and satisfy
    \[
        |\alpha(x)|=\widetilde{|\alpha|}(x)
        \qquad\text{for every }\alpha\in C.
    \]
    Consequently, for $m$-a.e.\ $x$, the assignment
    $
        \tilde\alpha(x)\mapsto\alpha(x),\;\alpha\in C,
    $
    extends uniquely to an isometric embedding
    \[
        J_x\colon \widetilde{\mathcal H}_x\longrightarrow \mathcal H_x
    \]
    whose range is the closed linear span of $\{\alpha_k(x)\mid k\geq1\}$.
    This range equals $\mathcal H_x$ for $m$-a.e.\ $x$.
    Indeed, $C$ is dense in $\L^2(m;\mathcal H)$.
    Let $(\eta_j)_{j\geq1}$ be an almost fundamental family of $\mathcal H$,
    and choose a strictly positive $w\in\L^2(m)$.
    For each $j$, the section
    \[
        \frac{w}{1+|\eta_j|}\eta_j
    \]
    belongs to $\L^2(m;\mathcal H)$ and can therefore be approximated in $\L^2(m;\mathcal H)$
    by elements of $C$.
    Passing to an almost everywhere convergent subsequence shows that
    $\eta_j(x)$ belongs to the closed linear span of
    $\{\alpha_k(x)\mid k\geq1\}$ for $m$-a.e.\ $x$.
    Countability and the generating property of $(\eta_j)$ prove the assertion.
    Define $J_x$ to be zero on the remaining $m$-null set.
    Since
    \[
        J\tilde\alpha_k=\alpha_k
        \quad m\text{-a.e.},
        \qquad k\geq1,
    \]
    and both families are almost fundamental for their respective almost
    Hilbert bundle structures, $J$ and its almost everywhere defined inverse
    are almost measurable.
    Hence $J$ is an isometric almost Hilbert bundle isomorphism.

    \item We then construct a quasi-continuous representative of each element of $\D(\breve{\mathcal E})$.
    Fix $\alpha\in\D(\breve{\mathcal E})$.
    Choose $(\beta_n)_{n\geq1}\subset C$ converging to $\alpha$ in the
    $\breve{\mathcal E}_1$-norm.
    The Kato estimate gives
    \[
        \mathcal E_1(|\beta_n-\alpha|)
        \leq
        \breve{\mathcal E}_1(\beta_n-\alpha)
        \longrightarrow0.
    \]
    By the subsequence convergence theorem for quasi-regular Dirichlet forms,
    after passing to a subsequence,
    \[
        \widetilde{|\beta_n-\alpha|}\longrightarrow0
        \quad\text{q.e.}
    \]
    For every $n,j\geq1$, we have
    \[
        |\tilde\beta_n-\tilde\beta_j|
        =
        \widetilde{|\beta_n-\beta_j|}
        \leq
        \widetilde{|\beta_n-\alpha|}
        +
        \widetilde{|\beta_j-\alpha|}
        \quad\text{q.e.}
    \]
    Taking a common exceptional set for these countably many relations,
    we see that $(\tilde\beta_n(x))_{n\geq1}$ is Cauchy in
    $\widetilde{\mathcal H}_x$ for q.e.\ $x$.
    Define
    $
        \tilde\alpha(x):=\lim_{n\to\infty}\tilde\beta_n(x)
    $
    where this limit exists, and set $\tilde\alpha(x)=0$ on the remaining
    exceptional set.
    We verify that $\tilde\alpha$ is quasi-continuous.
    For every $\gamma\in C$,
    \[
        \bigl|
            \widetilde{|\beta_n+\gamma|}
            -
            \widetilde{|\alpha+\gamma|}
        \bigr|
        \leq
        \widetilde{|\beta_n-\alpha|}
        \quad\text{q.e.}
    \]
    Passing to the limit, first with $\gamma=0$ and then with
    $\gamma=\pm\alpha_k$, gives
    \[
        |\tilde\alpha|=\widetilde{|\alpha|},
        \qquad
        |\tilde\alpha\pm\tilde\alpha_k|
        =
        \widetilde{|\alpha\pm\alpha_k|}
        \quad\text{q.e.}
    \]
    Therefore,
    \[
        \langle\tilde\alpha,\tilde\alpha_k\rangle
        =
        \frac14\bigl(
            \widetilde{|\alpha+\alpha_k|}^{\,2}
            -
            \widetilde{|\alpha-\alpha_k|}^{\,2}
        \bigr)
        \quad\text{q.e.}
    \]
    Both $|\tilde\alpha|$ and all these inner products are quasi-continuous.
    Thus $\tilde\alpha\in\Gamma^{\mathrm{qc}}(\widetilde{\mathcal H})$.
    Since $\widetilde{|\beta_n-\alpha|}=|\beta_n-\alpha|$ $m$-a.e.,
    we also have $\beta_n(x)\to\alpha(x)$ for $m$-a.e.\ $x$.
    The continuity of $J_x$ consequently yields
    \[
        J_x\tilde\alpha(x)
        =
        \lim_{n\to\infty}J_x\tilde\beta_n(x)
        =
        \lim_{n\to\infty}\beta_n(x)
        =
        \alpha(x)
        \quad m\text{-a.e.}
    \]

    \item We show that this representative is unique up to q.e.\ equality.
    Finally, suppose that
    $\xi\in\Gamma^{\mathrm{qc}}(\widetilde{\mathcal H})$
    also satisfies $J\xi=\alpha$ $m$-a.e.
    For every $k\geq1$, the quasi-continuous functions
    $\langle\xi,\tilde\alpha_k\rangle$ and
    $\langle\tilde\alpha,\tilde\alpha_k\rangle$
    agree $m$-a.e., hence q.e.
    By countability and the generating property of
    $(\tilde\alpha_k)_{k\geq1}$, it follows that $\xi=\tilde\alpha$ q.e.
    This also proves independence of the approximating sequence.
    Property~1 is now established, and the enumeration
    $(\alpha_k)_{k\geq1}$ of $C$ satisfies Property~2 by construction.
\end{itemize}
\end{proof}
For a Kato form $\breve{\mathcal E}$, we write
$(\mathcal H^{\breve{\mathcal E}},J)$ for a quasi-realization.
Theorems~\ref{thm:quasi-realization-existence} and~\ref{thm:quasi-realization-uniqueness} determine this pair up to a unique
isometric quasi-Hilbert bundle isomorphism intertwining
the realization maps.
In what follows, we assume that $\breve{\mathcal E}$ is a Kato form
on $(X,m,\mathcal E)$ and fix such a quasi-realization
$(\mathcal H^{\breve{\mathcal E}},J)$.
\begin{cor}
Let $(\alpha_k)_{k\geq1}\subset\D(\breve{\mathcal E})$
be a countable family which is dense with respect to
the $\breve{\mathcal E}_1$-norm.
Then $(\tilde\alpha_k)_{k\geq1}$ is a quasi-fundamental
family representing the quasi-structure of
$\mathcal H^{\breve{\mathcal E}}$.
\end{cor}
\begin{proof}
Let $(\beta_j)_{j\geq1}\subset\D(\breve{\mathcal E})$
be a family as in Property~2 of Definition~\ref{dfn:quasi-realization}.
For each $j\geq1$, choose a sequence
$(\alpha_{k(j,n)})_{n\geq1}$ converging to $\beta_j$
in the $\breve{\mathcal E}_1$-norm.
The Kato estimate gives
\[
    \mathcal E_1(|\alpha_{k(j,n)}-\beta_j|)
    \leq
    \breve{\mathcal E}_1(\alpha_{k(j,n)}-\beta_j)
    \longrightarrow0.
\]
Passing to a subsequence for each $j$, we obtain
\[
    |\tilde\alpha_{k(j,n)}-\tilde\beta_j|
    =
    \widetilde{|\alpha_{k(j,n)}-\beta_j|}
    \longrightarrow0
    \quad\text{q.e.}
\]
By countability, these convergences hold outside a common
exceptional set. Hence
\[
    \mathcal H_x^{\breve{\mathcal E}}
    =
    \overline{\operatorname{Span}
    \{\tilde\beta_j(x)\mid j\geq1\}}
    \subset
    \overline{\operatorname{Span}
    \{\tilde\alpha_k(x)\mid k\geq1\}}
    \subset \mathcal H_x^{\breve{\mathcal E}}
    \quad\text{for q.e.\ }x.
\]
Since each $\tilde\alpha_k$ is quasi-continuous,
their mutual inner products and their mixed inner products
with $(\tilde\beta_j)_{j\geq1}$ are quasi-continuous.
Thus $(\tilde\alpha_k)_{k\geq1}$ is a quasi-fundamental
family representing the same quasi-structure.
\end{proof}
\begin{cor}
The following statements hold:
\begin{enumerate}
    \item For every $\alpha,\beta\in\D(\breve{\mathcal E})$
    and $c\in\mathbb R$,
    \[
        \widetilde{\alpha+c\beta}
        =
        \tilde\alpha+c\tilde\beta
        \quad\text{q.e.}
    \]

    \item For every $\alpha\in\D(\breve{\mathcal E})$,
    \[
        |\tilde\alpha|
        =
        \widetilde{|\alpha|}
        \quad\text{q.e.}
    \]

    \item For every $\alpha\in\D(\breve{\mathcal E})$
    and $f\in\D(\mathcal E)_+$ with $f\leq|\alpha|$
    $m$-a.e.,
    \[
        \widetilde{\frac f{|\alpha|}\alpha}
        =
        \frac{\tilde f}{|\tilde\alpha|}\tilde\alpha
        \quad\text{q.e.}
    \]

    \item For every $\alpha\in\D(\breve{\mathcal E})^\infty$
    and $f\in\D(\mathcal E)^\infty$,
    \[
        \widetilde{f\alpha}
        =
        \tilde f\,\tilde\alpha
        \quad\text{q.e.}
    \]
\end{enumerate}
\end{cor}
\begin{proof}
Assertions $1$ and $2$ follow from the uniqueness of
quasi-continuous representatives.
For $3$, scalar quasi-continuous uniqueness gives
$0\leq\tilde f\leq|\tilde\alpha|$ q.e.
Thus
$
    \frac{\tilde f}{|\tilde\alpha|}\tilde\alpha
$
is quasi-continuous, with the stated convention at
$|\tilde\alpha|=0$.
Its image under $J$ agrees with
$\frac f{|\alpha|}\alpha$ $m$-a.e.,
so Property~1 gives the desired identity.
For $4$, the algebra property of $\D(\mathcal E)^\infty$
justifies the product
$f\alpha$.
Since $\tilde f\,\tilde\alpha$ is quasi-continuous and
\[
    J(\tilde f\,\tilde\alpha)=f\alpha
    \quad m\text{-a.e.},
\]
the identity follows again from Property~1.
\end{proof}
We conclude this section by defining the resolvent of a vector measure of finite energy integral with respect to a closed set. In Section~\ref{subsec:damped-transports}, we use this resolvent to identify the form associated with damped transport.
\begin{dfn}
\label{dfn:vector-measure-finite-energy}
Let $\mu$ be a smooth measure and let $\alpha$ be a
measurable section of $H^{\breve{\mathcal E}}$.
We say that $\alpha\mu$ is finite energy integral if $|\alpha|\mu$ is of finite energy integral, that is,
if there exists a constant $C\ge 0$ such that
\[
    \int_X |\tilde f|\,|\alpha|\,d\mu
    \le C\,\mathcal E_1(f)^{1/2}
    \qquad \text{for every } f\in \D(\mathcal E).
\]
\end{dfn}
Suppose that $\alpha\mu$ is of finite energy integral. Let $F\subset X$ be closed and let $\lambda>0$. Set
\[
    \D(\breve{\mathcal E})_F
    :=
    \left\{
        \xi\in\D(\breve{\mathcal E})
        \,\middle|\,
        \tilde\xi=0 \text{ q.e.\ on } X\setminus F
    \right\}.
\]
This subspace is closed in the form norm. By the finite energy integral assumption and the Kato property,
\[
    \left|\int_X \langle\alpha,\tilde\xi\rangle\,d\mu\right|
    \le \int_X |\tilde\xi|\,|\alpha|\,d\mu
    \le C\,\mathcal E_1(|\xi|)^{1/2}
    \le C\,\breve{\mathcal E}_1(\xi)^{1/2}
    \qquad \text{for every } \xi\in\D(\breve{\mathcal E})_F.
\]
Since the $\breve{\mathcal E}_\lambda^{1/2}$-norm is equivalent to the
form norm, the Riesz representation theorem gives a unique element
$\breve G_{\lambda,F}(\alpha\mu)\in\D(\breve{\mathcal E})_F$ such that
\[
    \breve{\mathcal E}_\lambda
    \bigl(\breve G_{\lambda,F}(\alpha\mu),\xi\bigr)
    =
    \int_X \langle\alpha,\tilde\xi\rangle\,d\mu
    \qquad \text{for every } \xi\in\D(\breve{\mathcal E})_F.
\]
We call $\breve G_{\lambda,F}(\alpha\mu)$ the $\lambda$-resolvent of
$\alpha\mu$ with respect to $F$. When $F=X$, we simply write
$\breve G_\lambda(\alpha\mu)$.
For $\xi\in L^2(m,H^{\breve{\mathcal E}})$, write
$\breve G_{\lambda,F}\xi:=\breve G_{\lambda,F}(\xi m)$.
Taking $\breve G_{\lambda,F}\xi$ as the test section gives
\[
    \bigl(\breve G_{\lambda,F}(\alpha\mu),\xi\bigr)_2
    =
    \breve{\mathcal E}_\lambda
    \bigl(\breve G_{\lambda,F}(\alpha\mu),\breve G_{\lambda,F}\xi\bigr)
    =
    \int_X
    \bigl\langle\alpha,\widetilde{\breve G_{\lambda,F}\xi}\bigr\rangle\,d\mu .
\]

\clearpage
\section{Basic operations on bundle forms}\label{sec:bundle-form-operations}
Let $\breve{\mathcal E}$ be a Kato form on $(X,m,\mathcal E)$.

\subsection{Potential perturbation}\label{subsec:potential-perturbation}
Let $\mu$ be a smooth measure.
Since $\mu$ does not charge exceptional sets,
the quasi-structure of $\mathcal H^{\breve{\mathcal E}}$ naturally
induces an almost Hilbert bundle structure over $(X,\mu)$.
We regard $\mathcal H^{\breve{\mathcal E}}$ as an almost Hilbert
bundle over $(X,\mu)$ when discussing almost measurability
with respect to $\mu$.
Let $R$ be an almost measurable section of
$\operatorname{End}(\mathcal H^{\breve{\mathcal E}})$ with respect to
$\mu$ such that
\[
    R(x)=R(x)^\dagger
    \qquad \mu\text{-a.e.}
\]
\begin{dfn}
    We call $R\mu$ a smooth potential if the measure
$\|R\|_{\mathrm{op}}\,\mu$ is smooth.
\end{dfn}

\begin{dfn}
Let $R\mu$ be a smooth potential.
The potential perturbation of $\breve{\mathcal E}$ by $R\mu$
is the symmetric form $\breve{\mathcal E}^{R\mu}$ with domain
\[
    \D(\breve{\mathcal E}^{R\mu})
    \coloneqq
    \left\{
        \alpha\in\D(\breve{\mathcal E})
        \,\middle|\,
        \int_X
        \langle |R|\tilde\alpha,\tilde\alpha\rangle\,d\mu
        <\infty
    \right\},
\]
defined by
\[
    \breve{\mathcal E}^{R\mu}(\alpha,\beta)
    \coloneqq
    \breve{\mathcal E}(\alpha,\beta)
    +
    \int_X
    \langle R\tilde\alpha,\tilde\beta\rangle\,d\mu,
    \qquad
    \alpha,\beta\in\D(\breve{\mathcal E}^{R\mu}).
\]
\end{dfn}
\begin{ex}
Assume there exist constants
$a\in[0,1)$ and $b\geq0$ such that
\[
    \int_X
    \langle R_-\tilde\alpha,\tilde\alpha\rangle\,d\mu
    \leq
    a\,\breve{\mathcal E}(\alpha)
    +b\,\|\alpha\|_2^2
    \qquad
    \text{for every }\alpha\in\D(\breve{\mathcal E}).
\]
Then $\breve{\mathcal E}^{R\mu}$
is a bundle form on $\L^2(m;\mathcal H)$.
Moreover,
\[
    \breve{\mathcal E}^{R\mu}(\alpha)
    \geq
    (1-a)\breve{\mathcal E}(\alpha)
    -b\|\alpha\|_2^2
    \geq
    -b\|\alpha\|_2^2
    \qquad
    \text{for every }
    \alpha\in\D(\breve{\mathcal E}^{R\mu}).
\]
\end{ex}
We record the following density property for later use in the proof of Theorem~\ref{thm:damped-transport-association}.
\begin{lem}\label{lem:perturbed-form-nest-density}
Suppose that $\breve{\mathcal E}^{R\mu}$ is a bundle form and let $(F_n)$ be a nest. Then
\[
    \bigcup_{n\geq1}\D(\breve{\mathcal E}^{R\mu})_{F_n}
    \quad\text{is dense in }\D(\breve{\mathcal E}^{R\mu})
    \text{ for the perturbed form norm},
\]
where $\D(\breve{\mathcal E}^{R\mu})_{F_n}$ consists of the sections in $\D(\breve{\mathcal E}^{R\mu})$ vanishing q.e.\ outside $F_n$.
\end{lem}
\begin{proof}
Fix $\beta\in\D(\breve{\mathcal E}^{R\mu})$.
By the nest property, we can approximate $|\beta|$ by functions supported on members of $(F_n)$ in the scalar form norm.
Taking their positive parts and then their minima with $|\beta|$ gives functions $f_j$ supported on members of the nest such that
\[
    0\leq f_j\leq|\beta|,\qquad
    f_j\longrightarrow|\beta|\text{ in }\L^2(m),\qquad
    \sup_j\mathcal E_1(f_j)<\infty .
\]
Set $\beta_j=(f_j/|\beta|)\beta$, with value zero on $\{|\beta|=0\}$ then $\tilde\beta_j$ is a scalar multiple of $\tilde\beta$
with multiplier in $[0,1]$ q.e., so
\[
    \int_X\langle |R|\tilde\beta_j,\tilde\beta_j\rangle\,d\mu
    \leq
    \int_X\langle |R|\tilde\beta,\tilde\beta\rangle\,d\mu .
\]
Thus $(\beta_j)$ is bounded in the perturbed form norm and converges to $\beta$ in $\L^2$. A weakly convergent subsequence in the form domain has limit $\beta$. This proves the asserted density.
\end{proof}

\subsection{Restriction}\label{subsec:restriction-forms}
Let $\widetilde{\mathcal K}$ be a quasi-Hilbert subbundle of
$\mathcal H^{\breve{\mathcal E}}$.
Using $J$, we identify $\L^2(m;\widetilde{\mathcal K})$ with the corresponding
closed subspace of $\L^2(m;\mathcal H)$.

\begin{dfn}
We say that $\widetilde{\mathcal K}$ is compatible with $\breve{\mathcal E}$
if there exists a countable family
$(\alpha_k)_{k\geq1}\subset\D(\breve{\mathcal E})$ such that
\[
\widetilde{\mathcal K}_x
=
\overline{\operatorname{Span}
\left\{
\tilde\alpha_k(x)
\;\middle|\;
k\geq1
\right\}}
\quad\text{q.e.}
\]
\end{dfn}

In what follows, we assume that $\widetilde{\mathcal K}$ is compatible with
$\breve{\mathcal E}$.

\begin{dfn}
The restriction form of $\breve{\mathcal E}$ on $\widetilde{\mathcal K}$
is defined by
\[
\D(\breve{\mathcal E}^{\widetilde{\mathcal K}})
:=
\left\{
\alpha\in\D(\breve{\mathcal E})
\;\middle|\;
\tilde\alpha(x)\in \widetilde{\mathcal K}_x
\quad\text{q.e.}
\right\},
\]
and
\[
\breve{\mathcal E}^{\widetilde{\mathcal K}}(\alpha,\beta)
:=
\breve{\mathcal E}(\alpha,\beta),
\qquad
\alpha,\beta\in\D(\breve{\mathcal E}^{\widetilde{\mathcal K}}).
\]
We regard $\breve{\mathcal E}^{\widetilde{\mathcal K}}$ as a form on
$\L^2(m;\widetilde{\mathcal K})$ through the above identification.
\end{dfn}

\begin{prp}\label{prp:restriction-kato-form}
$\breve{\mathcal E}^{\widetilde{\mathcal K}}$ is a Kato form with respect to $\mathcal E$ on $\L^2(m;\widetilde{\mathcal K})$ Moreover, $(\widetilde{\mathcal K},\operatorname{id})$ is a quasi-realization of $\breve{\mathcal E}^{\widetilde{\mathcal K}}$..
\end{prp}

\begin{proof}
We first prove density.
It suffices to show that $\eta=0$ whenever
$\eta\in\L^2(m;\widetilde{\mathcal K})$ satisfies
$\langle\eta,\alpha\rangle_2=0$ for every
$\alpha\in\D(\breve{\mathcal E}^{\widetilde{\mathcal K}})$.
Fix such an $\eta$.
Choose a generating family $(\alpha_k)_{k\geq1}$ as in the
compatibility assumption.
Replacing each $\alpha_k$ by
$\alpha_k/(1\vee|\alpha_k|)$, the Kato property allows us
to assume that $\alpha_k\in\D(\breve{\mathcal E})^\infty$
without changing the fiberwise closed span.
For every $f\in\D(\mathcal E)^\infty$, the multiplication
property and the identity
$\widetilde{f\alpha_k}=\tilde f\,\tilde\alpha_k$ q.e.\ give
$f\alpha_k\in\D(\breve{\mathcal E}^{\widetilde{\mathcal K}})$.
Hence
\[
    \int_X f\,\langle\eta,\alpha_k\rangle\,dm=0
    \qquad
    \text{for every }f\in\D(\mathcal E)^\infty,\ k\geq1.
\]
Since $\alpha_k$ is bounded,
$\langle\eta,\alpha_k\rangle\in\L^2(m)$.
The density of $\D(\mathcal E)^\infty$ in $\L^2(m)$
therefore yields $\langle\eta,\alpha_k\rangle=0$
$m$-a.e.\ for every $k$.
By countability and the generating property, $\eta=0$
$m$-a.e., proving density in $\L^2(m;\widetilde{\mathcal K})$.
Closedness follows by extracting a q.e.\ convergent
subsequence from a sequence converging in the
$\breve{\mathcal E}_1$-norm and using the closedness
of each fiber $\widetilde{\mathcal K}_x$.
Finally, for
$\alpha\in\D(\breve{\mathcal E}^{\widetilde{\mathcal K}})$
and $f\in\D(\mathcal E)_+$ with $f\leq|\alpha|$,
the quasi-continuous representative of
$\frac{f}{|\alpha|}\alpha$ still takes values in
$\widetilde{\mathcal K}$ q.e.
Thus the domain condition and the Kato inequality
are inherited from $\breve{\mathcal E}$. The quasi-realization assertion follows by adjoining the generators from the compatibility assumption to a countable form-norm dense subset of $\D(\breve{\mathcal E}^{\widetilde{\mathcal K}})$.
\end{proof}

\begin{thm}[Existence of a detecting measure]\label{thm:detecting-measure-existence}
There exists a smooth measure $\mu$ such that,
for every $\alpha\in\D(\breve{\mathcal E})$,
\[
    \tilde\alpha(x)\in\widetilde{\mathcal K}_x
    \quad\mu\text{-a.e.}
    \qquad\Longleftrightarrow\qquad
    \tilde\alpha(x)\in\widetilde{\mathcal K}_x
    \quad\text{q.e.}
\]
\end{thm}
\begin{proof}
    Set $V:=\D(\breve{\mathcal E}^{\widetilde{\mathcal K}})$.
For each finite smooth measure $\nu$, set
\[
    V_\nu
    :=
    \left\{
        \alpha\in\D(\breve{\mathcal E})
        \,\middle|\,
        \tilde\alpha(x)\in\widetilde{\mathcal K}_x
        \quad\text{for $\nu$-a.e.\ }x\in X
    \right\}.
\]
It suffices to find a finite smooth measure $\mu$ such that
$V_\mu=V$.

\begin{itemize}
    \item We first show that
    $
        V=\bigcap_\nu V_\nu.
    $
    Since smooth measures charge no exceptional sets, $V\subset V_\nu$
    for every $\nu$. Conversely, fix
    $\alpha\in\D(\breve{\mathcal E})\setminus V$. Then
    \[
        A_\alpha:=\{x\in X:\tilde\alpha(x)\notin\widetilde{\mathcal K}_x\}
    \]
    is nonexceptional and agrees with a Borel set up to an exceptional set.
    Indeed, q.e.\ the distance from $\tilde\alpha(x)$ to
    $\widetilde{\mathcal K}_x$ is the infimum of
    $|\tilde\alpha(x)-\tilde\beta(x)|$ over the countable family of finite
    rational linear combinations $\beta$ of the generating sections.
    Every nonexceptional Borel set is charged by a finite smooth
    measure \cite[Theorem~2.2.3]{FukushimaOshimaTakeda2011}%
    \cite[Chapter~VI, Theorem~1.2 and Corollary~1.4]{MaRockner1992}.
    Hence $\nu(A_\alpha)>0$ for some such $\nu$, so $\alpha\notin V_\nu$.

    \item We next show that $V=\bigcap_j V_{\nu_j}$ for a countable family
    of finite smooth measures $(\nu_j)$.
    Note that each $V_\nu$ is a closed linear subspace of $\D(\breve{\mathcal E})$
    in the form norm. (Indeed, suppose that $\alpha_k\in V_\nu$ and
    $\alpha_k\to\alpha$ in this norm. Passing to a subsequence, we have
    $\tilde\alpha_k(x)\to\tilde\alpha(x)$ q.e., hence also $\nu$-a.e.
    Since each fiber $\widetilde{\mathcal K}_x$ is closed, $\alpha\in V_\nu$.)
    By the first step,
    \[
        \D(\breve{\mathcal E})\setminus V
        =
        \bigcup_\nu\bigl(\D(\breve{\mathcal E})\setminus V_\nu\bigr).
    \]
    The sets on the right form an open
    cover of $\D(\breve{\mathcal E})\setminus V$. This subspace is second
    countable because $\D(\breve{\mathcal E})$ is a separable metric space.
    Choose a countable subcover and denote the corresponding measures by
    $(\nu_j)$. Taking complements gives $V=\bigcap_j V_{\nu_j}$.
\end{itemize}
Set
\[
    \mu:=\sum_j\frac{2^{-j}}{1+\nu_j(X)}\,\nu_j.
\]
Then $\mu(X)\leq\sum_{j\geq1}2^{-j}=1$, and $\mu$ charges no exceptional
set, so $\mu$ is a finite smooth measure. Since all coefficients
are positive, a Borel set is $\mu$-null if and only if it is $\nu_j$-null
for every $j$. Hence $V_\mu=\bigcap_j V_{\nu_j}=V$, as required.
\end{proof}

\begin{dfn}
    A smooth measure $\mu$ satisfying the equivalence in Theorem~\ref{thm:detecting-measure-existence} is called a
$\widetilde{\mathcal K}$-detecting measure for $\breve{\mathcal E}$.
\end{dfn}
Let $P$ denote the orthogonal projection
onto $\widetilde{\mathcal K}$. 
\begin{thm}[Restriction as a limit of potential perturbations]\label{thm:restriction-as-supremum}
Let $\mu$ be a $\widetilde{\mathcal K}$-detecting measure for
$\breve{\mathcal E}$. Then for every
$\alpha\in\D(\breve{\mathcal E})$, we have
\begin{align*}
    \breve{\mathcal E}^{\widetilde{\mathcal K}}(\alpha)
    &=
    \lim_{n\to\infty}
    \breve{\mathcal E}^{\,n(I-P)\mu}(\alpha)\\
    &=
    \breve{\mathcal E}(\alpha)
    +
    \lim_{n\to\infty}
    n\int_X |(I-P)\tilde\alpha|^2\,d\mu.
\end{align*}
Here each quadratic form is extended by $+\infty$ outside its domain.
With this convention, the equality also holds on all of
$\L^2(m;\mathcal H)$.
\end{thm}

\begin{proof}
Fix $\alpha\in\D(\breve{\mathcal E})$. We show that
$\breve{\mathcal E}^{\,n(I-P)\mu}(\alpha)$ increases to
$\breve{\mathcal E}^{\widetilde{\mathcal K}}(\alpha)$ as $n\to\infty$.
Since $\mu$ is $\widetilde{\mathcal K}$-detecting, the integral $\int_X |(I-P)\tilde\alpha|^2\,d\mu$ vanishes
if and only if $\tilde\alpha(x)\in\widetilde{\mathcal K}_x$ q.e.,
that is, if and only if
$\alpha\in\D(\breve{\mathcal E}^{\widetilde{\mathcal K}})$.
In this case, every term equals
$\breve{\mathcal E}(\alpha)
=\breve{\mathcal E}^{\widetilde{\mathcal K}}(\alpha)$.
Otherwise, the integral is strictly positive and the limit is $+\infty$.
Finally, if $\alpha\notin\D(\breve{\mathcal E})$, both the restriction
form and every perturbed form take the value $+\infty$ at $\alpha$.
\end{proof}

\subsection{Jump perturbation}\label{subsec:jump-perturbation}
In this subsection, we fix a symmetric $\sigma$-finite Borel
measure $j$ on $X\times X$ satisfying
\[
    j(\operatorname{diag}_X)=0,
    \qquad
    j(N\times X)=0
\]
for every exceptional Borel set $N\subset X$.
For $f\in\D(\mathcal E)$, define
\[
\mathcal J(f)
:=
\frac12\int_{X\times X}
|\tilde f(x)-\tilde f(y)|^2\,j(dx,dy).
\]
The jump-perturbed form $\mathcal E^j$
is defined by
\[
\D(\mathcal E^j)
:=
\left\{
f\in\D(\mathcal E)
\;\middle|\;
\mathcal J(f)<\infty
\right\},
\]
and
\[
\mathcal E^j(f,g)
:=
\mathcal E(f,g)
+
\mathcal J(f,g)
\qquad
f,g\in\D(\mathcal E^j).
\]
We write $\mathcal J(f,g)$ for
the polarization of $\mathcal J$.
\begin{dfn}\label{dfn:connection-rule}
A family $T=(T(x,y))_{(x,y)\in X\times X}$ of bounded
linear maps
\[
T(x,y)\colon
(\mathcal H^{\breve{\mathcal E}}_x)^*
\longrightarrow
(\mathcal H^{\breve{\mathcal E}}_y)^*
\]
is called a connection rule with respect to $j$ if:
\begin{enumerate}
    \item
    $\|T(x,y)\|_{\mathrm{op}}\leq1$
    for $j$-a.e.\ $(x,y)$.
    \item
    $T(y,x)=T(x,y)^\dagger$
    for $j$-a.e.\ $(x,y)$.
    \item
    For every $k,l\geq1$, the function
    \[
    (x,y)\longmapsto
    \langle\xi_k(x),T(x,y)^*\xi_l(y)\rangle
    \]
    is $j$-measurable, where $(\xi_k)_{k\geq1}$
    is a quasi-fundamental family representing the
    quasi-structure of $\mathcal H^{\breve{\mathcal E}}$.
\end{enumerate}
Here we regard the dual map $T(x,y)^*$ as a map from
$\mathcal H^{\breve{\mathcal E}}_y$ to $\mathcal H^{\breve{\mathcal E}}_x$
via the canonical identifications with the biduals.
The measurability condition is independent of the choice
of the representative.
We identify connection rules that agree $j$-a.e.
\end{dfn}
\begin{dfn}
For $\alpha\in\D(\breve{\mathcal E})$, define
\[
\breve{\mathcal J}^{T}(\alpha)
:=
\frac12\int_{X\times X}
\left(
|\tilde\alpha(x)|^2+|\tilde\alpha(y)|^2
-2\langle\tilde\alpha(x),T(x,y)^*\tilde\alpha(y)\rangle
\right)\,j(dx,dy).
\]
The jump-perturbed form $\breve{\mathcal E}^{j,T}$
is defined by
\[
\D(\breve{\mathcal E}^{j,T})
:=
\left\{
\alpha\in\D(\breve{\mathcal E})
\;\middle|\;
\breve{\mathcal J}^{T}(\alpha)<\infty
\right\},
\]
and
\[
\breve{\mathcal E}^{j,T}(\alpha,\beta)
:=
\breve{\mathcal E}(\alpha,\beta)
+
\breve{\mathcal J}^{T}(\alpha,\beta)
\qquad
\alpha,\beta\in\D(\breve{\mathcal E}^{j,T}).
\]
We write $\breve{\mathcal J}^{T}(\alpha,\beta)$ for
the polarization of $\breve{\mathcal J}^{T}$.
\end{dfn}

Both $\mathcal E^j$ and $\breve{\mathcal E}^{j,T}$ are
closed, symmetric, and nonnegative.
Indeed, a Cauchy sequence for either perturbed form norm
converges in the original form norm and admits a subsequence
converging q.e.\ in the corresponding scalar space or
quasi-Hilbert bundle.
Since $j$ does not charge exceptional sets in either
coordinate, Fatou's lemma applied to the nonnegative jump
integrand gives membership of the limit in the perturbed
domain.
Applying the same argument to differences gives convergence
in the perturbed form norm.
Their domains, however, need not be dense.

\begin{prp}\label{prp:jump-perturbation-kato}
For every $\alpha\in\D(\breve{\mathcal E}^{j,T})$
and $f\in\D(\mathcal E^j)_+$ with
$f\leq|\alpha|$ $m$-a.e., we have
\[
|\alpha|\in\D(\mathcal E^j),
\qquad
\frac{f}{|\alpha|}\alpha
\in\D(\breve{\mathcal E}^{j,T}),
\]
and
\[
\mathcal J(|\alpha|,f)
\leq
\breve{\mathcal J}^{T}
\left(\alpha,\frac{f}{|\alpha|}\alpha\right),
\qquad
\mathcal E^j(|\alpha|,f)
\leq
\breve{\mathcal E}^{j,T}
\left(\alpha,\frac{f}{|\alpha|}\alpha\right).
\]
Here $\frac{f}{|\alpha|}\alpha$ is set equal to zero
on $\{|\alpha|=0\}$.
Consequently, if $\D(\mathcal E^j)$ and
$\D(\breve{\mathcal E}^{j,T})$ are dense in
$\L^2(m)$ and $\L^2(m;\mathcal H)$, respectively, and
$\mathcal E^j$ is quasi-regular, then
$\breve{\mathcal E}^{j,T}$ is a Kato form with respect
to $\mathcal E^j$.
\end{prp}

\begin{proof}
The original Kato property gives
$|\alpha|\in\D(\mathcal E)$ and
$\frac{f}{|\alpha|}\alpha\in\D(\breve{\mathcal E})$.
The contractivity of $T(x,y)^*$ yields
\[
\mathcal J(|\alpha|)
\leq
\breve{\mathcal J}^{T}(\alpha),
\qquad
\breve{\mathcal J}^{T}
\left(\frac{f}{|\alpha|}\alpha\right)
\leq
\mathcal J(f)+\breve{\mathcal J}^{T}(\alpha),
\]
which proves the domain assertions.
The same pointwise calculation gives
\[
\mathcal J(|\alpha|,f)
\leq
\breve{\mathcal J}^{T}
\left(\alpha,\frac{f}{|\alpha|}\alpha\right).
\]
Adding the original Kato inequality proves the inequality
for the perturbed forms.
The Markov property of $\mathcal E^j$ follows from that
of $\mathcal E$ and the contraction property of the
scalar jump term.
The final assertion follows from Theorem~\ref{thm:quasi-realization-existence}.
\end{proof}

\begin{ex}\label{ex:jump-perturbation-equivalent-norms}
Suppose there exist constants $a,b\geq0$ such that
\[
\begin{aligned}
\mathcal J(f)
&\leq
a\,\mathcal E(f)+b\,\|f\|_2^2,
&& f\in\D(\mathcal E),\\
\breve{\mathcal J}^{T}(\alpha)
&\leq
a\,\breve{\mathcal E}(\alpha)+b\,\|\alpha\|_2^2,
&& \alpha\in\D(\breve{\mathcal E}).
\end{aligned}
\]
Then both perturbed domains coincide with the original
domains, and the corresponding form norms are equivalent.
In particular, $\mathcal E^j$ is quasi-regular and has
the same nests as $\mathcal E$.
Moreover, $(\mathcal H^{\breve{\mathcal E}},J)$ is also a
quasi-realization of $\breve{\mathcal E}^{j,T}$
with respect to $\mathcal E^j$.
\end{ex}

\clearpage
\section{On Riemannian manifolds}\label{sec:riemannian-bundle-forms}
\subsection{Connection energy with Neumann boundary condition}\label{subsec:neumann-connection-energy}
Let $M$ be a smooth Riemannian manifold with smooth boundary
$\partial M$, and let $\mathrm{vol}$ denote its Riemannian
volume measure.
Let $\mathcal E$ be the Neumann Dirichlet form on
$\L^2(\mathrm{vol})$.
We denote by $\nu$ the outward unit normal along $\partial M$
and by $\sigma$ the induced Riemannian measure on $\partial M$.
Smooth functions and sections are understood to be smooth
up to the boundary.
Let $E\to M$ be a smooth metric vector bundle of rank
$r\in\mathbb N$.
We fix a metric connection
\[
D:\Gamma^\infty(E)\longrightarrow\Gamma^\infty(T^*M\otimes E).
\]
We regard $D$ as a densely defined operator from
$\L^2(\mathrm{vol};E)$ to $\L^2(\mathrm{vol};T^*M\otimes E)$
with domain
\[
\D(D)
:=
\left\{
\alpha\in\Gamma^\infty(E)
\;\middle|\;
\alpha\in\L^2(\mathrm{vol};E),\
D\alpha\in\L^2(\mathrm{vol};T^*M\otimes E)
\right\}.
\]
Let $D^\dagger$ denote its formal adjoint with respect to
$\mathrm{vol}$ and the bundle metrics.
The integration-by-parts formula reads
\[
\int_M\langle D\alpha,\beta\rangle\,d\mathrm{vol}
=
\int_M\langle\alpha,D^\dagger\beta\rangle\,d\mathrm{vol}
+
\int_{\partial M}\langle\alpha,\beta(\nu)\rangle\,d\sigma,
\]
for $\alpha\in\Gamma^\infty(E)$ and
$\beta\in\Gamma^\infty_c(T^*M\otimes E)$.
Here compact supports are allowed to meet $\partial M$.
By testing this identity with $\beta$ supported in the interior
of $M$, we see that $D$ is closable.
We denote its closure again by $D$ and set
\[
W_C^{1,2}(E):=\D(D),
\qquad
\|\alpha\|_{W_C^{1,2}(E)}^2
:=\|\alpha\|_2^2+\|D\alpha\|_2^2.
\]
We define
\[
{\mathcal E}^C(\alpha,\beta)
:=
\frac12\int_M\langle D\alpha,D\beta\rangle\,d\mathrm{vol},
\qquad
\alpha,\beta\in
\D({\mathcal E}^C):=W_C^{1,2}(E).
\]
Then ${\mathcal E}^C$ is a densely defined, symmetric,
nonnegative closed form on $\L^2(\mathrm{vol};E)$.
We call ${\mathcal E}^C$ the connection energy with
Neumann boundary conditions.
\begin{prp}
${\mathcal E}^C$ is a Kato form on $(M,\text{vol},\mathcal E)$.
\end{prp}

\begin{proof}
We verify the Kato property.
First, let $\alpha\in\Gamma^\infty(E)\cap\D({\mathcal E}^C)$
and $f\in C^\infty\cap\D(\mathcal E)$ satisfy
$0\le f\le|\alpha|$.
For $\varepsilon>0$, set
$r_\varepsilon:=\sqrt{|\alpha|^2+\varepsilon^2}$.
Since $D$ is a metric connection,
\[
dr_\varepsilon
=
\frac{\langle D\alpha,\alpha\rangle}{r_\varepsilon},
\qquad
|dr_\varepsilon|\le |D\alpha|,
\]
and by the Leibniz rule,
\begin{align*}
D\left(\frac{f}{r_\varepsilon}\alpha\right)
&=
\frac{1}{r_\varepsilon}\,df\otimes\alpha
+
\frac{f}{r_\varepsilon}
\left(
D\alpha
-
\frac{1}{r_\varepsilon}\,dr_\varepsilon\otimes\alpha
\right).
\end{align*}
Since $0\le r_\varepsilon-\varepsilon\le|\alpha|$,
the above estimate implies that
$r_\varepsilon-\varepsilon\in\D(\mathcal E)$.
Moreover, the formula for $dr_\varepsilon$ and dominated
convergence show that $(r_\varepsilon-\varepsilon)_{\varepsilon>0}$
is Cauchy in the $\mathcal E_1$-norm as $\varepsilon\downarrow0$.
Since $r_\varepsilon-\varepsilon\to|\alpha|$
in $\L^2(\mathrm{vol})$, the closedness of $\mathcal E$ yields
$|\alpha|\in\D(\mathcal E)$.
Similarly,
$\left|\frac{f}{r_\varepsilon}\alpha\right|\le f$,
and the above formula gives
\[
\left|D\left(\frac{f}{r_\varepsilon}\alpha\right)\right|
\le |df|+2|D\alpha|.
\]
Thus
$\frac{f}{r_\varepsilon}\alpha\in\D({\mathcal E}^C)$.
The same formula and dominated convergence show that this family
is Cauchy in the ${\mathcal E}^C_1$-norm as
$\varepsilon\downarrow0$.
Since it converges to $\frac{f}{|\alpha|}\alpha$
in $\L^2(\mathrm{vol};E)$, the closedness of
${\mathcal E}^C$ yields
$\frac{f}{|\alpha|}\alpha\in\D({\mathcal E}^C)$ and
\begin{align*}
{\mathcal E}^C\left(\alpha,\frac{f}{|\alpha|}\alpha\right)
-\mathcal E(|\alpha|,f)
&=
\lim_{\varepsilon\downarrow0}
\left[
{\mathcal E}^C\left(
\alpha,\frac{f}{r_\varepsilon}\alpha
\right)
-\mathcal E(r_\varepsilon-\varepsilon,f)
\right]\\
&=
\lim_{\varepsilon\downarrow0}
\frac12\int_M\frac{f}{r_\varepsilon}
\left(
|D\alpha|^2-|dr_\varepsilon|^2
\right)\,d\mathrm{vol}
\ge0.
\end{align*}
Here $\frac{f}{|\alpha|}\alpha$ is understood to be zero on
$\{|\alpha|=0\}$.
Now let $\alpha\in\D({\mathcal E}^C)$ and
$f\in\D(\mathcal E)\cap\L^\infty(\mathrm{vol})$ satisfy
$0\le f\le|\alpha|$.
Choose $\alpha_n\in\Gamma^\infty(E)\cap\D({\mathcal E}^C)$
and a uniformly bounded sequence
$f_n\in C^\infty\cap\D(\mathcal E)$ such that
$
\alpha_n\to\alpha\;\text{in }{\mathcal E}^C_1,
\;
f_n\to f\;\text{in }\mathcal E_1.
$
Passing to a subsequence, we may also assume convergence
$\mathrm{vol}$-a.e.
Fix $\varepsilon>0$ and set
$r_{\varepsilon,n}:=\sqrt{|\alpha_n|^2+\varepsilon^2}$
and
$r_\varepsilon:=\sqrt{|\alpha|^2+\varepsilon^2}$.
The differential identities above apply to $\alpha_n$ and $f_n$
without requiring $f_n\le|\alpha_n|$.
Since $r_{\varepsilon,n}\ge\varepsilon$,
$|\alpha_n|/r_{\varepsilon,n}\le1$, and $(f_n)$ is uniformly
bounded, the form-norm convergence and dominated convergence
imply that these identities pass to the limit in the
corresponding $\L^2$ spaces.
By the closedness of $d$ and $D$, we obtain
$
r_\varepsilon-\varepsilon\in\D(\mathcal E),
\;
\frac{f}{r_\varepsilon}\alpha\in\D({\mathcal E}^C),
$
and
\begin{align*}
d(r_\varepsilon-\varepsilon)
&=
\frac{\langle D\alpha,\alpha\rangle}{r_\varepsilon},\\
D\left(\frac{f}{r_\varepsilon}\alpha\right)
&=
\frac{1}{r_\varepsilon}\,df\otimes\alpha
+
\frac{f}{r_\varepsilon}
\left(
D\alpha
-
\frac{1}{r_\varepsilon}\,dr_\varepsilon\otimes\alpha
\right).
\end{align*}
Since $0\le f\le|\alpha|$, the preceding estimates and
dominated convergence argument apply as $\varepsilon\downarrow0$.
Thus $|\alpha|\in\D(\mathcal E)$ and
$\frac{f}{|\alpha|}\alpha\in\D({\mathcal E}^C)$,
and the same calculation yields
\[
\mathcal E(|\alpha|,f)
\le
{\mathcal E}^C\left(
\alpha,\frac{f}{|\alpha|}\alpha
\right).
\]
For general $f\in\D(\mathcal E)_+$ with $f\le|\alpha|$,
the same conclusion follows by applying the above argument
to $f\wedge n$ and letting $n\to\infty$, using truncation,
dominated convergence, and the closedness of
${\mathcal E}^C$.
\end{proof}

\begin{rem}
The inequality $|d|\alpha||\le|D\alpha|$ obtained in the proof is called Kato's inequality.
\end{rem}
We equip $E$ with the natural quasi-Hilbert bundle structure
induced by its smooth bundle structure and the
capacity.
\begin{prp}
$(E,\operatorname{Id}_E)$ is a quasi-realization of
${\mathcal E}^C$, where $\operatorname{Id}_E$ is regarded
as a map to the underlying almost Hilbert bundle.
\end{prp}
\begin{proof}
    Smooth local trivializations and the quasi-continuity
of Sobolev functions show that every
$\alpha\in\D({\mathcal E}^C)$ admits a quasi-continuous
representative as a section of $E$.
Moreover, by adjoining compactly supported multiples of local
frame sections to a countable smooth
${\mathcal E}^C_1$-dense family, we obtain a
${\mathcal E}^C_1$-dense family whose smooth representatives
form a quasi-fundamental family of $E$.
Thus the defining properties of a quasi-realization hold.
\end{proof}
\subsection{Connection energy with mixed boundary conditions}\label{subsec:mixed-connection-energy}
We equip $\otimes^kT^*M$ with the metric connection induced by
the Levi--Civita connection, again denoted by $D$.
Along $\partial M$, we identify $T^*\partial M$ with the
subspace of $T^*M$ consisting of covectors that annihilate
$\nu$.
For every integer $k\geq0$, let $H_C^{1,2}(\otimes^kT^*M)$
denote the closure in $W_C^{1,2}(\otimes^kT^*M)$ of
\[
\left\{
\alpha\in\Gamma^\infty(\otimes^kT^*M)
\cap W_C^{1,2}(\otimes^kT^*M)
\;\middle|\;
\alpha_x\in\otimes^kT_x^*\partial M
\text{ for every }x\in\partial M
\right\}.
\]
Define $T_{\mathrm{tan}}M\subset TM$ fiberwise by
\[
(T_{\mathrm{tan}}M)_x
:=
\begin{cases}
T_xM,&x\in M^\circ,\\
T_x\partial M,&x\in\partial M.
\end{cases}
\]
The family $T_{\mathrm{tan}}M$ is a quasi-Hilbert subbundle
of $TM$ and is equipped with the induced quasi-Hilbert
bundle structure.
Its dual is denoted by $T_{\mathrm{tan}}^*M$ and is identified
with the corresponding family of subspaces of $T^*M$ by
the Riemannian metric.
With domain $H_C^{1,2}(\otimes^kT^*M)$, the connection
energy $\mathcal E^C$ coincides with the restriction of the
Neumann connection energy to the quasi-Hilbert subbundle
$\otimes^kT_{\mathrm{tan}}^*M$ of $\otimes^kT^*M$.
In particular, it is a Kato form with respect to $\mathcal E$.
Moreover, the pair
$
\left(\otimes^kT_{\mathrm{tan}}^*M,\iota\right)
$
is a quasi-realization of this form,
where $\iota$ is the natural inclusion into the underlying
almost Hilbert bundle of $\otimes^kT^*M$.
The boundary measure $\sigma$ is a $T_{\mathrm{tan}}M$-detecting measure
for the Neumann connection energy on $TM$. More generally, $\sigma$ is
$\otimes^kT_{\mathrm{tan}}^*M$-detecting for the Neumann
connection energy on $\otimes^kT^*M$.
Let $P(x)$ denote the orthogonal projection onto
$\otimes^k T_x^*\partial M$ at $x\in\partial M$.
By Theorem~\ref{thm:restriction-as-supremum}, the mixed connection energy,
extended by $+\infty$ outside $H_C^{1,2}(\otimes^kT^*M)$, satisfies
\[
    \mathcal E^C(\alpha)
    =\frac12\int_M |D\alpha|^2\,d\mathrm{vol}
        +
    \lim_{n\to\infty}\left(
        n\int_{\partial M}|(I-P)\tilde\alpha|^2\,d\sigma
    \right).
\]
Replacing tensor powers by exterior powers defines
$H_C^{1,2}(\Lambda^kT^*M)$ for $0\leq k\leq\dim M$.
The same assertions hold, with quasi-realization
$(\Lambda^kT_{\mathrm{tan}}^*M,\iota)$.

\subsection{Hodge energy with mixed boundary conditions}\label{subsec:hodge-energy-mixed-boundary}
For each $0\le k\le\dim M$, we regard the exterior derivative
$d$ as a densely defined operator from
$\L^2(\mathrm{vol};\Lambda^kT^*M)$ to
$\L^2(\mathrm{vol};\Lambda^{k+1}T^*M)$
with initial domain
\[
\D(d)\coloneqq
\left\{
\alpha\in\Gamma^\infty(\Lambda^kT^*M)
\;\middle|\;
\alpha\in\L^2(\mathrm{vol};\Lambda^kT^*M),\
d\alpha\in\L^2(\mathrm{vol};\Lambda^{k+1}T^*M)
\right\}.
\]
Let $d^\dagger$ denote the formal adjoint of $d$ with respect
to $\mathrm{vol}$ and the induced bundle metrics.
The integration-by-parts formula reads
\[
\int_M\langle d\alpha,\beta\rangle\,d\mathrm{vol}
=
\int_M\langle\alpha,d^\dagger\beta\rangle\,d\mathrm{vol}
+
\int_{\partial M}\langle\alpha,\nu(\beta)\rangle\,d\sigma,
\]
for $\alpha\in\Gamma^\infty(\Lambda^kT^*M)$ and
$\beta\in\Gamma_c^\infty(\Lambda^{k+1}T^*M)$, where
$\nu(\beta)$ denotes contraction of $\beta$ with $\nu$.
By testing against $\beta$ supported in $M^\circ$, we see
that $d$ is closable.
We denote its closure again by $d$ and its Hilbert space
adjoint by $d^\dagger$.
We set
\[
W_H^{1,2}(\Lambda^kT^*M)
:=\D(d)\cap\D(d^\dagger),
\]
equipped with the norm
\[
\|\alpha\|_{W_H^{1,2}(\Lambda^kT^*M)}^2
:=
\|\alpha\|_2^2+\|d\alpha\|_2^2+\|d^\dagger\alpha\|_2^2.
\]
We define the Hodge energy on
$\L^2(\mathrm{vol};\Lambda^kT^*M)$ with domain
\[
\D(\mathcal E^H):=W_H^{1,2}(\Lambda^kT^*M)
\]
by
\[
\mathcal E^H(\alpha,\beta)
:=
\frac12\int_M
\left(
\langle d\alpha,d\beta\rangle
+
\langle d^\dagger \alpha,d^\dagger\beta\rangle
\right)\,d\mathrm{vol},
\qquad
\alpha,\beta\in\D(\mathcal E^H),
\]
with the convention that $d^\dagger=0$ on $0$-forms.
Then $\mathcal E^H$ is a bundle form on $\Lambda^k T^*M$.

\begin{prp}
The pair
$
\left(
\Lambda^kT_{\mathrm{tan}}^*M,
\iota
\right)
$
is a quasi-realization of $\mathcal E^H$, where $\iota$ denotes the natural inclusion into the underlying
almost Hilbert bundle of $\Lambda^kT^*M$.
\end{prp}
\begin{proof}
    For $k\ge1$, the adjoint identity and the integration-by-parts
formula give
\[
\int_{\partial M}\langle\beta,\nu(\alpha)\rangle\,d\sigma
=
\int_M
\left(
\langle d\beta,\alpha\rangle
-
\langle\beta,d^\dagger\alpha\rangle
\right)\,d\mathrm{vol}
=0
\]
for every smooth $\alpha\in\D(d^\dagger)$ and
$\beta\in\Gamma_c^\infty(\Lambda^{k-1}T^*M)$.
Since the boundary values of $\beta$ are arbitrary,
$\nu(\alpha)=0$ on $\partial M$.
For general elements of $\D(d^\dagger)$, the same identity
expresses vanishing of the weak normal trace.
Local boundary regularity shows that
every $\alpha\in\D(\mathcal E^H)$ belongs to
$W^{1,2}_{\mathrm{loc}}$ up to the boundary.
Consequently, it admits a quasi-continuous representative
satisfying $\nu(\tilde\alpha)=0$ q.e. on $\partial M$.
Finally, elliptic regularization provides a smooth
$\mathcal E^H_1$-dense family.
Adjoining compactly supported smooth local generators
with vanishing normal component gives a quasi-fundamental
family of $\Lambda^kT_{\mathrm{tan}}^*M$.
Thus the defining properties of a quasi-realization hold.
\end{proof}
We equip $\Lambda^kT^*M$ with the connection induced by
the Levi--Civita connection.
Let $\mathscr R_k$ denote the curvature endomorphism
in the Weitzenb\"ock formula for $k$-forms.
Along $\partial M$, define
\[
S:T\partial M\longrightarrow T\partial M,
\qquad
S(w):=D_w\nu,
\]
and let $S_k$ be the induced symmectric endomorphism of
$\Lambda^kT^*\partial M$, given by
\[
(S_k\alpha)(w_1,\ldots,w_k)
:=
\sum_{j=1}^k
\alpha(w_1,\ldots,S(w_j),\ldots,w_k).
\]
For
$\alpha,\beta\in\Gamma_c^\infty(\Lambda^kT^*M)$
satisfying $\nu(\alpha)=\nu(\beta)=0$ on $\partial M$,
the integrated Weitzenb\"ock formula yields
\[
\begin{aligned}
\mathcal E^H(\alpha,\beta)
={}&
\mathcal E^C(\alpha,\beta)
+
\frac12\int_M
\langle\mathscr R_k\alpha,\beta\rangle\,d\mathrm{vol}+
\frac12\int_{\partial M}
\langle S_k\alpha,\beta\rangle\,d\sigma.
\end{aligned}
\]
Thus, on this class of forms, the Hodge energy is the
potential perturbation of the connection energy by the
symmetric endomorphism-valued measure
$
\frac12\mathscr R_k\,\mathrm{vol}+\frac12 S_k\,\sigma.
$
\begin{ex}
For $k=1$, $\mathscr R_1$ and $S_1$ are respectively the Ricci endomorphism and the endomorphism associated with the second fundamental form.
Assume either that $k=1$ and $M$ is complete with Ricci curvature bounded from below
and convex boundary (possibly empty),
or that $M$ is compact, with arbitrary $k$.
Then the two form domains coincide,
\[
W_H^{1,2}(\Lambda^kT^*M)
=
H_C^{1,2}(\Lambda^kT^*M),
\]
and the preceding identity extends to the whole common
domain, with the boundary integral evaluated using
quasi-continuous representatives.
Thus the potential perturbation identity holds
as an equality of closed forms.
\end{ex}

\clearpage
\chapter{Multiplicative Stochastic Transport}\label{chap:multiplicative-stochastic-transports}
\section*{Overview}
This chapter constructs stochastic transport from the Kato forms of Chapter~\ref{chap:bundle-forms-quasi-realization}. Theorems~\ref{thm:multiplicative-transport-existence} and~\ref{thm:multiplicative-transport-uniqueness} show that every Kato form admits a unique associated contractive quasi-right transport on its quasi-realization, up to indistinguishability for q.e.\ starting point. Conversely, a symmetric contractive quasi-right transport determines a Kato form (Theorem~\ref{thm:transport-kato-form}). Thus, Kato forms and symmetric contractive quasi-right transports are in one-to-one correspondence. Section~\ref{sec:transport-operations} begins with the resolvent formula for vector-valued measures (Theorem~\ref{thm:vector-revuz-resolvent}) and then realizes the operations of Section~\ref{sec:bundle-form-operations} as damped, restriction, and jump transports. To allow for applications to nonsmooth geometry, we construct damped transport independently of Hsu's approach to multiplicative functionals. Section~\ref{sec:riemannian-transports} applies these constructions on Riemannian manifolds to stochastic parallel transport with Neumann or mixed boundary conditions and to Weitzenb\"ock Ricci-damped transport with mixed boundary conditions.

\section*{Prior work}
It\^o \cite{Ito1963} introduced parallel transport along Brownian paths on a Riemannian manifold. For the construction by horizontal lifting, see \cite{Elworthy1982,Hsu2002}. Meyer \cite{Meyer1962} established the correspondence between scalar multiplicative functionals and subordinate semigroups. On a Hermitian bundle, Must\u{a}\c{t}ea \cite{Mustatea2022} constructed stochastic parallel transport using the connection heat kernel. Our resolvent formula for vector-valued measures builds on the scalar Revuz correspondence \cite{Revuz1970}. G\"uneysu \cite{Guneysu2012} proved a Feynman--Kac formula for Schr\"odinger operators on Riemannian vector bundles, and G\"uneysu, Milatovic and Truc \cite{GuneysuMilatovicTruc2014} constructed transport on weighted graphs from unitary maps along edges. For reflecting Brownian motion, Ikeda and Watanabe \cite{IkedaWatanabe1978} used excursion theory to represent solutions of the heat equation on forms. Hsu \cite{Hsu2002Boundary} constructed the corresponding multiplicative functional using an approximation suggested by Airault \cite{Airault1976}. Arnaudon and Li \cite{ArnaudonLi2017} described how the damped transport loses its normal component at excursion endpoints.

\section{Multiplicative stochastic transport}\label{sec:multiplicative-transports}
\subsection{Stochastic lift}
Let $\widetilde{\mathcal H}$ be a quasi-Hilbert bundle over $(X,m, \mathcal E)$.
\begin{dfn}
An $\widetilde{\mathcal H}^*$-valued stochastic process $(\mathsf V_t)_{t\geq0}$
is called a lift of $(\mathsf X_t)_{t\geq0}$ to $\widetilde{\mathcal H}$ if
\[
    \mathsf V_t(\omega)\in \widetilde{\mathcal H}^*_{\mathsf X_t(\omega)}
    \qquad
    \text{for every }(t,\omega)\in[0,\infty)\times\Omega.
\]
\end{dfn}
\begin{dfn}
Let $(\xi_k)_{k\geq1}$ be a quasi-fundamental family of $\widetilde{\mathcal H}$,
and let $(F_n)_{n\geq1}$ be a nest adapted to $(\xi_k)_{k\geq1}$.
A lift $(\mathsf V_t)_{t\geq0}$ of $(\mathsf X_t)_{t\geq0}$
to $\widetilde{\mathcal H}$ is called quasi-right-continuous with respect to
$(\xi_k,F_n)$ if, for q.e.\ $x\in X$ and every $n\geq1$,
the map
\[
    [0,\tau_n)\ni t\mapsto\mathsf V_t\in \widetilde{\mathcal H}^*|_{F_n},
    \quad
    \tau_n
    :=
    \inf\{t\geq0\mid\mathsf X_t\notin F_n\}.
\]
is right-continuous $\mathbb P_x$-a.s.,
where $\widetilde{\mathcal H}^*|_{F_n}$ carries the topology induced by
$(\xi_k^*)_{k\geq1}$.
\end{dfn}

\begin{prp}\label{prp:quasi-right-lift-invariance}
Let $(\xi_k)_{k\geq1}$ and $(\xi'_k)_{k\geq1}$ be two equivalent
quasi-fundamental families of $\widetilde{\mathcal H}$, and let $(F_n)_{n\geq1}$
and $(F'_n)_{n\geq1}$ be nests adapted to these families,
respectively.
Then a lift $(\mathsf V_t)_{t\geq0}$ is quasi-right-continuous
with respect to $(\xi_k,F_n)$ if and only if it is
quasi-right-continuous with respect to $(\xi'_k,F'_n)$.
\end{prp}
\begin{proof}
Assume that $(\mathsf V_t)_{t\geq0}$ is quasi-right-continuous
with respect to $(\xi_k,F_n)$.
By Proposition~\ref{prp:compatible-nest-topologies}, there exists a common
refinement $(G_n)_{n\geq1}$ such that the topologies on
$\widetilde{\mathcal H}^*|_{G_n}$ induced by $(\xi_k^*)_{k\geq1}$ and
$((\xi'_k)^*)_{k\geq1}$ coincide.
Denote by $\tau_n$, $\tau'_n$, and $\sigma_n$ the exit times
from $F_n$, $F'_n$, and $G_n$, respectively.
Since $(G_n)_{n\geq1}$ is a nest,
\[
    \sigma_n\uparrow\zeta
    \qquad
    \mathbb P_x\text{-a.s.\ for q.e.\ }x\in X,
\]
where $\zeta$ is the lifetime.
Fix such an $x$ and work on a common event of probability one
on which this convergence and the assumed right-continuity
hold for every $n$.
Since $G_n\subset F_n$ and $\sigma_n\leq\tau_n$,
the restriction of $\mathsf V$ to $[0,\sigma_n)$ is
right-continuous as an $\widetilde{\mathcal H}^*|_{G_n}$-valued map.
By the choice of $(G_n)_{n\geq1}$, this holds also for the
topology induced by $((\xi'_k)^*)_{k\geq1}$.
Fix $m\geq1$ and $t<\tau'_m$.
Choose $n\geq m$ such that $t<\sigma_n$.
On a sufficiently small right neighborhood of $t$,
$\mathsf V$ is right-continuous as an $\widetilde{\mathcal H}^*|_{F'_n}$-valued
map and takes values in $\widetilde{\mathcal H}^*|_{F'_m}$.
The topology on $\widetilde{\mathcal H}^*|_{F'_m}$ is the subspace topology
inherited from $\widetilde{\mathcal H}^*|_{F'_n}$, so $\mathsf V$ is
right-continuous at $t$ as an $\widetilde{\mathcal H}^*|_{F'_m}$-valued map.
Thus $(\mathsf V_t)_{t\geq0}$ is quasi-right-continuous
with respect to $(\xi'_k,F'_n)$.
The converse follows by symmetry.
\end{proof}
We therefore simply call such a lift quasi-right-continuous.

\subsection{Multiplicative stochastic transport}\label{subsec:multiplicative-transports}
Let $\widetilde{\mathcal H}$ be a quasi-Hilbert bundle over $(X,m,\mathcal{E})$.
\begin{dfn}\label{dfn:stochastic-transport}
Let $(\mathsf U_t)_{t\geq0}$ be a family of bounded linear maps
\[
    \mathsf U_t(\omega)\colon
    \widetilde{\mathcal H}^*_{\mathsf X_0(\omega)}
    \longrightarrow \widetilde{\mathcal H}^*_{\mathsf X_t(\omega)}.
\]
We call $(\mathsf U_t)_{t\geq0}$ a stochastic transport
of $(\mathsf X_t)_{t\geq0}$ on $\widetilde{\mathcal H}$
if the following conditions hold\;:
\begin{enumerate}
    \item $\mathsf U_0(\omega)=I$ for every $\omega\in\Omega$.

    \item For every measurable section $\xi$ of $\widetilde{\mathcal H}^*$,
    the lift
    $(\mathsf U_t\xi(\mathsf X_0))_{t\geq0}$
    is $(\mathcal F_{t+}^{\mathsf X})_{t\geq0}$-progressively
    measurable.
    \item For every $T>0$, there exists a constant $C_T<\infty$
    such that
    \[
        \mathbb E_x\left[
            \sup_{0\leq t\leq T}\|\mathsf U_t\|_{\mathrm{op}}
        \right]
        \leq C_T
        \qquad\text{for every }x\in X.
    \]
\end{enumerate}
A stochastic transport $(\mathsf U_t)_{t\geq0}$
is called contractive if condition~3 holds with $C_T=1$
for every $T>0$.
Equivalently, for every $x\in X$,
$\|\mathsf U_t\|_{\mathrm{op}}\leq1$ simultaneously for all $t\geq0$,
$\mathbb P_x$-a.s.
\end{dfn}
Let $(\mathsf U_t)_{t\geq0}$ be a stochastic transport
of $(\mathsf X_t)_{t\geq0}$ on $\widetilde{\mathcal H}$.
For every bounded measurable section $\alpha$ of $\widetilde{\mathcal H}$
and every $t\geq0$, we set
\[
    \breve{p}_t^{\mathsf U}\alpha(x)
    \coloneqq
    \mathbb E_x\left[
        \mathsf U_t^*\alpha(\mathsf X_t)
    \right],
    \qquad x\in X,
\]
where the expectation is understood as a Bochner integral
in $\widetilde{\mathcal H}_x$.
Conditions~2 and~3 ensure that $\breve{p}_t^{\mathsf U}\alpha$
is a bounded measurable section, with
\[
    |\breve{p}_t^{\mathsf U}\alpha(x)|
    \leq C_T\|\alpha\|_\infty,
    \qquad 0\leq t\leq T,\quad x\in X.
\]
\begin{dfn}
Let $(\mathsf U_t)_{t\geq0}$ be a stochastic transport.
\begin{enumerate}
    \item We call $(\mathsf U_t)_{t\geq0}$ multiplicative if,
    for q.e.\ $x\in X$ and every deterministic $t,h\geq0$,
    \[
        \mathsf U_{t+h}
        =
        (\mathsf U_h\circ\theta_t)\mathsf U_t
        \qquad\mathbb P_x\text{-a.s.}
    \]

    \item We call $(\mathsf U_t)_{t\geq0}$ strongly multiplicative if,
    for q.e.\ $x\in X$, every bounded
    $(\mathcal F_{t+}^{\mathsf X})_{t\geq0}$-stopping time
    $\tau$, and every $h\geq0$,
    \[
        \mathsf U_{\tau+h}
        =
        (\mathsf U_h\circ\theta_\tau)\mathsf U_\tau
        \qquad\mathbb P_x\text{-a.s.}
    \]
\end{enumerate}
\end{dfn}
For $0\leq s\leq t$, we set
\[
    \mathsf U_{s,t}(\omega)
    :=
    \mathsf U_{t-s}(\theta_s\omega)
    \colon
    \widetilde{\mathcal H}_{\mathsf X_s(\omega)}^*
    \longrightarrow
    \widetilde{\mathcal H}_{\mathsf X_t(\omega)}^*.
\]
Multiplicativity is equivalent to
\[
    \mathsf U_{0,t}
    =
    \mathsf U_{s,t}\mathsf U_{0,s}
    \qquad\mathbb P_x\text{-a.s.}
\]
for q.e.\ $x\in X$ and every deterministic $0\leq s\leq t$.
Strong multiplicativity is equivalent to
\[
    \mathsf U_{0,\tau+h}
    =
    \mathsf U_{\tau,\tau+h}\mathsf U_{0,\tau}
    \qquad\mathbb P_x\text{-a.s.}
\]
for q.e.\ $x\in X$, every bounded
$(\mathcal F_{t+}^{\mathsf X})_{t\geq0}$-stopping time $\tau$,
and every $h\geq0$.
\begin{ex}
Multiplicativity does not imply strong multiplicativity.
Let $(\mathsf X_t)_{t\geq0}$ be reflecting Brownian motion on
$[0,\infty)$, equipped with Lebesgue measure, and consider the quasi-Hilbert bundle
\[
    \widetilde{\mathcal H}_x=
    \begin{cases}
        \mathbb R,&x>0,\\
        \{0\},&x=0,
    \end{cases}
\]
Identifying the fibers with subspaces of $\mathbb R$, set
\[
    \mathsf U_t
    :=
    \mathbf1_{\{\mathsf X_0>0,\ \mathsf X_t>0\}}.
\]
This defines a contractive multiplicative stochastic transport
(because for every $x\geq0$ and every deterministic $t>0$,
$\mathbb P_x(\mathsf X_t>0)=1$).
On the other hand, fix $x>0$ and let
\[
    \tau:=
    \inf\{t\geq0\mid \mathsf X_t=0\}\wedge1.
\]
The event $\{\tau<1\}$ has positive $\mathbb P_x$-probability.
On this event, $\mathsf U_\tau=0$, while the strong Markov
property gives, for every fixed $h>0$,
\[
    \mathsf U_{\tau+h}=1,
    \qquad
    (\mathsf U_h\circ\theta_\tau)\mathsf U_\tau=0
    \qquad\mathbb P_x\text{-a.s.}
\]
Thus strong multiplicativity fails at the bounded stopping time $\tau$.
\end{ex}
If $(\mathsf U_t)_{t\geq0}$ is multiplicative, the Markov property of
$(\mathsf X_t)_{t\geq0}$ implies
\[
    \breve{p}_0^{\mathsf U}=I,
    \qquad
    \breve{p}_{s+t}^{\mathsf U}
    =\breve{p}_s^{\mathsf U}\breve{p}_t^{\mathsf U},
    \qquad s,t\geq0.
\]
Thus $(\breve{p}_t^{\mathsf U})_{t\geq0}$ is a semigroup
on bounded measurable sections of $\widetilde{\mathcal H}$.
We call it the semigroup induced by
$(\mathsf U_t)_{t\geq0}$. 
Multiplicativity, the Markov property, and condition~3 yield
constants $M\geq1$ and $c\geq0$ such that
\[
    \mathbb E_x[\|\mathsf U_t\|_{\mathrm{op}}]\leq M e^{ct},
    \qquad t\geq0,\quad\text{for q.e.\ }x\in X.
\]
For every $\lambda>c$ and every bounded measurable section
$\alpha$ of $\widetilde{\mathcal H}$, we therefore define
\[
    \breve R_\lambda^{\mathsf U}\alpha(x)
    \coloneqq
    \int_0^\infty e^{-\lambda t}
        \breve p_t^{\mathsf U}\alpha(x)\,dt
    =
    \mathbb E_x\left[
        \int_0^\infty e^{-\lambda t}
            \mathsf U_t^*\alpha(\mathsf X_t)\,dt
    \right],
    \qquad\text{for q.e.\ }x\in X,
\]
where the integrals are Bochner integrals in
$\widetilde{\mathcal H}_x$.
We call $(\breve R_\lambda^{\mathsf U})_{\lambda>c}$
the resolvent induced by $(\mathsf U_t)_{t\geq0}$.
In particular, if $(\mathsf U_t)_{t\geq0}$ is contractive,
we may take $M=1$ and $c=0$; the above formula then defines
$\breve R_\lambda^{\mathsf U}\alpha(x)$ for every $\lambda>0$
and $x\in X$, with
$|\breve R_\lambda^{\mathsf U}\alpha(x)|
\leq\lambda^{-1}\|\alpha\|_\infty$.
\begin{dfn}
A multiplicative stochastic transport $(\mathsf U_t)_{t\geq0}$
is said to be symmetric if, for every $t\geq0$ and all bounded measurable
sections $\alpha,\beta$ of $\widetilde{\mathcal H}$ belonging to
$\mathrm L^1(m;\widetilde{\mathcal H})$,
\[
    \int_X
        \langle\breve p_t^{\mathsf U}\alpha(x),\beta(x)\rangle
        \,dm(x)
    =
    \int_X
        \langle\alpha(x),\breve p_t^{\mathsf U}\beta(x)\rangle
        \,dm(x).
\]
\end{dfn}
Equivalently, the same identity holds with
$\breve p_t^{\mathsf U}$ replaced by $\breve R_\lambda^{\mathsf U}$
for every sufficiently large $\lambda$ and all such $\alpha,\beta$.
\begin{dfn}\label{dfn:quasi-right-transport}
A multiplicative stochastic transport $(\mathsf U_t)_{t\geq0}$
of $(\mathsf X_t)_{t\geq0}$ on $\widetilde{\mathcal H}$
is called a quasi-right transport if, for every measurable
section $\xi$ of $\widetilde{\mathcal H}^*$, the lift 
$(\mathsf U_t\xi(\mathsf X_0))_{t\geq0}$
is quasi-right-continuous.
\end{dfn}

\begin{rem}
Equivalently, let $(\eta_k)_{k\geq1}$ be a quasi-fundamental
family of $\widetilde{\mathcal H}$.
Then a multiplicative stochastic transport $(\mathsf U_t)_{t\geq0}$
is a quasi-right transport if and only if, for every measurable section $\xi$ of $\widetilde{\mathcal H}^*$
and for q.e.\ $x\in X$, the scalar processes
\begin{align*}
    &t\longmapsto
    \langle\eta_k(\mathsf X_t),
        \mathsf U_t\xi(\mathsf X_0)\rangle,
    \qquad k\geq1,\\
    &t\longmapsto|\mathsf U_t\xi(\mathsf X_0)|
\end{align*}
are right-continuous $\mathbb P_x$-a.s.
This follows from the definition of the topologies on
$\widetilde{\mathcal H}^*|_{F_n}$ and the fact that the exit times of a nest
increase to the lifetime.
\end{rem}

\begin{rem}
If $(\mathsf U_t)_{t\geq0}$ is a quasi-right transport, then,
for q.e.\ $x\in X$, the paths $t\mapsto\mathsf U_t v$ are
quasi-right-continuous simultaneously for all
$v\in\widetilde{\mathcal H}_x^*$ on a common event of
$\mathbb P_x$-probability one.
Indeed, apply the definition to a countable fiberwise dense
family of measurable initial sections and extend by the almost sure
local boundedness of $\|\mathsf U_t\|_{\mathrm{op}}$
(Definition~\ref{dfn:stochastic-transport}, condition~3).
\end{rem}

\begin{prp}\label{prp:quasi-right-strong-multiplicativity}
Let $(\mathsf U_t)_{t\geq0}$ be a quasi-right transport.
Then $(\mathsf U_t)_{t\geq0}$ is strongly multiplicative.
\end{prp}
\begin{proof}\;
\begin{itemize}
\item
We first show that the short-time error between the original
transport and the transport restarted at a stopping time
tends to zero in $\mathrm L^1$.
Choose a properly exceptional Borel set $N$ outside which
deterministic-time multiplicativity and the conclusion of
the preceding remark hold.
By the strong Markov property, the latter also holds after
restarting at a bounded stopping time, with any initial vector
measurable at that time.
Fix $x\notin N$, a bounded
$(\mathcal F_{t+}^{\mathsf X})_{t\geq0}$-stopping time
$\tau\leq T$, an increment $h>0$, and
$v\in\widetilde{\mathcal H}_x^*$.
To approximate $\tau$ by finite-valued stopping times and
record the intervening time interval, set
\[
    \tau_n:=2^{-n}\lceil2^n\tau\rceil,
    \qquad
    \delta_n:=\tau_n-\tau.
\]
Then $\tau_n\downarrow\tau$ and $\delta_n\downarrow0$.
Both $\mathsf U_{\tau_n}v$ and
$(\mathsf U_{\delta_n}\circ\theta_\tau)\mathsf U_\tau v$
belong to the fiber over $\mathsf X_{\tau_n}$.
By quasi-right-continuity, the strong Markov property, and
$\mathsf U_0=I$, they converge to the same vector
$\mathsf U_\tau v$ in the bundle topology on a common adapted
nest. Hence their difference in fiber norm tends to zero.
To measure this error, put
\[
    \varepsilon_n:=
    \left|
        \mathsf U_{\tau_n}v
        -(\mathsf U_{\delta_n}\circ\theta_\tau)
            \mathsf U_\tau v
    \right|.
\]
Thus $\varepsilon_n\to0$ almost surely.
We use the zero extension after the lifetime throughout.
Moreover,
\[
\begin{aligned}
    \varepsilon_n
    &\leq
    |v|\sup_{0\leq r\leq T+1}\|\mathsf U_r\|_{\mathrm{op}}+
    |\mathsf U_\tau v|
    \left(
        \sup_{0\leq r\leq1}\|\mathsf U_r\|_{\mathrm{op}}
        \circ\theta_\tau
    \right).
\end{aligned}
\]
The right-hand side is integrable: condition~3 bounds the
expectation of its first term, and the strong Markov property
bounds that of its second term by
$C_1\mathbb E_x[|\mathsf U_\tau v|]$.
Dominated convergence therefore gives
\[
    \mathbb E_x[\varepsilon_n]\longrightarrow0.
\]

\item
We next obtain two identities using deterministic-time
multiplicativity.
Since $\tau_n$ takes finitely many deterministic values,
\[
    \mathsf U_{\tau_n+h}
    =
    (\mathsf U_h\circ\theta_{\tau_n})\mathsf U_{\tau_n}
    \qquad\mathbb P_x\text{-a.s.}
\]
Also, $\delta_n$ is $\mathcal F_\tau^{\mathsf X,x}$-measurable.
Conditional on this sigma-field, it is a fixed time parameter
for the shifted path. The strong Markov property and
deterministic-time multiplicativity under
$\mathbb P_{\mathsf X_\tau}$ consequently give
\[
    \mathsf U_{h+\delta_n}\circ\theta_\tau
    =
    (\mathsf U_h\circ\theta_{\tau_n})
    (\mathsf U_{\delta_n}\circ\theta_\tau)
    \qquad\mathbb P_x\text{-a.s.}
\]
Here we used $\tau+\delta_n=\tau_n$ and the composition rule
for the shifts.

\item
We now propagate the short-interval error through the
remaining time increment $h$.
Subtracting the preceding identities after applying them to
the relevant initial vectors yields
\[
\begin{aligned}
    &\mathsf U_{\tau_n+h}v
    -(\mathsf U_{h+\delta_n}\circ\theta_\tau)\mathsf U_\tau v=
    (\mathsf U_h\circ\theta_{\tau_n})
    \left[
        \mathsf U_{\tau_n}v
        -(\mathsf U_{\delta_n}\circ\theta_\tau)\mathsf U_\tau v
    \right].
\end{aligned}
\]
Progressive measurability makes $\varepsilon_n$
$\mathcal F_{\tau_n}^{\mathsf X,x}$-measurable.
The strong Markov property and condition~3 give
\[
    \mathbb E_x\!\left[
        \|\mathsf U_h\circ\theta_{\tau_n}\|_{\mathrm{op}}
        \,\middle|\,\mathcal F_{\tau_n}^{\mathsf X,x}
    \right]
    \leq C_h.
\]
Taking norms and conditional expectations therefore gives
\[
    \mathbb E_x\!\left[
        \left|
            \mathsf U_{\tau_n+h}v
            -(\mathsf U_{h+\delta_n}\circ\theta_\tau)
                \mathsf U_\tau v
        \right|
    \right]
    \leq C_h\mathbb E_x[\varepsilon_n]
    \longrightarrow0.
\]

\item
Finally, we pass to the limit by quasi-right-continuity.
Since $\tau_n+h\downarrow\tau+h$ and $h+\delta_n\downarrow h$,
the norm inside the last expectation converges almost surely to
\[
    \left|
        \mathsf U_{\tau+h}v
        -(\mathsf U_h\circ\theta_\tau)\mathsf U_\tau v
    \right|.
\]
On $\{\tau+h<\zeta\}$ this follows in the bundle topology on
a common adapted nest; on the complementary event both vectors
and their approximants vanish.
Fatou's lemma and the preceding estimate show that this limit
is zero almost surely.
Taking a common event of probability one for a countable dense
set of initial vectors and using boundedness of the operators,
we conclude that
\[
    \mathsf U_{\tau+h}
    =
    (\mathsf U_h\circ\theta_\tau)\mathsf U_\tau
    \qquad\mathbb P_x\text{-a.s.}
\]
The case $h=0$ follows from $\mathsf U_0=I$.
The exceptional set $N$ was chosen independently of $\tau$
and $h$, so this proves strong multiplicativity.
\end{itemize}
\end{proof}

\subsection{From multiplicative transport to bundle forms}\label{subsec:transport-bundle-forms}
Let $\widetilde{\mathcal H}$ be a quasi-Hilbert bundle over $(X,m,\mathcal{E})$.
\begin{dfn}\label{dfn:transport-association}
    Let $(\mathsf U_t)_{t\geq0}$ be a multiplicative stochastic transport
of $(\mathsf X_t)_{t\geq0}$ on $\widetilde{\mathcal H}$,
and let $\breve{\mathcal E}$ be a bundle form on $\mathrm L^2(m;\widetilde{\mathcal H})$. 
\begin{itemize}
    \item $(\mathsf U_t)_{t\geq0}$ is said to be associated
    with $\breve{\mathcal E}$ if, for every bounded
    measurable section $\alpha$ of $\widetilde{\mathcal H}$ belonging to
    $\mathrm L^2(m;\widetilde{\mathcal H})$ and every $t\geq0$,
    \[
        \breve p_t^{\mathsf U}\alpha
        =
        \breve T_t\alpha
        \qquad m\text{-a.e.}
    \]
    \item $(\mathsf U_t)_{t\geq0}$ is said to be properly
    associated with $\breve{\mathcal E}$ if it is associated
    with $\breve{\mathcal E}$ and, for every bounded measurable
    section $\alpha$ of $\widetilde{\mathcal H}$ belonging to
    $\mathrm L^2(m;\widetilde{\mathcal H})$ and every $t>0$, 
    $\breve p_t^{\mathsf U}\alpha$ is a quasi-continuous representative of $\breve T_t\alpha$.
\end{itemize}
\end{dfn}
\begin{rem}
    Association is equivalently characterized by
$\breve R_\lambda^{\mathsf U}\alpha=\breve G_\lambda\alpha$
$m$-a.e.\ for every sufficiently large $\lambda$ and every bounded
measurable section $\alpha\in\mathrm L^2(m;\widetilde{\mathcal H})$. If $(\mathsf U_t)_{t\geq0}$ is contractive, Association is equivalently characterized by $\breve R_\lambda^{\mathsf U}\alpha$ is a quasi-continuous representative of $\breve G_\lambda\alpha$ for every sufficiently large $\lambda$ and every bounded measurable section $\alpha\in\L^2(m;\widetilde{\mathcal H})$.
\end{rem}
Clearly, every multiplicative stochastic transport associated
with a bundle form is symmetric.

\begin{thm}[Kato form associated with a transport]\label{thm:transport-kato-form}
Let $(\mathsf U_t)_{t\geq0}$ be a symmetric contractive
quasi-right transport of $(\mathsf X_t)_{t\geq0}$ on $\widetilde{\mathcal H}$.
Then there exists a unique Kato form $\breve{\mathcal E}$ 
such that $(\mathsf U_t)_{t\geq0}$ is properly associated
with $\breve{\mathcal E}$.
\end{thm}
\begin{proof} Lemmas~\ref{lem:transport-bundle-form} and~\ref{lem:association-is-proper} give existence and proper association of the bundle form. Contractivity and Theorem~\ref{thm:kato-domination} give its Kato property.
\begin{lem}\label{lem:transport-bundle-form}
Let
$(\mathsf U_t)_{t\geq0}$ be a symmetric quasi-right
transport of $(\mathsf X_t)_{t\geq0}$ on $\widetilde{\mathcal H}$.
Then there exists a unique bundle form
$\breve{\mathcal E}$ on $\mathrm L^2(m;\widetilde{\mathcal H})$
such that $(\mathsf U_t)_{t\geq0}$ is associated
with $\breve{\mathcal E}$.
\end{lem}
\begin{proof}
For every $T>0$, the local integrability bound gives
\[
    \|\breve p_t^{\mathsf U}\alpha\|_\infty
    \leq C_T\|\alpha\|_\infty,
    \qquad 0\leq t\leq T,
\]
for every bounded measurable section $\alpha$.
By symmetry and duality, we also have
\[
    \|\breve p_t^{\mathsf U}\alpha\|_1
    \leq C_T\|\alpha\|_1,
    \qquad 0\leq t\leq T,
\]
for every bounded measurable section
$\alpha\in\mathrm L^1(m;\widetilde{\mathcal H})$.
Interpolation and density therefore yield a unique
bounded extension $\breve T_t$ to $\mathrm L^2(m;\widetilde{\mathcal H})$,
with
\[
    \|\breve T_t\|_{\mathrm{op}}\leq C_T,
    \qquad 0\leq t\leq T.
\]
Multiplicativity and symmetry imply that
$(\breve T_t)_{t\geq0}$ is a semigroup of self-adjoint
operators. It remains to prove strong continuity.
The bounded quasi-continuous sections belonging to
$\mathrm L^1(m;\widetilde{\mathcal H})$ form a dense subspace of
$\mathrm L^2(m;\widetilde{\mathcal H})$.
Indeed, normalize a quasi-fundamental family to obtain
bounded sections, multiply them by functions in
$C_b(X)\cap\mathrm L^1(m)$, and take finite linear combinations.
Fix such a section $\alpha$, and denote its fiberwise
Riesz dual by $\alpha^\flat$.
Choose a common nest for the quasi-continuity of $\alpha$
and the quasi-right-continuity of
$(\mathsf U_t\alpha^\flat(\mathsf X_0))_{t\geq0}$.
Since the nest exit times increase to the lifetime,
continuity of the pairing and $\mathsf U_0=I$ yield
\[
    \left\langle
        \alpha(\mathsf X_t),
        \mathsf U_t\alpha^\flat(\mathsf X_0)
    \right\rangle
    \longrightarrow|\alpha(x)|^2
    \qquad(t\downarrow0),
\]
$\mathbb P_x$-a.s.\ for q.e.\ $x$.
Moreover, for $0\leq t\leq T$,
\[
    \left|
        \left\langle
            \alpha(\mathsf X_t),
            \mathsf U_t\alpha^\flat(\mathsf X_0)
        \right\rangle
    \right|
    \leq
    \|\alpha\|_\infty|\alpha(x)|
    \sup_{0\leq r\leq T}\|\mathsf U_r\|_{\mathrm{op}}
    \qquad\mathbb P_x\text{-a.s.}
\]
The right-hand side is integrable with respect to
$m(dx)\mathbb P_x(d\omega)$, since its integral is at most
$C_T\|\alpha\|_\infty\|\alpha\|_1$.
Dominated convergence thus gives
\begin{align*}
    \langle\breve T_t\alpha,\alpha\rangle_{\mathrm L^2(m;\widetilde{\mathcal H})}
    &=
    \int_X\mathbb E_x\left[
        \left\langle
            \alpha(\mathsf X_t),
            \mathsf U_t\alpha^\flat(\mathsf X_0)
        \right\rangle
    \right]\,dm(x)\\
    &\longrightarrow\|\alpha\|_2^2.
\end{align*}
By self-adjointness and the semigroup property,
\[
    \|\breve T_t\alpha-\alpha\|_2^2
    =
    \langle\breve T_{2t}\alpha,\alpha\rangle_{\mathrm L^2(m;\widetilde{\mathcal H})}
    -2\langle\breve T_t\alpha,\alpha\rangle_{\mathrm L^2(m;\widetilde{\mathcal H})}
    +\|\alpha\|_2^2
    \longrightarrow0.
\]
Density and local uniform boundedness extend this
convergence to every $\alpha\in\mathrm L^2(m;\widetilde{\mathcal H})$.
The correspondence between strongly continuous
self-adjoint semigroups and densely defined,
lower-semibounded, closed symmetric forms \cite[Section~1.3]{FukushimaOshimaTakeda2011}, after a constant shift, now yields
a bundle form $\breve{\mathcal E}$ associated with
$(\mathsf U_t)_{t\geq0}$.
Uniqueness follows from the uniqueness of the
$\mathrm L^2$-extensions and this correspondence.
\end{proof}
\begin{lem}\label{lem:association-is-proper}
Let $(\mathsf U_t)_{t\geq0}$ be a contractive
quasi-right transport of $(\mathsf X_t)_{t\geq0}$ on $\widetilde{\mathcal H}$.
If $(\mathsf U_t)_{t\geq0}$ is associated
with $\breve{\mathcal E}$, then it is properly associated
with $\breve{\mathcal E}$.
\end{lem}
\begin{proof}\;
It suffices to show that, for every $\lambda>0$ and every bounded
measurable section $\alpha\in\L^2(m;\widetilde{\mathcal H})$, the section
$\breve R_\lambda^{\mathsf U}\alpha$ is quasi-continuous in
$\widetilde{\mathcal H}$. By association, $\breve R_\lambda^{\mathsf U}\alpha=\breve G_\lambda \alpha$
$m$-a.e., and by contractivity, $|\breve R_\lambda^{\mathsf U}\alpha|\leq\|\alpha\|_\infty/\lambda$.
We verify quasi-continuity of its pairings and its norm separately.
\begin{itemize}
    \item We first show that $\langle \breve R_\lambda^{\mathsf U}\alpha,\xi\rangle$ is quasi-continuous
    for every bounded quasi-continuous dual section $\xi$.
    By the scalar fine-continuity criterion, it suffices to show that,
    for q.e.\ $x$,
    \[
        \langle (\breve R_\lambda^{\mathsf U}\alpha)(\mathsf X_s),\xi(\mathsf X_s)\rangle
        \longrightarrow \langle (\breve R_\lambda^{\mathsf U}\alpha)(x),\xi(x)\rangle
        \quad(s\downarrow0),\qquad \mathbb P_x\text{-a.s.}
    \]
    We use here the equivalence between quasi-continuity and q.e.\ fine
    continuity in \cite[Theorems~4.6.1 and~A.2.7]{FukushimaOshimaTakeda2011},
    applied in a regular representation
    \cite[Chapter~VI]{MaRockner1992}.
    To compare $\xi(\mathsf X_s)$ with the initial value of $\xi$, use
    the lift $(\mathsf U_s\xi(\mathsf X_0))_{s\geq0}$ from the definition
    of a quasi-right transport. For q.e.\ $x$, splitting the resolvent
    integral at any bounded stopping time $\tau$ and applying strong multiplicativity (Proposition~\ref{prp:quasi-right-strong-multiplicativity}) and the strong Markov property gives
    \begin{align*}
        &\mathbb E_x\!\left[
            \int_0^\infty e^{-\lambda r}
            \langle \alpha(\mathsf X_r),\mathsf U_r\xi(\mathsf X_0)\rangle\,dr
            \,\middle|\,\mathcal F_\tau^{\mathsf X,x}
        \right]\\
        &\qquad =
        \int_0^\tau e^{-\lambda r}
            \langle \alpha(\mathsf X_r),\mathsf U_r\xi(\mathsf X_0)\rangle\,dr        
        +e^{-\lambda\tau}
            \langle (\breve R_\lambda^{\mathsf U}\alpha)(\mathsf X_\tau),\mathsf U_\tau\xi(\mathsf X_0)\rangle .
    \end{align*}
    The conditional expectation has a right-continuous martingale version.
    Since the identity holds at every bounded stopping time, the optional
    section theorem \cite[Theorem~IV.84]{DellacherieMeyer1978} identifies the optional process on the right with this
    version up to indistinguishability. The integral term is continuous
    and $\mathsf U_0=I$, so
    \[
        \langle (\breve R_\lambda^{\mathsf U}\alpha)(\mathsf X_s),\mathsf U_s\xi(\mathsf X_0)\rangle
        \longrightarrow
        \langle (\breve R_\lambda^{\mathsf U}\alpha)(\mathsf X_0),\xi(\mathsf X_0)\rangle
        \qquad(s\downarrow0),\quad\mathbb P_x\text{-a.s.}
    \]
    By quasi-right-continuity of this lift and quasi-continuity of $\xi$,
    both $\mathsf U_s\xi(\mathsf X_0)$ and $\xi(\mathsf X_s)$ converge
    to $\xi(\mathsf X_0)$ on a common nest as $s\downarrow0$. Hence
    \[
        |\mathsf U_s\xi(\mathsf X_0)-\xi(\mathsf X_s)|
        \longrightarrow 0
        \qquad(s\downarrow0),\quad\mathbb P_x\text{-a.s.}
    \]
    The asserted convergence follows since $\breve R_\lambda^{\mathsf U}\alpha$ is bounded.

    \item We next show that
    $|\breve R_\lambda^{\mathsf U}\alpha|
    \leq\widetilde{|\breve G_\lambda\alpha|}$ q.e.
    Contractivity and association give semigroup domination by $T_t$.
    Hence Theorem~\ref{thm:kato-domination} gives
    $|\breve G_\lambda\alpha|\in\D(\mathcal E)$, so this scalar function
    has a quasi-continuous representative. Splitting the resolvent at a
    deterministic time $s>0$ yields
    \[
        |(\breve R_\lambda^{\mathsf U}\alpha)(x)|
        \leq s\|\alpha\|_\infty
        +e^{-\lambda s}
        p_s^{\mathsf X}\bigl|\breve R_\lambda^{\mathsf U}\alpha\bigr|(x).
    \]
    By association, $|\breve R_\lambda^{\mathsf U}\alpha|
    =\widetilde{|\breve G_\lambda\alpha|}$ $m$-a.e.
    Proper association of the underlying process therefore allows us to
    replace the function inside $p_s^{\mathsf X}$ by this quasi-continuous
    representative for q.e.\ $x$, simultaneously for $s=1/n$, $n\geq1$.
    Since the representative is bounded q.e., its quasi-continuity gives
    \[
        p_{1/n}^{\mathsf X}\bigl(\widetilde{|\breve G_\lambda\alpha|}\bigr)(x)
        \longrightarrow \widetilde{|\breve G_\lambda\alpha|}(x)
        \qquad\text{for q.e.\ }x.
    \]
    Letting $n\to\infty$ gives the desired inequality.

    \item We finally prove the reverse inequality by stopping at a set on
    which $|\breve R_\lambda^{\mathsf U}\alpha|$ has a prescribed upper
    bound. Fix $a\geq0$ and a compact set
    $F\subset\{|\breve R_\lambda^{\mathsf U}\alpha|\leq a\}$.
    Let $\sigma_F$ be its hitting time. The discounted hitting probability
    $e_F(x):=\mathbb E_x[e^{-\lambda\sigma_F};\,\sigma_F<\zeta]$
    is quasi-continuous and equals $1$ q.e.\ on $F$
    \cite[Sections~4.2 and~4.6]{FukushimaOshimaTakeda2011}.
    Split the resolvent at $\sigma_F$, first stopping at $\sigma_F\wedge n$
    and then letting $n\to\infty$. By contractivity, the contribution
    before the hit is bounded by
    $\lambda^{-1}\|\alpha\|_\infty(1-e_F(x))$. At the hit,
    $|(\breve R_\lambda^{\mathsf U}\alpha)(\mathsf X_{\sigma_F})|\leq a$, so
    \[
        |(\breve R_\lambda^{\mathsf U}\alpha)(x)|
        \leq\lambda^{-1}\|\alpha\|_\infty(1-e_F(x))+a\,e_F(x).
    \]
    With $\widetilde{|\breve G_\lambda\alpha|}$ on the left, this inequality
    holds $m$-a.e., hence q.e., since both sides are quasi-continuous.
    In particular, $\widetilde{|\breve G_\lambda\alpha|}\leq a$ q.e.\ on $F$.
    Thus $\{|\breve R_\lambda^{\mathsf U}\alpha|\leq a
    <\widetilde{|\breve G_\lambda\alpha|}\}$ is exceptional, since otherwise
    it would contain a compact subset of positive capacity, contradicting
    the bound on $F$. Taking the union over rational $a\geq0$ gives the
    reverse inequality q.e.
\end{itemize}
Thus $|\breve R_\lambda^{\mathsf U}\alpha|=\widetilde{|\breve G_\lambda\alpha|}$
q.e., and the first step gives quasi-continuity of the pairings with a
bounded quasi-fundamental family of dual sections. Hence
$\breve R_\lambda^{\mathsf U}\alpha$ is a quasi-continuous representative
of $\breve G_\lambda\alpha$ in $\widetilde{\mathcal H}$.
\end{proof}
\end{proof}

\begin{prp}\label{prp:transport-quasi-realization}
Let $(\mathsf U_t)_{t\geq0}$ be a contractive
quasi-right transport on $\widetilde{\mathcal H}$ that is
properly associated with a bundle form $\breve{\mathcal E}$.
Then $(\widetilde{\mathcal H},\operatorname{id})$ is a
quasi-realization of $\breve{\mathcal E}$.
\end{prp}
\begin{proof}
By proper association, the sections $\breve T_t\alpha$,
with $t>0$ and $\alpha\in\mathrm L^2\cap\mathrm L^\infty$,
admit quasi-continuous representatives in
$\widetilde{\mathcal H}$.
These sections are dense in $\D(\breve{\mathcal E})$
in the form norm. The estimate
\[
    \mathcal E_1(|u-v|)
    \leq\breve{\mathcal E}_1(u-v)
\]
and the scalar capacity estimate allow us to extract
a subsequence of any such form approximation converging
uniformly in fiber norm on each member of a common nest.
Thus every element of $\D(\breve{\mathcal E})$ has a
quasi-continuous representative. Uniqueness follows by
applying scalar quasi-continuous uniqueness to the norm
of the difference.
Next, choose a bounded quasi-fundamental family $(\eta_j)_{j\geq1}$
with $\eta_j\in\mathrm L^2(m;\widetilde{\mathcal H})$,
by normalization and multiplication by a bounded,
strictly positive quasi-continuous $\mathrm L^2$ function.
Set $\beta_{j,n}:=\breve p_{1/n}^{\mathsf U}\eta_j$.
Each $\beta_{j,n}$ is a quasi-continuous representative
of an element of $\D(\breve{\mathcal E})$.
Quasi-right-continuity, $\mathsf U_0=I$, contractivity,
and dominated convergence give, for q.e.\ $x$,
\[
\begin{aligned}
    |\beta_{j,n}(x)-\eta_j(x)|^2
    &\leq
    p_{1/n}^{\mathsf X}|\eta_j|^2(x)
    -2\langle\beta_{j,n}(x),\eta_j(x)\rangle
    +|\eta_j(x)|^2\\
    &\longrightarrow0.
\end{aligned}
\]
Hence $(\beta_{j,n})_{j,n\geq1}$ has dense linear span
in every fiber outside an exceptional set.
Adjoining the quasi-continuous representatives of a
countable form-dense subset of $\D(\breve{\mathcal E})$
gives a quasi-fundamental family representing the original
quasi-structure and satisfying the second condition
of a quasi-realization.
\end{proof}

\subsection{Uniqueness of multiplicative transport}
\begin{thm}[Uniqueness of associated transport]\label{thm:multiplicative-transport-uniqueness}
Let $(\mathsf U_t)_{t\geq0}$ and $(\mathsf U'_t)_{t\geq0}$
be two quasi-right transports of $(\mathsf X_t)_{t\geq0}$
on $\widetilde{\mathcal H}$ associated with the same bundle form
$\breve{\mathcal E}$ on $\L^2(m;\widetilde{\mathcal H})$.
Then
\[
    \mathbb P_x\!\left(
        \mathsf U_t=\mathsf U'_t,\;
        \forall t\geq0
    \right)=1
    \qquad\text{for q.e.\ }x\in X.
\]
\end{thm}
\begin{proof}
    Combine Lemmas~\ref{lem:transport-uniqueness-fixed-time} and~\ref{lem:transport-uniqueness-qe} below.
\begin{lem}\label{lem:transport-uniqueness-fixed-time}
Let $(\mathsf U_t)_{t\geq0}$ and $(\mathsf U'_t)_{t\geq0}$
be two multiplicative stochastic transports of $(\mathsf X_t)_{t\geq0}$
on $\widetilde{\mathcal H}$ associated with the same bundle form
$\breve{\mathcal E}$ on $\L^2(m;\widetilde{\mathcal H})$.
Then for all $t>0$,
\[
    \mathbb P_x\!\left(
        \mathsf U_t=\mathsf U'_t
    \right)=1
    \qquad\text{for }m\text{-a.e. }x\in X.
\]
\end{lem}
\begin{proof}
Fix $t>0$, a bounded measurable section
$\alpha\in\L^2(m;\widetilde{\mathcal H})$, and a cylinder function
\[
    F=\prod_{j=1}^n f_j(\mathsf X_{t_j}),
    \qquad
    0\leq t_1<\cdots<t_n\leq t,
\]
where $f_1,\ldots,f_n$ are bounded Borel functions on $X$.
By multiplicativity, the Markov property, and the fact that
both transports are associated with $\breve{\mathcal E}$,
we have
\[
\begin{aligned}
    \mathbb E_x\!\left[
        F\mathsf U_t^*\bigl(\alpha(\mathsf X_t)\bigr)
    \right]
    &=
    \breve T_{t_1}\!\left(
        f_1\breve T_{t_2-t_1}\!\left(
            f_2\cdots
            \breve T_{t_n-t_{n-1}}\!\left(
                f_n\breve T_{t-t_n}\alpha
            \right)\cdots
        \right)
    \right)(x)
\end{aligned}
\]
for $m$-a.e.\ $x\in X$. The same identity holds for $\mathsf U'$, and hence
\[
    \mathbb E_x\!\left[
        F\bigl(\mathsf U_t^*-(\mathsf U'_t)^*\bigr)
        \bigl(\alpha(\mathsf X_t)\bigr)
    \right]
    =0
    \qquad\text{for }m\text{-a.e. }x\in X. \tag{$\ast$}
\]
Fix $h>0$ and a bounded cylinder function $G$ up to time $h$.
By the Markov property at time $t$, we have
\[
\begin{aligned}
    &\mathbb E_x\!\left[F(G\circ\theta_t)
        \bigl(\mathsf U_t^*-(\mathsf U'_t)^*\bigr)
        \bigl(\alpha(\mathsf X_t)\bigr)
    \right]
    \\
    &\quad=
    \mathbb E_x\!\left[F\mathbb E_{\mathsf X_t}[G]
        \bigl(\mathsf U_t^*-(\mathsf U'_t)^*\bigr)
        \bigl(
        \alpha(\mathsf X_t)\bigr)
    \right]
    =0
\end{aligned}
\]
for $m$-a.e.\ $x\in X$.
Here the last equality follows from the preceding identity
applied to $(\mathbb E_\cdot[G])\alpha$, which is again bounded
and belongs to $\L^2(m;\widetilde{\mathcal H})$.
Consequently, by the monotone class theorem, $(\ast)$ holds
for every bounded $\mathcal F^\mathsf X_{t+}$-measurable function $F$. Choose a quasi-fundamental family $(\alpha_k)_{k\in\mathbb N}$
of $\widetilde{\mathcal H}$ consisting of bounded
$\L^2(m;\widetilde{\mathcal H})$-sections.
For $j,k\in\mathbb N$, set
\[
    F_{j,k}
    :=
    \operatorname{sgn}\!\left(
        \left\langle
            \bigl(\mathsf U_t^*-(\mathsf U'_t)^*\bigr)
            \alpha_k(\mathsf X_t),
            \alpha_j(\mathsf X_0)
        \right\rangle
    \right),
\]
where $\operatorname{sgn}(0)=0$.
This is a bounded $\mathcal F^\mathsf X_{t+}$-measurable function.
Applying $(\ast)$ with $\alpha=\alpha_k$ and $F=F_{j,k}$,
and pairing the resulting identity with $\alpha_j(x)$,
we obtain
\[
    \mathbb E_x\!\left[
        \left|
            \left\langle
                \bigl(\mathsf U_t^*-(\mathsf U'_t)^*\bigr)
                \alpha_k(\mathsf X_t),
                \alpha_j(\mathsf X_0)
            \right\rangle
        \right|
    \right]
    =0
    \qquad\text{for }m\text{-a.e. }x\in X.
\]
By countability, it follows that
\[
    \left\langle
        \bigl(\mathsf U_t^*-(\mathsf U'_t)^*\bigr)
        \alpha_k(\mathsf X_t),
        \alpha_j(\mathsf X_0)
    \right\rangle
    =0
    \qquad\text{for all }j,k\in\mathbb N,
    \quad\mathbb P_x\text{-a.s.},
\]
for $m$-a.e.\ $x\in X$.
Since the quasi-fundamental family is total outside an
exceptional set, its values are total both
at $x$ and at $\mathsf X_t$, $\mathbb P_x$-a.s.,
for $m$-a.e.\ $x\in X$.
The boundedness of the two operators therefore yields
\[
    \mathsf U_t^*=(\mathsf U'_t)^*
    \qquad\mathbb P_x\text{-a.s.},
    \qquad\text{for }m\text{-a.e. }x\in X.
\]
Taking duals proves the assertion.
\end{proof}

\begin{lem}\label{lem:transport-uniqueness-qe}
Let $(\mathsf U_t)_{t\geq0}$ and $(\mathsf U'_t)_{t\geq0}$
be two quasi-right transports of $(\mathsf X_t)_{t\geq0}$
on $\widetilde{\mathcal H}$.
Suppose
\[
    \mathbb P_x\!\left(
        \mathsf U_q=\mathsf U'_q,\;
        \forall q\in\mathbb Q_+
    \right)=1
    \qquad\text{for }m\text{-a.e. }x\in X. 
\]
Then
\[
    \mathbb P_x\!\left(
        \mathsf U_t=\mathsf U'_t,\;
        \forall t\geq0
    \right)=1
    \qquad\text{for q.e.\ }x\in X.
\]
\end{lem}
\begin{proof}
    Set $r_n:=2^{-n}$ for $n\ge1$.
We first observe that, for q.e.\ $x\in X$,
\[
    \mathsf U_q\circ\theta_{r_n}
    =
    \mathsf U'_q\circ\theta_{r_n}
    \qquad
    \text{for all }q\in\mathbb Q_+
    \text{ and }n\ge1,
    \quad\mathbb P_x\text{-a.s.} \tag{$\star$}
\]
Indeed, choose a Borel $m$-null set $B$ outside which the hypothesis of the lemma
holds.
By the proper association of the underlying process,
$p_{r_n}^{\mathsf X}\mathbf1_B$ is a quasi-continuous
representative of zero. Hence, for q.e.\ $x$,
$p_{r_n}^{\mathsf X}(x,B)=0$ for every $n\geq1$.
The Markov property and the hypothesis therefore give
\[
    \mathsf U_q\circ\theta_{r_n}
    =
    \mathsf U'_q\circ\theta_{r_n},
    \qquad n\geq1,\quad q\in\mathbb Q_+,
\]
simultaneously $\mathbb P_x$-a.s.
Next, fix a bounded measurable section $\xi$ of $\widetilde{\mathcal H}^*$ and set
\[
    d_t\coloneqq
    \left|
        (\mathsf U_t-\mathsf U'_t)\xi(\mathsf X_0)
    \right|.
\]
By quasi-right-continuity and exhaustion up to the lifetime
by the nest exit times, $(d_t)_{t\geq0}$ has right-continuous
paths $\mathbb P_x$ a.s. for q.e.\ $x$.
Here both lifts vanish after the lifetime, and $d_0=0$
by the initial conditions.
Fix such an $x$ for which the preceding shifted identities
($\star$) hold.
Then $d_{r_n}\to0$, $\mathbb P_x$-a.s., and
\[
    d_{r_n}
    \leq
    \|\xi\|_\infty
    \sup_{0\leq s\leq1}
    \bigl(\|\mathsf U_s\|_{\mathrm{op}}+\|\mathsf U'_s\|_{\mathrm{op}}\bigr).
\]
The local operator-norm bounds and dominated convergence
therefore yield $\mathbb E_x[d_{r_n}]\to0$.
Fix a positive rational time $q$.
For $n$ sufficiently large that $r_n<q$, multiplicativity
and ($\star$) imply
\[
\begin{aligned}
    d_q
    &=
    \left|
        \bigl(\mathsf U_{q-r_n}\circ\theta_{r_n}\bigr)
        \bigl(\mathsf U_{r_n}-\mathsf U'_{r_n}\bigr)
        \xi(\mathsf X_0)
    \right|
    \\
    &\leq
    \|\mathsf U_{q-r_n}\circ\theta_{r_n}\|_{\mathrm{op}}\,d_{r_n}.
\end{aligned}
\]
Since $d_{r_n}$ is $\mathcal F_{r_n+}^{\mathsf X}$-measurable,
the Markov property gives
\begin{align*}
    \mathbb E_x[d_q]
    &\leq
    \mathbb E_x\left[
        d_{r_n}\,
        \mathbb E_{\mathsf X_{r_n}}
        [\|\mathsf U_{q-r_n}\|_{\mathrm{op}}]
    \right]\\
    &\leq C_q\,\mathbb E_x[d_{r_n}]
    \longrightarrow0.
\end{align*}
Thus $d_q=0$ simultaneously for all $q\in\mathbb Q_+$,
$\mathbb P_x$-a.s.
Right-continuity then gives $d_t=0$ simultaneously for
all $t\geq0$.
Finally, apply this argument to a countable family of bounded
measurable sections whose values have dense linear span
in $\widetilde{\mathcal H}_x^*$ for q.e.\ $x$.
Taking the union of their exceptional sets and the countable
intersection of the corresponding full-probability events,
boundedness of the operators yields
$\mathsf U_t=\mathsf U'_t$ simultaneously for all $t\geq0$,
$\mathbb P_x$-a.s., for q.e.\ $x$.

\end{proof}

\end{proof}

\subsection{Existence of multiplicative transport}
\begin{thm}[Existence of associated transport]\label{thm:multiplicative-transport-existence}
Let $\breve{\mathcal E}$ be a Kato form on $(X,m,\mathcal E)$.
Then there exists a contractive
quasi-right transport of $(\mathsf X_t)_{t\geq0}$ on $\mathcal H^{\breve{\mathcal E}}$ that is properly associated with $\breve{\mathcal E}$.
\end{thm}
\begin{proof}\;
Lemma~\ref{lem:transport-rational-times} constructs the transport at rational times, and Lemma~\ref{lem:transport-pairing-regularization} extends it to all times. Lemma~\ref{lem:transport-strong-multiplicativity-association} proves strong multiplicativity and association. Finally, Lemma~\ref{lem:transport-quasi-right-continuity} establishes quasi-right-continuity and proper association.
\begin{lem}\label{lem:transport-rational-times}
There exists a family $(\mathsf U_q)_{q\in\mathbb Q_+}$
of linear contractions
$
    \mathsf U_q(\omega)\colon
    (\mathcal H_{\mathsf X_0(\omega)}^{\breve{\mathcal E}})^*
    \rightarrow (\mathcal H_{\mathsf X_q(\omega)}^{\breve{\mathcal E}})^*
$
with the following properties:
\begin{enumerate}
    \item $\mathsf U_0(\omega)=I$ for every $\omega\in\Omega$.

    \item For every bounded measurable section $\xi$ of $(\mathcal H^{\breve{\mathcal E}})^*$,
    the family
    $(\mathsf U_q\xi(\mathsf X_0))_{q\in\mathbb Q_+}$
    is $(\mathcal F_q^{\mathsf X})_{q\in\mathbb Q_+}$-adapted.

    \item For q.e.\ $x\in X$ and every $q,h\in\mathbb Q_+$,
    \[
        \mathsf U_{q+h}
        =
        (\mathsf U_h\circ\theta_q)\mathsf U_q
        \quad\mathbb P_x\text{-a.s.}
    \]

    \item For every $\alpha\in\D(\breve{\mathcal E})^\infty$,
    for q.e.\ $x\in X$ and every $q\in\mathbb Q_+$,
    \[
        \mathbb E_x\left[
            \mathsf U_q^*\tilde\alpha(\mathsf X_q)
        \right]
        =
        (\widetilde{\breve T_q\alpha})(x).
    \]
\end{enumerate}
\end{lem}
\begin{proof}\;
\begin{itemize}
\item We construct a contraction from cylinder sections observed up to each rational time $q$.
Choose a countable $\mathscr A\subset 2^X$ generating
the Borel $\sigma$-algebra $\mathcal B(X)$ which includes $X$ and is closed under taking unions and intersections, and choose a countable $\mathbb Q$-vector space
$\mathscr C$ of bounded measurable sections belonging to
$\L^2(m;\mathcal H^{\breve{\mathcal E}})$, containing a bounded quasi-fundamental family.
Enlarge $\mathscr C$ so that it is stable under multiplication
by $\mathbf1_B$, $B\in\mathscr A$, and under taking
quasi-continuous representatives of $\breve T_r u$,
$r\in\mathbb Q_+\setminus\{0\}$.
By countability, semigroup domination, and proper association
of the underlying process, we may choose a properly
exceptional Borel set $N$ such that the relevant semigroup
identities and domination inequalities hold on $X\setminus N$,
and the chosen fundamental family is fiberwise total there.
Fix $q\in\mathbb Q_+\setminus\{0\}$.
For $x\in X\setminus N$ and a cylinder section
\[
    \Phi=F\alpha(\mathsf X_q),
    \qquad
    F=\prod_{j=1}^n\mathbf1_{B_j}(\mathsf X_{t_j}),
\]
where $\alpha\in\mathscr C$, $B_j\in\mathscr A$, and
$0<t_1<\cdots<t_n<q$ are rational, set
\[
    \Lambda_{x,q}\Phi
    :=
    \widetilde{
        \breve T_{t_1}\!\left(
            \mathbf1_{B_1}\breve T_{t_2-t_1}\!\left(
                \mathbf1_{B_2}\cdots
                \breve T_{t_n-t_{n-1}}\!\left(
                    \mathbf1_{B_n}\breve T_{q-t_n}\alpha
                \right)\cdots
            \right)
        \right)
    }(x)
    \in \mathcal H_x^{\breve{\mathcal E}}.
\]
When $n=0$, define
\[
    \Lambda_{x,q}\bigl(\alpha(\mathsf X_q)\bigr)
    :=
    (\widetilde{\breve T_q\alpha})(x).
\]
Extend this assignment formally by $\mathbb Q$-linearity
to finite sums. Passing to a common time grid and refining the sets from
$\mathscr A$, we may therefore write
\[
    \Phi
    =
    \sum_{\ell=1}^{L}
    \mathbf1_{C_\ell}\alpha_\ell(\mathsf X_q),
\]
where the $C_\ell$ are pairwise disjoint cylinder sets and
$\alpha_\ell\in\mathscr C$.
Repeated application of semigroup domination gives
\[
    \left|
        \Lambda_{x,q}
        \bigl(\mathbf1_{C_\ell}\alpha_\ell(\mathsf X_q)\bigr)
    \right|
    \leq
    \mathbb E_x\!\left[
        \mathbf1_{C_\ell}|\alpha_\ell(\mathsf X_q)|
    \right].
\]
Consequently, disjointness yields
\[
\begin{aligned}
    |\Lambda_{x,q}\Phi|
    &\leq
    \sum_{\ell=1}^{L}
    \mathbb E_x\!\left[
        \mathbf1_{C_\ell}|\alpha_\ell(\mathsf X_q)|
    \right]
    =\mathbb E_x[|\Phi|].
\end{aligned}
\]
Applying this estimate to the difference of two
representations shows that the assigned value depends
only on the $\mathbb P_x$-equivalence class of $\Phi$.
Thus $\Lambda_{x,q}$ is well-defined and
$\mathbb Q$-linear on the rational span of the cylinder
sections. Since the rational span of these cylinder sections is dense in
$\L^1(\Omega,\mathcal F_q^{\mathsf X},\mathbb P_x;\mathcal H_{\mathsf X_q}^{\breve{\mathcal E}})$,
$\Lambda_{x,q}$ extends uniquely to a linear contraction
\[
    \Lambda_{x,q}\colon
    \L^1(\Omega,\mathcal F_q^{\mathsf X},\mathbb P_x;\mathcal H_{\mathsf X_q}^{\breve{\mathcal E}})
    \longrightarrow \mathcal H_x^{\breve{\mathcal E}}.
\]

\item We represent this contraction by a random operator and use its dual to define $\mathsf U_q$.
For every bounded $\mathcal F_q^{\mathsf X}$-measurable
section $\Phi$ of $\mathcal H_{\mathsf X_q}^{\breve{\mathcal E}}$, the section
$x\mapsto\Lambda_{x,q}\Phi$ is measurable by approximation
with measurably chosen cylinder sections and the contraction
estimate. For $\alpha\in\mathscr C$ and $x\in X\setminus N$, define
the $\mathcal H_x^{\breve{\mathcal E}}$-valued measure
\[
    \nu_{\alpha,x}(A)
    :=
    \Lambda_{x,q}\bigl(\mathbf1_A\alpha(\mathsf X_q)\bigr),
    \qquad A\in\mathcal F_q^{\mathsf X}.
\]
Continuity of $\Lambda_{x,q}$ gives countable additivity,
and the contraction estimate yields
\[
    |\nu_{\alpha,x}|(A)
    \leq
    \mathbb E_x\!\left[
        \mathbf1_A|\alpha(\mathsf X_q)|
    \right].
\]
Moreover, $x\mapsto\nu_{\alpha,x}(A)$ is measurable for
every $A\in\mathcal F_q^{\mathsf X}$.
The parameter-dependent Hilbert-valued Radon--Nikodym
theorem therefore provides a jointly measurable section
$b_\alpha(x,\omega)\in \mathcal H_x^{\breve{\mathcal E}}$ such that
\[
    \nu_{\alpha,x}(A)
    =
    \mathbb E_x\!\left[
        \mathbf1_A b_\alpha(x,\cdot)
    \right],
    \qquad A\in\mathcal F_q^{\mathsf X},
\]
and
\[
    |b_\alpha(x,\omega)|
    \leq
    |\alpha(\mathsf X_q(\omega))|
    \qquad\mathbb P_x\text{-a.s.}
\]
By uniqueness of the densities and countability of
$\mathscr C$, outside a jointly measurable set whose
section at each $x$ is $\mathbb P_x$-null, the map
$\alpha\mapsto b_\alpha(x,\omega)$ is $\mathbb Q$-linear
and the preceding bound holds simultaneously for all
$\alpha\in\mathscr C$.
We also exclude the set where $\mathsf X_q\in N$, which
has $\mathbb P_x$-measure zero for every $x\in X\setminus N$.
On the remaining set, the assignment
\[
    \alpha(\mathsf X_q(\omega))
    \longmapsto b_\alpha(x,\omega),
    \qquad \alpha\in\mathscr C,
\]
is well-defined and contractive.
By fiberwise totality, it extends uniquely to a linear
contraction
\[
    B_q(x,\omega)\colon
    \mathcal H_{\mathsf X_q(\omega)}^{\breve{\mathcal E}}
    \longrightarrow \mathcal H_x^{\breve{\mathcal E}}.
\]
Extend $B_q$ by zero elsewhere. $B_q$ is a measurable family of operators since its values on the countable generating family are jointly measurable.
The defining identities for the densities, followed by
density and continuity, give
\[
    \Lambda_{x,q}\Phi
    =
    \mathbb E_x\!\left[B_q(x,\cdot)\Phi\right],
    \qquad
    \Phi\in
    \L^1(\Omega,\mathcal F_q^{\mathsf X},
         \mathbb P_x;\mathcal H_{\mathsf X_q}^{\breve{\mathcal E}}),
\]
for every $x\in X\setminus N$.
Perform this construction for every positive rational $q$
and set
\[
    \mathsf U_q(\omega)
    :=
    B_q(\mathsf X_0(\omega),\omega)^*,
    \qquad q>0,
    \qquad
    \mathsf U_0(\omega):=I.
\]
These are linear contractions and satisfy properties
$1$ and $2$. We next prove property $3$.

\item We prove multiplicativity at rational times.
Fix positive rational numbers $q,h$.
Let $F$ and $G$ be cylinder functions up to times $q$ and $h$, respectively, formed from rational observation times
and indicators of sets in $\mathscr A$.
For $\alpha\in\mathscr C$ and $x\in X\setminus N$,
splitting the defining iterated semigroup expression
at time $q$ gives
\[
\begin{aligned}
    &\mathbb E_x\!\left[
        F(G\circ\theta_q)\,
        \mathsf U_{q+h}^*\alpha(\mathsf X_{q+h})
    \right]
    \\
    &\quad=
    \mathbb E_x\!\left[
        F\mathsf U_q^*
        \Lambda_{\mathsf X_q,h}
        \bigl(G\alpha(\mathsf X_h)\bigr)
    \right]
    \\
    &\quad=
    \mathbb E_x\!\left[
        F(G\circ\theta_q)\,
        \mathsf U_q^*(\mathsf U_h^*\circ\theta_q)
        \alpha(\mathsf X_{q+h})
    \right].
\end{aligned}
\]
The first equality uses the representation of
$\Lambda$ by $\mathsf U^*$ and the stability
of $\mathscr C$ under the operations appearing in
the iterated semigroup expression.
The second follows from the Markov property at time $q$
and the identity
\[
    \Lambda_{y,h}\bigl(G\alpha(\mathsf X_h)\bigr)
    =\mathbb E_y\!\left[
        G\mathsf U_h^*\alpha(\mathsf X_h)
    \right],
    \qquad y\in X\setminus N.
\]
Since the products $F(G\circ\theta_q)$ generate
$\mathcal F_{q+h}^{\mathsf X}$, the monotone class theorem extends the equality of the first and last expectations to every bounded $\mathcal F_{q+h}^{\mathsf X}$-measurable test function.
Both vector-valued random variables are
$\mathcal F_{q+h}^{\mathsf X}$-measurable and integrable,
so
\[
    \mathsf U_{q+h}^*\alpha(\mathsf X_{q+h})
    =
    \mathsf U_q^*(\mathsf U_h^*\circ\theta_q)
    \alpha(\mathsf X_{q+h})
    \qquad\mathbb P_x\text{-a.s.}
\]
By countability of $\mathscr C$, this holds simultaneously
for all $\alpha\in\mathscr C$.
Since $N$ is properly exceptional, fiberwise totality
at $\mathsf X_{q+h}$ and continuity of the operators yield
\[
    \mathsf U_{q+h}
    =
    (\mathsf U_h\circ\theta_q)\mathsf U_q
    \qquad\mathbb P_x\text{-a.s.}
\]
Taking a countable intersection over positive rational
$q,h$ gives simultaneous validity.
The cases $q=0$ or $h=0$ follow from $\mathsf U_0=I$.
This proves property $3$.

\item We extend the semigroup representation to every bounded section in the form domain.
Fix $\alpha\in\D(\breve{\mathcal E})^\infty$.
After modifying $\tilde\alpha$ on an exceptional set
if necessary, we may assume that it is bounded.
Choose a countable enlargement of $\mathscr C$
containing $\tilde\alpha$ and having the same closure
properties.
Repeating the construction of $\Lambda$ for this enlarged
family gives linear contractions $\Lambda'_{x,q}$
for every $x\in X\setminus N_\alpha$ and every positive
rational $q$, where $N_\alpha\supset N$ is a properly
exceptional Borel set.
The original and enlarged constructions agree on the
rational span of the original cylinder sections.
Since this span is dense in
$\L^1(\Omega,\mathcal F_q^{\mathsf X},
\mathbb P_x;\mathcal H_{\mathsf X_q}^{\breve{\mathcal E}})$,
uniqueness of the continuous extension gives
\[
    \Lambda'_{x,q}=\Lambda_{x,q},
    \qquad x\in X\setminus N_\alpha.
\]
Consequently,
\[
\begin{aligned}
    \mathbb E_x\!\left[
        \mathsf U_q^*\tilde\alpha(\mathsf X_q)
    \right]
    =
    \Lambda_{x,q}\bigl(\tilde\alpha(\mathsf X_q)\bigr)
    =
    \Lambda'_{x,q}\bigl(\tilde\alpha(\mathsf X_q)\bigr)
    =
    (\widetilde{\breve T_q\alpha})(x)
\end{aligned}
\]
for every $x\in X\setminus N_\alpha$ and every positive
rational $q$. This proves property $4$.
\end{itemize}
\end{proof}
\begin{lem}\label{lem:transport-pairing-regularization}
Let $(\mathsf U_q)_{q\in\mathbb Q_+}$ be the family
constructed in Lemma~\ref{lem:transport-rational-times}.
Then there exists a family $(\mathsf U_t)_{t\geq0}$
of linear contractions
$
    \mathsf U_t(\omega)\colon
    (\mathcal H_{\mathsf X_0(\omega)}^{\breve{\mathcal E}})^*
    \rightarrow (\mathcal H_{\mathsf X_t(\omega)}^{\breve{\mathcal E}})^*,
$
extending $(\mathsf U_q)_{q\in\mathbb Q_+}$ with the
following properties:
\begin{enumerate}
    \item For every bounded measurable section $\xi$
    of $(\mathcal H^{\breve{\mathcal E}})^*$, the lift
    $(\mathsf U_t\xi(\mathsf X_0))_{t\geq0}$
    is $(\mathcal F_{t+}^{\mathsf X})_{t\geq0}$-progressively
    measurable.

    \item For every $\alpha\in\D(\breve{\mathcal E})^\infty$,
    for q.e.\ $x\in X$ and every bounded measurable section
    $\xi$ of $(\mathcal H^{\breve{\mathcal E}})^*$,
    the process
    $
        t\mapsto
        \left\langle
            \tilde\alpha(\mathsf X_t),
            \mathsf U_t\xi(\mathsf X_0)
        \right\rangle
    $
    has right-continuous paths $\mathbb P_x$-a.s.
\end{enumerate}
\end{lem}
\begin{proof}\;
\begin{itemize}
\item We first obtain right-continuous versions of a countable family of scalar pairings by submartingale regularization.
By resolvent approximation, we can choose a countable family $(\alpha_k)_{k\geq1}\subset\D(\breve L)$,
dense in $\D(\breve{\mathcal E})$ in the form norm, such that
$\alpha_k,\breve L\alpha_k\in
\L^\infty(m;\mathcal H^{\breve{\mathcal E}})$ for every $k\geq1$.
By the quasi-realization property,
$(\tilde\alpha_k)_{k\geq1}$ is a quasi-fundamental family;
in particular,
\[
    \mathcal H_x^{\breve{\mathcal E}}
    =
    \overline{\operatorname{Span}_{\mathbb R}
    \{\tilde\alpha_k(x)\mid k\geq1\}}
    \quad\text{for q.e.\ }x\in X.
\]
Choose a measurable sections $(\eta_j)_{j\geq1}$ of $(\mathcal H^{\breve{\mathcal E}})^*$ which make the orthonormal basis of fibers outside an exceptional set.
For $q\in\mathbb Q_+$, set
\[
    Z_q^{k,j}
    \coloneqq
    \left\langle
        \tilde\alpha_k(\mathsf X_q),
        \mathsf U_q\eta_j(\mathsf X_0)
    \right\rangle.
\]
By multiplicativity, the Markov property, and property $4$
of Lemma~\ref{lem:transport-rational-times}, there exists a properly exceptional
Borel set $N$ such that, for every $x\in X\setminus N$,
$k,j\geq1$, and $q,h\in\mathbb Q_+$,
\[
\begin{aligned}
    \mathbb E_x\!\left[
        Z_{q+h}^{k,j}\mid\mathcal F_q^{\mathsf X,x}
    \right]
    &=
    \mathbb E_x\!\left[
        \left\langle
            \bigl(\mathsf U_h^*
            \tilde\alpha_k(\mathsf X_h)\bigr)\circ\theta_q,
            \mathsf U_q\eta_j(\mathsf X_0)
        \right\rangle
        \,\middle|\,\mathcal F_q^{\mathsf X,x}
    \right]
    \\
    &=
    \left\langle
        \mathbb E_{\mathsf X_q}\!\left[
            \mathsf U_h^*\tilde\alpha_k(\mathsf X_h)
        \right],
        \mathsf U_q\eta_j(\mathsf X_0)
    \right\rangle
    \\
    &=
    \left\langle
        (\widetilde{\breve T_h\alpha_k})(\mathsf X_q),
        \mathsf U_q\eta_j(\mathsf X_0)
    \right\rangle
    \qquad\mathbb P_x\text{-a.s.}
\end{aligned}
\]
Consequently, enlarging $N$ if necessary, we obtain
\[
\begin{aligned}
    \left|
        \mathbb E_x\!\left[
            Z_{q+h}^{k,j}-Z_q^{k,j}
            \,\middle|\,\mathcal F_q^{\mathsf X,x}
        \right]
    \right|
    &=
    \left|
        \left\langle
            (\widetilde{\breve T_h\alpha_k})(\mathsf X_q)
            -\tilde\alpha_k(\mathsf X_q),
            \mathsf U_q\eta_j(\mathsf X_0)
        \right\rangle
    \right|
    \\
    &\leq
    \|\breve T_h\alpha_k-\alpha_k\|_\infty
    \\
    &\leq
    \int_0^h
        \|\breve T_s\breve L\alpha_k\|_\infty\,ds
    \leq
    h\|\breve L\alpha_k\|_\infty
    \qquad\mathbb P_x\text{-a.s.},
\end{aligned}
\]
where we used contractivity of $\mathsf U_q$,
the generator identity in $\L^2$, and semigroup domination.
In particular,
\[
    \left|
        \mathbb E_x[Z_{q+h}^{k,j}]
        -
        \mathbb E_x[Z_q^{k,j}]
    \right|
    \leq h\|\breve L\alpha_k\|_\infty.
\]
Fix $x\in X\setminus N$.
Consequently, both processes
\[
    \left(
        Z_q^{k,j}+q\|\breve L\alpha_k\|_\infty
    \right)_{q\in\mathbb Q_+},
    \qquad
    \left(
        -Z_q^{k,j}+q\|\breve L\alpha_k\|_\infty
    \right)_{q\in\mathbb Q_+}
\]
are $(\mathcal F_q^{\mathsf X,x})$-submartingales
under $\mathbb P_x$.
The preceding expectation estimate shows
that their expectations are right-continuous along
rational times.
Since $(\mathcal F_t^{\mathsf X,x})_{t\geq0}$ satisfies
the usual conditions, the submartingale regularization theorem \cite{CohenElliott2015} shows
that, $\mathbb P_x$-a.s., the rational right limits
$
    Z_t^{k,j}
    :=
    \lim_{\substack{q\downarrow t\\q\in\mathbb Q_+}}
    Z_q^{k,j},
    \; t\geq0,
$
exist simultaneously for all $t\geq0$, have c\`adl\`ag
paths, and agree with the original variables at every
rational time.
\item We next construct a family $(\mathsf U_t)_{t\geq0}$
of linear contractions such that, for every $x\in X\setminus N$,
$\mathbb P_x$-a.s.,
\[
    \left\langle
        \tilde\alpha_k(\mathsf X_t),
        \mathsf U_t\eta_j(\mathsf X_0)
    \right\rangle
    =Z_t^{k,j},
    \qquad t\geq0,\quad k,j\geq1.
\]
For $t\geq0$, fix a finite sum
\[
    v=\sum_j b_j\eta_j(\mathsf X_0),
    \qquad b_j\in\mathbb R,
\]
and formally set
\[
    \ell_{t,v}\!\left(
        \sum_k a_k\tilde\alpha_k(\mathsf X_t)
    \right)
    :=
    \sum_{k,j}a_kb_jZ_t^{k,j}
\]
for every finitely supported real sequence $(a_k)$.
We verify that this assignment is well-defined and bounded.
By quasi-continuity, after enlarging $N$ if necessary,
for every finitely supported rational sequence $(a_k)$,
the process
\[
    t\longmapsto
    \left|\sum_k a_k\tilde\alpha_k(\mathsf X_t)\right|
\]
has right-continuous paths $\mathbb P_x$-a.s. for every
$x\in X\setminus N$.
We also arrange that the families $(\tilde\alpha_k)$
and $(\eta_j)$ are fiberwise total outside $N$.
By countability, we may work on a common event of full
$\mathbb P_x$-measure on which these right-continuity
properties and all the preceding conclusions concerning
$Z^{k,j}$ hold, and the process avoids $N$.
For finitely supported rational sequences $(a_k)$ and
$(b_j)$, contractivity at rational times gives
\[
    \left|\sum_{k,j}a_kb_jZ_q^{k,j}\right|
    =
    \left|
        \left\langle
            \sum_k a_k\tilde\alpha_k(\mathsf X_q),
            \mathsf U_q v
        \right\rangle
    \right|
    \leq
    \left|\sum_k a_k\tilde\alpha_k(\mathsf X_q)\right|
    |v|.
\]
Letting $q\downarrow t$ through rational times yields
\[
    \left|\sum_{k,j}a_kb_jZ_t^{k,j}\right|
    \leq
    \left|\sum_k a_k\tilde\alpha_k(\mathsf X_t)\right|
    |v|,
    \qquad t\geq0.
\]
By continuity in the coefficients, this estimate also
holds for finitely supported real sequences $(a_k)$ and
$(b_j)$.
In particular, the defining value of $\ell_{t,v}$ is
independent of the representations of both
$v$ and $\sum_k a_k\tilde\alpha_k(\mathsf X_t)$.
Thus $\ell_{t,v}$ is a well-defined linear functional on
$\operatorname{Span}_{\mathbb R}
\{\tilde\alpha_k(\mathsf X_t)\mid k\geq1\}$ satisfying
\[
    |\ell_{t,v}(\alpha)|\leq |v|\,|\alpha|.
\]
Since this span is dense in
$\mathcal H_{\mathsf X_t}^{\breve{\mathcal E}}$,
$\ell_{t,v}$ extends uniquely to an element of
$(\mathcal H_{\mathsf X_t}^{\breve{\mathcal E}})^*$
with norm at most $|v|$.
Define
$
    \mathsf U_t v:=\ell_{t,v}.
$
This assignment is linear in $v$ and satisfies
$|\mathsf U_t v|\leq |v|$.
Since the linear span of $(\eta_j(\mathsf X_0))_{j\geq1}$
is dense in $(\mathcal H_{\mathsf X_0}^{\breve{\mathcal E}})^*$,
it extends uniquely to a linear contraction
\[
    \mathsf U_t\colon
    (\mathcal H_{\mathsf X_0}^{\breve{\mathcal E}})^*
    \longrightarrow (\mathcal H_{\mathsf X_t}^{\breve{\mathcal E}})^*,
\]
satisfying
\[
    \left\langle
        \tilde\alpha_k(\mathsf X_t),
        \mathsf U_t\eta_j(\mathsf X_0)
    \right\rangle
    =Z_t^{k,j},
    \qquad k,j\geq1.
\]
\item We prove right-continuity of the pairings.
Enlarge $N$ if necessary so that
\[
    |\tilde\alpha_k(y)|
    \leq \|\alpha_k\|_\infty,
    \qquad y\in X\setminus N,\quad k\geq1.
\]
For a bounded measurable section $\xi$ of $(\mathcal H^{\breve{\mathcal E}})^*$, set
\[
    \xi^{(n)}(y)
    :=
    \sum_{j=1}^n
    \langle\xi(y),\eta_j(y)\rangle\,\eta_j(y).
\]
Fix $x\in X\setminus N$ and work on the common
full-measure event constructed above.
Then for every $k,n\geq1$,
\[
    \left\langle
        \tilde\alpha_k(\mathsf X_t),
        \mathsf U_t\xi^{(n)}(\mathsf X_0)
    \right\rangle
    =
    \sum_{j=1}^n
    \langle\xi(\mathsf X_0),\eta_j(\mathsf X_0)\rangle
    Z_t^{k,j},
\]
so this process has c\`adl\`ag paths.
Moreover, contractivity gives, for every $T>0$,
\[
    \sup_{0\leq t\leq T}
    \left|
        \left\langle
            \tilde\alpha_k(\mathsf X_t),
            \mathsf U_t
            (\xi-\xi^{(n)})(\mathsf X_0)
        \right\rangle
    \right|
    \leq
    \|\alpha_k\|_\infty
    |(\xi-\xi^{(n)})(\mathsf X_0)|
    \longrightarrow0.
\]
Thus the lift paired with $\tilde\alpha_k(\mathsf X_t)$
is a uniform limit on bounded time intervals of
c\`adl\`ag functions, and hence is itself c\`adl\`ag so the fundamental pairings are c\`adl\`ag simultaneously
for all initial vectors.
Now fix $\alpha\in\D(\breve{\mathcal E})^\infty$.
By quasi-continuity and countability, there exists a properly
exceptional Borel set $N_\alpha\supset N$ such that, for every
$x\in X\setminus N_\alpha$, the processes
\[
    t\longmapsto
    \left|
        \tilde\alpha(\mathsf X_t)
        -\sum_k a_k\tilde\alpha_k(\mathsf X_t)
    \right|
\]
are right-continuous $\mathbb P_x$-a.s., simultaneously
for all finitely supported rational sequences $(a_k)$.
Fix $x\in X\setminus N_\alpha$ and work on the common
probability 1 event on which these properties and
the preceding conclusions hold.
For every $v\in(\mathcal H_x^{\breve{\mathcal E}})^*$,
contractivity and right-continuity of the fundamental
pairings give
\[
\begin{aligned}
    \varlimsup_{r\downarrow t}
    \left|
        \langle\tilde\alpha(\mathsf X_r),\mathsf U_rv\rangle
        -
        \langle\tilde\alpha(\mathsf X_t),\mathsf U_tv\rangle
    \right|
\leq
    2|v|
    \left|
        \tilde\alpha(\mathsf X_t)
        -\sum_k a_k\tilde\alpha_k(\mathsf X_t)
    \right|.
\end{aligned}
\]
Fiberwise density makes the right-hand side arbitrarily
small.
Thus the pairing with $\tilde\alpha(\mathsf X_t)$ is
right-continuous, simultaneously for all initial vectors. This proves property $2$.
\end{itemize}
\end{proof}
\begin{lem}\label{lem:transport-strong-multiplicativity-association}
The family $(\mathsf U_t)_{t\geq0}$ constructed in
Lemma~\ref{lem:transport-pairing-regularization} is a contractive strongly multiplicative transport
of $(\mathsf X_t)_{t\geq0}$ on $\mathcal H^{\breve{\mathcal E}}$
associated with $\breve{\mathcal E}$.
\end{lem}
\begin{proof}\;
\begin{itemize}
\item We first extend the conditional semigroup identity from rational times to bounded stopping times. Fix $u\in\D(\breve{\mathcal E})^\infty$ and a positive
rational number $h$.
For rational $q\geq0$, rational-time multiplicativity
and the semigroup representation give
\[
    \mathbb E_x\left[
        \mathsf U_{q+h}^*\tilde u(\mathsf X_{q+h})
        \,\middle|\,\mathcal F_{q+}^{\mathsf X}
    \right]
    =
    \mathsf U_q^*
    (\widetilde{\breve T_hu})(\mathsf X_q)
    \quad\mathbb P_x\text{-a.s.}
\]
for q.e.\ $x\in X$.
Let $\tau$ be a bounded
$(\mathcal F_{t+}^{\mathsf X})_{t\geq0}$-stopping time,
and set
$\tau_n\coloneqq2^{-n}\lceil2^n\tau\rceil$.
The preceding identity applies to $\tau_n$, since
$\tau_n$ takes finitely many rational values.
Test it against a bounded
$\mathcal F_{\tau+}^{\mathsf X}$-measurable random variable
and let $n\to\infty$.
Since $u,\breve T_hu\in\D(\breve{\mathcal E})^\infty$,
Lemma~\ref{lem:transport-pairing-regularization} and dominated convergence yield
\[
    \mathbb E_x\left[
        \mathsf U_{\tau+h}^*\tilde u(\mathsf X_{\tau+h})
        \,\middle|\,\mathcal F_{\tau+}^{\mathsf X}
    \right]
    =
    \mathsf U_\tau^*
    (\widetilde{\breve T_hu})(\mathsf X_\tau)
    \quad\mathbb P_x\text{-a.s.}
\]
The exceptional set may depend on $u$ and $h$, but can be chosen independently of $\tau$: only countably many rational-time identities and the pathwise right-continuity of the two pairings were used.

\item We test the identity against cylinder functions to obtain multiplicativity for rational time increments after a stopping time.
For a cylinder function
\[
    G=\prod_{i=1}^n\tilde f_i(\mathsf X_{r_i}),
    \qquad
    0<r_1<\cdots<r_n<h,
\]
with rational $r_i$ and $f_i\in\D(\mathcal E)^\infty$,
iteration of the preceding identity at $\tau+r_i$ gives
\[
    \mathbb E_x\left[
        (G\circ\theta_\tau)
        \mathsf U_{\tau+h}^*\tilde u(\mathsf X_{\tau+h})
        \,\middle|\,\mathcal F_{\tau+}^{\mathsf X}
    \right]
    =
    \mathsf U_\tau^*
    \mathbb E_{\mathsf X_\tau}\left[
        G\mathsf U_h^*\tilde u(\mathsf X_h)
    \right].
\]
Indeed, all intermediate sections remain in
$\D(\breve{\mathcal E})^\infty$, and the expectation on
the right has the same iterated semigroup expression.
Multiplying by a bounded
$\mathcal F_{\tau+}^{\mathsf X}$-measurable random variable $F$,
taking expectations, and using the strong Markov property
at $\tau$, we obtain
\[
\begin{aligned}
    \mathbb E_x\left[
        F(G\circ\theta_\tau)
        \mathsf U_{\tau+h}^*\tilde u(\mathsf X_{\tau+h})
    \right]
    &=
    \mathbb E_x\left[
        F\mathsf U_\tau^*
        \mathbb E_{\mathsf X_\tau}\left[
            G\mathsf U_h^*\tilde u(\mathsf X_h)
        \right]
    \right]
    \\
    &=
    \mathbb E_x\left[
        F(G\circ\theta_\tau)
        \mathsf U_\tau^*(\mathsf U_h^*\circ\theta_\tau)
        \tilde u(\mathsf X_{\tau+h})
    \right].
\end{aligned}
\]
The cylinder-test argument used to prove property $3$
of Lemma~\ref{lem:transport-rational-times} now applies to this equality.
Factors at $h$ are absorbed into $u$, and factors after
$h$ are first conditioned at $\tau+h$ and then absorbed
in the same way.
Using countable determining scalar functions and a
countable quasi-fundamental family of terminal sections,
include the exceptional sets for their finite products,
rational-time semigroup expressions, and right-continuity
assertions in a common properly exceptional Borel set $N$.
This set is independent of $\tau$ and rational $h$.
By the same monotone-class argument and fiberwise totality,
for every $x\in X\setminus N$,
\[
    \mathsf U_{\tau+h}
    =
    (\mathsf U_h\circ\theta_\tau)\mathsf U_\tau
    \quad\mathbb P_x\text{-a.s.},
    \qquad h\in\mathbb Q_+.
\]
The case $h=0$ follows from $\mathsf U_0=I$.

\item We extend multiplicativity to every real time increment by right-continuity of the scalar pairings.
For $h\geq0$, choose rational
$h_n\downarrow h$.
By Lemma~\ref{lem:transport-pairing-regularization} and its proof, the pairings with
the chosen terminal sections are right-continuous
simultaneously for all initial vectors.
The strong Markov property gives the same conclusion
along the shifted path $\theta_\tau$.
We may therefore pass to the limit in the scalar pairings
of both sides of
$\mathsf U_{\tau+h_n}
=(\mathsf U_{h_n}\circ\theta_\tau)\mathsf U_\tau$.
Fiberwise totality gives the desired identity at $h$.
This proves strong multiplicativity.

\item We show that the transport represents the bundle semigroup on $\L^2(m;\mathcal H^{\breve{\mathcal E}})$.
Contractivity gives
$
    |\breve p_t^{\mathsf U}\alpha|
    \leq p_t^{\mathsf X}|\alpha|
$
for every bounded measurable section $\alpha$.
Consequently, each $\breve p_t^{\mathsf U}$ descends to
a contraction on $\L^2(m;\mathcal H^{\breve{\mathcal E}})$.
Fix $u\in\D(\breve{\mathcal E})^\infty$ and $t\geq0$,
and choose rational $q_n\downarrow t$.
For every bounded integrable measurable section $\xi$
of $(\mathcal H^{\breve{\mathcal E}})^*$, Lemma~\ref{lem:transport-pairing-regularization}
and the bound
\[
    \left|
        \left\langle
            \tilde u(\mathsf X_r),
            \mathsf U_r\xi(\mathsf X_0)
        \right\rangle
    \right|
    \leq \|u\|_\infty|\xi(\mathsf X_0)|
\]
allow us to apply dominated convergence.
Using the rational-time semigroup representation and
the strong continuity of $(\breve T_t)_{t\geq0}$, we obtain
\begin{align*}
    \int_X
        \left\langle
            \breve p_t^{\mathsf U}\tilde u(x),\xi(x)
        \right\rangle\,dm(x)
    &=
    \lim_{n\to\infty}
    \int_X
        \left\langle
            (\breve T_{q_n}u)(x),\xi(x)
        \right\rangle\,dm(x)
    \\
    &=
    \int_X
        \left\langle
            (\breve T_tu)(x),\xi(x)
        \right\rangle\,dm(x).
\end{align*}
Hence
$\breve p_t^{\mathsf U}\tilde u=\breve T_tu$
$m$-a.e.
Density of $\D(\breve{\mathcal E})^\infty$ in
$\L^2(m;\mathcal H^{\breve{\mathcal E}})$ and contractivity
of both operators extend this equality to every bounded
measurable section belonging to $\L^2(m;\mathcal H^{\breve{\mathcal E}})$.
Thus $(\mathsf U_t)_{t\geq0}$ is associated with
$\breve{\mathcal E}$.
\end{itemize}
\end{proof}
\begin{lem}\label{lem:transport-quasi-right-continuity}
The strongly multiplicative transport $(\mathsf U_t)_{t\geq0}$
constructed in Lemma~\ref{lem:transport-strong-multiplicativity-association} is quasi-right and properly associated
with $\breve{\mathcal E}$.
\end{lem}
\begin{proof}\;
\begin{itemize}
\item We construct a right-continuous candidate for the norm of the transported section from its values at rational times.
Fix a measurable section $\xi$ of
$(\mathcal H^{\breve{\mathcal E}})^*$.
Choose a properly exceptional Borel set $N$ outside
which strong multiplicativity holds.
By multiplicativity and contractivity, for every
$x\in X\setminus N$,
\[
    |\mathsf U_q\xi(\mathsf X_0)|
    \leq|\mathsf U_r\xi(\mathsf X_0)|,
    \qquad r,q\in\mathbb Q_+,\quad r\leq q,
\]
simultaneously $\mathbb P_x$-a.s.
Thus the rational right limits
\[
    R_t\coloneqq
    \lim_{\substack{q\downarrow t\\q\in\mathbb Q_+}}
    |\mathsf U_q\xi(\mathsf X_0)|
\]
define, up to an evanescent set, an adapted,
right-continuous, nonincreasing process under each
$\mathbb P_x$, $x\in X\setminus N$.
Adaptedness follows from the right-continuity of the
filtration.
It suffices to show that
$(|\mathsf U_t\xi(\mathsf X_0)|)_{t\geq0}$ and $R$
are indistinguishable under $\mathbb P_x$ for q.e.\ $x$.

\item We show that the norm process is optional and is bounded above by its right-limit process.
For q.e.\ $x\in X$, we claim that
\[
    |\mathsf U_t\xi(\mathsf X_0)|\leq R_t,
    \qquad t\geq0,
\]
simultaneously $\mathbb P_x$-a.s.
Choose a countable family
$(\alpha_k)_{k\geq1}\subset\D(\breve{\mathcal E})^\infty$
whose quasi-continuous representatives form a
quasi-fundamental family, and include all finite rational
linear combinations in this family.
By Lemma~\ref{lem:transport-pairing-regularization}, applied to bounded truncations
of $\xi$, and quasi-continuity, we may enlarge $N$ so that,
for every $x\in X\setminus N$,
the processes
\[
    t\longmapsto
    \left\langle
        \tilde\alpha_k(\mathsf X_t),
        \mathsf U_t\xi(\mathsf X_0)
    \right\rangle,
    \qquad
    t\longmapsto|\tilde\alpha_k(\mathsf X_t)|
\]
have right-continuous paths simultaneously for all $k$,
$\mathbb P_x$-a.s.
We also arrange that this family is fiberwise total
outside $N$.
Fix $x\in X\setminus N$ and work with the usual
$\mathbb P_x$-augmentation of
$(\mathcal F_{t+}^{\mathsf X})_{t\geq0}$.
Every bounded stopping time of this augmentation has an
almost surely equal stopping time for the original
right-continuous filtration, so strong multiplicativity
remains available.
Fiberwise totality gives
\[
    |\mathsf U_t\xi(\mathsf X_0)|
    =
    \sup_{k,j\geq1}
    \frac{
        |\langle\tilde\alpha_k(\mathsf X_t),
            \mathsf U_t\xi(\mathsf X_0)\rangle|
    }{
        |\tilde\alpha_k(\mathsf X_t)|+j^{-1}
    },
    \qquad t\geq0,
\]
simultaneously $\mathbb P_x$-a.s.
Each quotient is adapted and has right-continuous paths.
Consequently, the norm process is optional, up to an
evanescent set, and has right lower semicontinuous paths.
In particular,
\[
    |\mathsf U_t\xi(\mathsf X_0)|
    \leq
    \varliminf_{s\downarrow t}
    |\mathsf U_s\xi(\mathsf X_0)|
    \leq R_t,
    \qquad t\geq0,
\]
simultaneously $\mathbb P_x$-a.s.

\item We prove the reverse inequality at bounded stopping times and use the optional section theorem to deduce quasi-right-continuity and proper association.
Fix a bounded stopping
time $\tau$.
Strong multiplicativity and contractivity give
\[
    |\mathsf U_{\tau+h}\xi(\mathsf X_0)|
    =
    |(\mathsf U_h\circ\theta_\tau)
        \mathsf U_\tau\xi(\mathsf X_0)|
    \leq|\mathsf U_\tau\xi(\mathsf X_0)|,
    \qquad h\in\mathbb Q_+,
\]
simultaneously $\mathbb P_x$-a.s.
By right lower semicontinuity and density of the rational
numbers, the inequality holds simultaneously for all
real $h\geq0$.
Letting $q\downarrow\tau$ through rational times therefore
gives
\[
    R_\tau\leq|\mathsf U_\tau\xi(\mathsf X_0)|
    \qquad\mathbb P_x\text{-a.s.}
\]
Combining the two inequalities yields
$|\mathsf U_\tau\xi(\mathsf X_0)|=R_\tau$ almost surely
for every bounded stopping time $\tau$.
Since both processes are optional, the optional section
theorem \cite[Theorem~IV.84]{DellacherieMeyer1978} implies that they are indistinguishable.
Hence $t\mapsto|\mathsf U_t\xi(\mathsf X_0)|$
has right-continuous paths $\mathbb P_x$-a.s.
Together with the right-continuity of the fundamental
pairings, this proves that
$(\mathsf U_t\xi(\mathsf X_0))_{t\geq0}$
is quasi-right-continuous.
Since $\xi$ was arbitrary, $(\mathsf U_t)_{t\geq0}$
is a quasi-right transport.
It is already associated with $\breve{\mathcal E}$,
so Lemma~\ref{lem:association-is-proper} implies that it is properly
associated with $\breve{\mathcal E}$.
\end{itemize}
\end{proof}
\end{proof}

\begin{ex}\label{ex:transport-no-left-limit}
Left limits of the associated transport need not exist even when the base process has continuous paths. As in the counterexample in Section~\ref{subsec:multiplicative-transports}, let $\mathsf X$ be reflecting Brownian motion on $X=[0,\infty)$ with Lebesgue measure. Set $\widetilde{\mathcal K}_x=\mathbb R$ for $x>0$ and $\widetilde{\mathcal K}_0=\{0\}$. The restriction of the Neumann energy to this quasi-subbundle is the Kato form
\[
    \breve{\mathcal E}(\alpha)=\frac12\int_0^\infty|\alpha'(x)|^2\,dx,
    \qquad
    \D(\breve{\mathcal E})=\{\alpha\in W^{1,2}([0,\infty)):\alpha(0)=0\},
\]
with quasi-realization $\widetilde{\mathcal K}$. For $x>0$, its associated transport, viewed in the ambient line bundle, is
\[
    \mathsf U_t^{\mathcal K}v=v\,\mathbf1_{\{t<\sigma_0\}},
    \qquad
    \sigma_0:=\inf\{t\geq0:\mathsf X_t=0\}.
\]
For $v\ne0$, the norm is $|v|$ before $\sigma_0$, whereas the fiber at $\mathsf X_{\sigma_0}=0$ is $\{0\}$. Thus the lift has no left limit at $\sigma_0$ in $\widetilde{\mathcal K}^*$, although its ambient scalar left limit exists. This is also an example of the restriction transport constructed in Section~\ref{subsec:restriction-transports}.
\end{ex}
\clearpage
\section{Basic operations on multiplicative transport}\label{sec:transport-operations}
Throughout this section, we fix a Kato form
$\breve{\mathcal E}$ on $(X,m,\mathcal E)$ and a properly associated
contractive quasi-right transport $(\mathsf U_t)_{t\geq0}$.
\begin{thm}[Resolvent formula]\label{thm:vector-revuz-resolvent}
Let $\mu$ be a smooth measure with corresponding PCAF $A^\mu$. Let $F\subset X$ be closed and let $\tau_F$ be the first exit time
of $\mathsf X$ from $F$, with $\tau_F\le\zeta$.
Suppose that $\alpha\mu$ is of finite energy integral. Then, for every
$\lambda>0$,
\[
    \widetilde{\breve G_{\lambda,F}(\alpha\mu)}(x)
    =
    \mathbb E_x\!\left[
        \int_0^{\tau_F} e^{-\lambda t}
        \mathsf U_t^*\alpha(\mathsf X_t)\,dA_t^\mu
    \right]
    \qquad \text{for q.e.\ } x\in X.
\]
\end{thm}
\begin{proof}
\;
\begin{itemize}
\item We first show that, for every $\lambda>0$ and every
$g\in\L^2(m,H^{\breve{\mathcal E}})$,
\[
    \mathbb E_x\!\left[
        \int_0^{\tau_F} e^{-\lambda t}
        \mathsf U_t^*g(\mathsf X_t)\,dt
    \right]
    =
    \widetilde{\breve G_{\lambda,F}g}(x)
    \qquad \text{for q.e.\ } x\in X.
\]
For each integer $k\ge1$, define
\[
    u_k(x):=
    \mathbb E_x\!\left[
        \int_0^\zeta e^{-\lambda t}
        \exp\!\left(
            -k\int_0^t \mathbf 1_{X\setminus F}(\mathsf X_r)\,dr
        \right)
        \mathsf U_t^*g(\mathsf X_t)\,dt
    \right].
\]
By contractivity, $|u_k|\le\widetilde{G_\lambda|g|}$ q.e., so $u_k$
is well defined q.e.\ and belongs to $\L^2(m,H^{\breve{\mathcal E}})$.
Using
\[
    1-\exp\!\left(
        -k\int_0^t \mathbf 1_{X\setminus F}(\mathsf X_r)\,dr
    \right)
    =
    k\int_0^t \mathbf 1_{X\setminus F}(\mathsf X_s)
    \exp\!\left(
        -k\int_s^t \mathbf 1_{X\setminus F}(\mathsf X_r)\,dr
    \right)\,ds
\]
and conditioning at time $s$, multiplicativity and the Markov property
give
\[
    \widetilde{\breve G_\lambda g}(x)-u_k(x)
    =
    k\,\widetilde{\breve G_\lambda(\mathbf 1_{X\setminus F}u_k)}(x)
    \qquad\text{for q.e.\ }x\in X.
\]
Here the full-space resolvent representation extends from bounded
sections to $\L^2$-sections by truncation and domination by the scalar
resolvent. Thus $u_k$ is the quasi-continuous representative of an
element of $\D(\breve{\mathcal E})$ satisfying
\[
    \breve{\mathcal E}_\lambda(u_k,\eta)
    +k\int_{X\setminus F}\langle u_k,\eta\rangle\,dm
    =(g,\eta)_2
    \qquad\text{for every }\eta\in\D(\breve{\mathcal E}).
\]
We next show that $u_k\to u:=\breve G_{\lambda,F}g$ in the form norm.
Taking $\eta=u_k$ gives
\[
    \breve{\mathcal E}_\lambda(u_k)
    +k\|\mathbf 1_{X\setminus F}u_k\|_2^2
    =(g,u_k)_2.
\]
Hence $(u_k)$ is bounded in the form norm and
$\|\mathbf 1_{X\setminus F}u_k\|_2\to0$. Every weak limit in the form domain therefore vanishes $m$-a.e.\ on the open set $X\setminus F$,
so its quasi-continuous representative vanishes there q.e.
Taking $\eta\in\D(\breve{\mathcal E})_F$ in the variational identity
eliminates the integral over $X\setminus F$ and identifies every weak
limit with $u$. Thus $u_k$ converges weakly to $u$ in the form domain.
Moreover,
\[
    \varlimsup_{k\to\infty}\breve{\mathcal E}_\lambda(u_k)
    \le \lim_{k\to\infty}(g,u_k)_2
    =(g,u)_2
    =\breve{\mathcal E}_\lambda(u).
\]
It follows that $u_k\to u$ in the form norm and hence, along a
subsequence, $u_k(x)\to\tilde u(x)$ for q.e.\ $x$.
Finally, we pass to the limit in the probabilistic expression for
$u_k$. For $t<\tau_F$, the exponential factor equals one. If
$\tau_F<t<\zeta$, the path visits the open set $X\setminus F$ before
$t$ and, by right continuity, spends a positive amount of time there
before $t$. Thus the exponential factor tends to zero. The single
time $t=\tau_F$ has no effect on the $dt$-integral. Dominated
convergence, using the scalar potential $\widetilde{G_\lambda|g|}$,
therefore gives
\[
    \widetilde{\breve G_{\lambda,F}g}(x)
    =
    \mathbb E_x\!\left[
        \int_0^{\tau_F}e^{-\lambda t}
        \mathsf U_t^*g(\mathsf X_t)\,dt
    \right]
    \qquad\text{for q.e.\ }x\in X.
\]
\item We next prove 
\[
    \widetilde{\breve G_{\lambda,F}(\alpha\mu)}(x)
    =\mathbb E_x\!\left[
        \int_0^{\tau_F} e^{-\lambda t}
        \mathsf U_t^*\alpha(\mathsf X_t)\,dA_t^\mu
    \right]
    \qquad \text{for } m\text{-a.e. } x\in X.
\]
It suffices to show that, for every bounded $\xi\in L^2(m,H^{\breve{\mathcal E}})$,
\[
    \mathbb E_m\!\left[
        \int_0^{\tau_F} e^{-\lambda t}
        \bigl\langle
            \mathsf U_t\xi(\mathsf X_0),\alpha(\mathsf X_t)
        \bigr\rangle\,dA_t^\mu
    \right]
    =\int_X \bigl\langle\widetilde{\breve G_{\lambda,F}\xi},\alpha\bigr\rangle\,d\mu . \tag{{$*$}}
\]
Suppose first that $\mu$ is finite and that $\alpha$ is a bounded
quasi-continuous section. To approximate $dA_t^\mu$ by
ordinary time integrals, put
\[
    h_\varepsilon(x)
    :=
    \frac{1}{\varepsilon}\mathbb E_x[A_\varepsilon^\mu].
\]
By the Markov property, association, and symmetry,
\[
\begin{aligned}
    &\mathbb E_m\!\left[
        \int_0^{\tau_F} e^{-\lambda t}
        \bigl\langle
            \mathsf U_t\xi(\mathsf X_0),\alpha(\mathsf X_t)
        \bigr\rangle
        \frac{A_{t+\varepsilon}^\mu-A_t^\mu}{\varepsilon}\,dt
    \right] =
    \int_X 
    \bigl\langle\widetilde{\breve G_{\lambda,F}\xi},\alpha\bigr\rangle
    h_\varepsilon\,dm .
\end{aligned}
\]
The left-hand side converges to the left-hand side of $(*)$.
(Indeed, after interchanging the time integrals, the integrand against
$dA_s^\mu$ is the left average of the discounted pairing multiplied by
$\mathbf 1_{\{t<\tau_F\}}$ over $[(s-\varepsilon)^+,s]$. Since $\alpha$ is quasi-continuous, the
quasi-right continuity of $\mathsf U$ gives right continuity of the
pairing. Its discontinuity times are therefore at most countable.
At every continuity time $s>0$, the left average converges to the value at $s$, and the continuity of $A^\mu$ ensures that the remaining times have zero $dA^\mu$-measure.)
For the right-hand side, the scalar Revuz formula \cite[Section~5.1]{FukushimaOshimaTakeda2011} gives, for every
bounded quasi-continuous $f$,
\[
    \int_X h_\varepsilon \langle\widetilde{\breve G_{\lambda,F}\xi},\alpha\rangle\,dm
    =
    \frac{1}{\varepsilon}\int_0^\varepsilon
    \int_X P_s (\langle\widetilde{\breve G_{\lambda,F}\xi},\alpha\rangle)\,d\mu\,ds
    \longrightarrow
    \int_X \langle\widetilde{\breve G_{\lambda,F}\xi},\alpha\rangle\,d\mu .
\]
This proves $(*)$ in this case. We now prove $(*)$ for a measurable section $\alpha$. Put
$\nu:=|\alpha|\mu$ and define $\beta:=\alpha/|\alpha|$ on
$\{|\alpha|>0\}$ and $\beta:=0$ elsewhere. Since
$dA_t^\nu=|\alpha|(\mathsf X_t)\,dA_t^\mu$, it suffices to prove
$(*)$ with $\beta$ and $\nu$ in place of $\alpha$ and $\mu$.
Choose measurable sets $F_n\uparrow X$ with $\nu(F_n)<\infty$ and a
sequence of bounded quasi-continuous sections $(\eta_j)_{j\ge1}$ whose
values are dense in the closed unit ball of each fiber q.e., with
$|\eta_j|\le1$. On $F_n$, let $\beta_n$ take the value of a section
among $\eta_1,\ldots,\eta_n$ closest to $\beta$, choosing the least
index in case of a tie, and set $\beta_n:=0$ outside $F_n$. Then
$|\beta_n|\le1$ and $\beta_n\to\beta$ $\nu$-a.e. For each measurable
set $E\subset F_n$ on which a given $\eta_j$ is selected, the
preceding case applies to $\eta_j$ and the finite measure
$\mathbf 1_E\nu$. Summing over these sets gives $(*)$ for $\beta_n$
and $\nu$.
By the contractivity of $\mathsf U$ and the scalar Revuz formula on
the left-hand side, and by resolvent domination on the right-hand
side, the error on either side is bounded by
\[
    \int_X \widetilde{G_\lambda|\xi|}\,|\beta_n-\beta|\,d\nu .
\]
Since $\nu$ is of finite energy integral,
\[
    \int_X \widetilde{G_\lambda|\xi|}\,d\nu
    =\bigl(G_\lambda\nu,|\xi|\bigr)_2<\infty .
\]
Dominated convergence therefore gives $(*)$ for $\beta$ and $\nu$,
and hence for $\alpha$ and $\mu$.

\item
It remains to identify the quasi-continuous representative.
For every integer $n\ge1$, multiplicativity, the Markov property,
and proper association give, for q.e.\ $x$,
\[
\begin{aligned}
    \mathbb E_x\!\left[
        \int_0^{\tau_F} e^{-\lambda t}
        \mathsf U_t^*\alpha(\mathsf X_t)\,dA_t^\mu\right]
    &-
    n\,\widetilde{
        \breve G_{\lambda+n,F}\bigl(\breve G_{\lambda,F}(\alpha\mu)\bigr)
    }(x)=
    \mathbb E_x\!\left[
        \int_0^{\tau_F} e^{-(\lambda+n)t}
        \mathsf U_t^*\alpha(\mathsf X_t)\,dA_t^\mu
    \right].
\end{aligned}
\]
The right-hand side tends to zero q.e.\ by the scalar potential bound and the continuity of $A^\mu$ at zero. On the other hand,
\[
    n\breve G_{\lambda+n,F}\bigl(\breve G_{\lambda,F}(\alpha\mu)\bigr)
    \longrightarrow
    \breve G_{\lambda,F}(\alpha\mu)
    \qquad\text{in the form norm}.
\]
Along a subsequence the quasi-continuous representatives therefore converge q.e. The claim follows.
\end{itemize}
\end{proof}

\subsection{Damped transport}\label{subsec:damped-transports}
Let $R\mu$ be a smooth potential on $\mathcal H^{\breve{\mathcal E}}$ as in Section~\ref{subsec:potential-perturbation}, and assume that $\breve{\mathcal E}^{R\mu}$ is a bundle form. The construction below is an operator-valued version of the scalar Feynman--Kac construction; see \cite[Section~6.1]{FukushimaOshimaTakeda2011} for killing by additive functionals.
\begin{thm}[Damped transport]\label{thm:damped-transport-existence}
There exists a family of bounded linear maps
\[
    \mathsf U_t^{R\mu}(\omega)\colon
    (\mathcal H_{\mathsf X_0(\omega)}^{\breve{\mathcal E}})^*
    \longrightarrow
    (\mathcal H_{\mathsf X_t(\omega)}^{\breve{\mathcal E}})^*,
\]
with $\mathsf U_t^{R\mu}=0$ for $t\geq\zeta$, such that,
for q.e.\ $x\in X$, $\mathbb P_x$-a.s.,
this family is locally bounded in operator norm on
$[0,\zeta)$ and satisfies
\[
    \mathsf U_t^{R\mu}
    =
    \mathsf U_t
    -
    \int_0^t
    \mathsf U_{s,t}\,
    R(\mathsf X_s)^*\mathsf U_s^{R\mu}\,dA_s^\mu
    \qquad\text{for all }0\leq t<\zeta,
\]
where $A^\mu$ is the PCAF associated with $\mu$ by the Revuz correspondence \cite[Section~5.1]{FukushimaOshimaTakeda2011}.
The integral is understood in the strong operator sense.
Such a family is unique up to
$\mathbb P_x$-indistinguishability for q.e.\ $x\in X$.
\end{thm}
\begin{proof}
Define
\[
    a_t:=\int_0^t
    \|R(\mathsf X_s)\|_{\mathrm{op}}\,dA_s^\mu,
    \qquad 0\leq t<\zeta.
\]
Since $\|R\|_{\mathrm{op}}\mu$ is smooth, $a$ is the PCAF
with this Revuz measure. In particular, for q.e.\ $x\in X$,
it is finite and continuous on $[0,\zeta)$, with $a_0=0$,
$\mathbb P_x$-a.s.
We choose a contractive version of $\mathsf U$ by setting
$\mathsf U_t=0$ wherever its operator norm exceeds one;
this does not change its indistinguishability class.
Thus $\|\mathsf U_{s,t}\|_{\mathrm{op}}\leq1$.
Fix such a path and define the Picard iterates by
\[
    \mathsf U_t^{R\mu,(0)}:=\mathsf U_t,
    \qquad
    \mathsf U_t^{R\mu,(n+1)}
    :=
    \mathsf U_t-
    \int_0^t
    \mathsf U_{s,t}R(\mathsf X_s)^*
    \mathsf U_s^{R\mu,(n)}\,dA_s^\mu.
\]
After application to a vector in the initial fiber,
these are Bochner integrals in the fiber at time $t$.
Strong measurability follows in measurable fiber coordinates,
and the following estimates ensure integrability.
Using contractivity and the continuity of $a$, induction gives
\[
    \|\mathsf U_t^{R\mu,(n+1)}-\mathsf U_t^{R\mu,(n)}\|_{\mathrm{op}}
    \leq \frac{a_t^{n+1}}{(n+1)!},
    \qquad n\geq0,\quad t<\zeta.
\]
Indeed, the first difference is bounded by $a_t$, and the
induction step follows from
\[
    \int_0^t\frac{a_s^n}{n!}\,da_s
    =\frac{a_t^{n+1}}{(n+1)!}.
\]
Consequently, for every $T<\zeta$, the iterates converge
uniformly in operator norm on $[0,T]$ to a family $\mathsf U^{R\mu}$,
and
\[
    \|\mathsf U_t^{R\mu}\|_{\mathrm{op}}\leq e^{a_t},
    \qquad t<\zeta.
\]
Dominated convergence in the strong operator integrals
shows that $\mathsf U^{R\mu}$ satisfies the stated equation
for every $t<\zeta$.
The construction is measurable in fiber coordinates;
we extend $\mathsf U^{R\mu}$ by zero after the lifetime.
For uniqueness, let $\widehat{\mathsf U}^{R\mu}$ be another such family
and fix $T<\zeta$ on a path where both equations hold.
Set
\[
    M_T:=
    \sup_{0\leq s\leq T}
    \|\mathsf U_s^{R\mu}-\widehat{\mathsf U}_s^{R\mu}\|_{\mathrm{op}}
    <\infty.
\]
Subtracting the equations and iterating the resulting estimate yields
\[
    \|\mathsf U_t^{R\mu}-\widehat{\mathsf U}_t^{R\mu}\|_{\mathrm{op}}
    \leq M_T\frac{a_t^n}{n!},
    \qquad 0\leq t\leq T,\quad n\geq0.
\]
Letting $n\to\infty$ proves equality for all $t\leq T$.
As $T<\zeta$ is arbitrary and both families vanish after the lifetime,
the required indistinguishability follows.
\end{proof}
\begin{dfn}
The family $(\mathsf U_t^{R\mu})_{t\geq0}$ given by Theorem~\ref{thm:damped-transport-existence}
is called the damped transport of $(\mathsf U_t)_{t\geq 0}$ associated with $R\mu$.
\end{dfn}
The defining equation and Gronwall's inequality yield, for q.e.\ $x\in X$,
$\mathbb P_x$-a.s.,
\[
    \|\mathsf U_t^{R\mu}\|_{\mathrm{op}}
    \leq
    \exp\left(
        \int_0^t\|R(\mathsf X_s)\|_{\mathrm{op}}\,dA_s^\mu
    \right),
    \qquad 0\leq t<\zeta.
\]
In the rest of this subsection, we assume that
\[
    \sup_{x\in X}
    \mathbb E_x\left[
        \sup_{0\leq t\leq T}\|\mathsf U_t^{R\mu}\|_{\mathrm{op}}
    \right]
    <\infty
    \qquad\text{for every }T>0.
\]
\begin{rem}
This assumption holds if $\nu:=\|R\|_{\mathrm{op}}\mu$ belongs to the Kato class. Indeed, the preceding estimate bounds $\sup_{0\leq t\leq T}\|\mathsf U_t^{R\mu}\|_{\mathrm{op}}$ by $e^{A_T^\nu}$, and Khasminskii's lemma gives $\sup_{x\in X\setminus N}\mathbb E_x[e^{A_T^\nu}]<\infty$ for every $T>0$ outside a properly exceptional set $N$ \cite[Lemma~2.1(b)]{Ying1997}. On paths starting in $N$, we take $\mathsf U_t^{R\mu}=0$ for $t>0$.
\end{rem}
\begin{prp}\label{prp:damped-transport-identities}
The family $(\mathsf U_t^{R\mu})_{t\geq0}$ is a quasi-right transport
of $(\mathsf X_t)_{t\geq0}$ on $\mathcal H^{\breve{\mathcal E}}$.
Moreover, for q.e.\ $x\in X$, $\mathbb P_x$-a.s.,
\[
    \mathsf U_t^{R\mu}
    =
    \mathsf U_t
    -
    \int_0^t
    \mathsf U_{s,t}^{R\mu}R(\mathsf X_s)^*\mathsf U_s\,dA_s^\mu,
    \qquad 0\leq t<\zeta.
\]
\end{prp}
\begin{proof}\;
\begin{itemize}
    \item We first prove multiplicativity.
Fix $t\geq0$.
For q.e.\ $x\in X$, we work $\mathbb P_x$-a.s.\ on $\{t<\zeta\}$.
Splitting the defining integral at $t$ and using the
multiplicativity of $\mathsf U$, we obtain
\[
\begin{aligned}
    \mathsf U_{t+h}^{R\mu}
    &=
    \mathsf U_{t,t+h}
    \left(
        \mathsf U_t
        -
        \int_0^t
        \mathsf U_{s,t}R(\mathsf X_s)^*\mathsf U_s^{R\mu}\,dA_s^\mu
    \right)
    -
    \int_t^{t+h}
    \mathsf U_{s,t+h}R(\mathsf X_s)^*\mathsf U_s^{R\mu}\,dA_s^\mu\\
    &=
    \mathsf U_{t,t+h}\mathsf U_t^{R\mu}
    -
    \int_0^h
    \mathsf U_{t+s,t+h}R(\mathsf X_{t+s})^*
    \mathsf U_{t+s}^{R\mu}\,dA_{t+s}^\mu,
    \qquad 0\leq h<\zeta-t.
\end{aligned}
\]
On the other hand, shifting the defining equation by $\theta_t$, using the additivity of $A^\mu$, and composing on the right with $\mathsf U_t^{R\mu}$, we obtain
\[
    (\mathsf U_h^{R\mu}\circ\theta_t)\mathsf U_t^{R\mu}
    =
    \mathsf U_{t,t+h}\mathsf U_t^{R\mu}
    -
    \int_0^h
    \mathsf U_{t+s,t+h}R(\mathsf X_{t+s})^*
    (\mathsf U_s^{R\mu}\circ\theta_t)\mathsf U_t^{R\mu}\,dA_{t+s}^\mu.
\]
Thus both
$h\mapsto \mathsf U_{t+h}^{R\mu}$ and $h\mapsto (\mathsf U_h^{R\mu}\circ\theta_t)\mathsf U_t^{R\mu}$
satisfy the same linear integral equation with initial value $\mathsf U_t^{R\mu}$.
The uniqueness argument applies to this initial value as well, and hence
\[
    \mathsf U_{t+h}^{R\mu}
    =
    (\mathsf U_h^{R\mu}\circ\theta_t)\mathsf U_t^{R\mu},
    \qquad t+h<\zeta.
\]
Since both sides vanish after the lifetime, this proves
multiplicativity.
\item The quasi-right-continuity follows from the quasi-right continuity of $(\mathsf U_t)_{t\geq0}$ and the locally uniform converge of Picard iterates $(\mathsf U_t^{R\mu,(n)})_{t\geq 0}$ in operator norm before $\zeta$.
\item Finally, we prove the stated integral identity.
For $t<\zeta$, set
\[
    \widehat{\mathsf U}_t^{R\mu}
    :=
    \mathsf U_t
    -
    \int_0^t
    \mathsf U_{s,t}^{R\mu}R(\mathsf X_s)^*\mathsf U_s\,dA_s^\mu.
\]
The local integrability of
$\|R(\mathsf X_s)\|_{\mathrm{op}}\,dA_s^\mu$
and the Gronwall bound for $\mathsf U_{s,t}^{R\mu}$
ensure that $\widehat{\mathsf U}^{R\mu}$ is locally bounded in operator norm. By using Fubini theorem and 
substituting the defining equation shifted by $\theta_s$, we obtain
\[
\begin{aligned}
    \widehat{\mathsf U}_t^{R\mu}
    ={}&
    \mathsf U_t
    -
    \int_0^t
    \mathsf U_{s,t}R(\mathsf X_s)^*\mathsf U_s\,dA_s^\mu+
    \int_0^t
    \mathsf U_{r,t}R(\mathsf X_r)^*
    \left(
        \int_0^r
        \mathsf U_{s,r}^{R\mu}R(\mathsf X_s)^*\mathsf U_s\,dA_s^\mu
    \right)dA_r^\mu\\
    ={}&
    \mathsf U_t
    -
    \int_0^t
    \mathsf U_{r,t}R(\mathsf X_r)^*
    \widehat{\mathsf U}_r^{R\mu}\,dA_r^\mu.
\end{aligned}
\]
Thus $\widehat{\mathsf U}^{R\mu}$ satisfies the same integral equation
as $\mathsf U^{R\mu}$.
Uniqueness gives $\widehat{\mathsf U}^{R\mu}=\mathsf U^{R\mu}$, proving the identity.
\end{itemize}
\end{proof}

\begin{thm}[Association of damped transport]\label{thm:damped-transport-association}
The damped transport
\[
    \mathsf U_t^{R\mu}(\omega)\colon
    (\mathcal H_{\mathsf X_0(\omega)}^{\breve{\mathcal E}})^*
    \longrightarrow
    (\mathcal H_{\mathsf X_t(\omega)}^{\breve{\mathcal E}})^*,
\]
is properly associated with $\breve{\mathcal E}^{R\mu}$.
\end{thm}
\begin{proof}\;
\begin{itemize}
    \item Let $F\subset X$ be closed, let $\tau_F$ be its first exit time, and assume that $\mathbf 1_F\|R\|_{\mathrm{op}}\mu$ is of finite energy integral.
    For a bounded measurable section $\alpha\in\L^2(m;\mathcal H)$ and $\lambda>0$ sufficiently large, set
\[
    \omega_F(x):=
    \mathbb E_x\!\left[
        \int_0^{\tau_F} e^{-\lambda t}
        (\mathsf U_t^{R\mu})^*\alpha(\mathsf X_t)\,dt
    \right].
\]
We claim that $\omega_F$ is the quasi-continuous representative of the unique $u\in\D(\breve{\mathcal E}^{R\mu})_F$ satisfying
\[
    \breve{\mathcal E}^{R\mu}_\lambda(u,\beta)
    =(\alpha,\beta)_2
    \qquad\text{for every }\beta\in\D(\breve{\mathcal E}^{R\mu})_F.
\]
To this end, it suffices to prove the stopped resolvent identity
\[
    \omega_F
    =
    \widetilde{\breve G_{\lambda,F}(\alpha m)}
    -
    \widetilde{\breve G_{\lambda,F}
        \bigl((\mathbf 1_F R\omega_F)\mu\bigr)}
    \qquad\text{q.e.}
\]
Taking adjoints in the integral identity of Proposition~\ref{prp:damped-transport-identities} gives
\[
    (\mathsf U_t^{R\mu})^*
    =
    \mathsf U_t^*
    -
    \int_0^t
        \mathsf U_s^*R(\mathsf X_s)(\mathsf U_{s,t}^{R\mu})^*
        \,dA_s^\mu .
\]
Substituting this identity into the definition of $\omega_F$ and applying Theorem~\ref{thm:vector-revuz-resolvent} to the first term, we obtain, for q.e.\ $x\in X$,
\begin{align*}
    \omega_F(x)
    &=
    \mathbb E_x\!\left[
        \int_0^{\tau_F} e^{-\lambda t}
        (\mathsf U_t^{R\mu})^*\alpha(\mathsf X_t)\,dt
    \right] \\
    &=
    \widetilde{\breve G_{\lambda,F}(\alpha m)}(x)
    -
    \mathbb E_x\!\left[
        \int_0^{\tau_F} e^{-\lambda t}
        \left(
            \int_0^t
            \mathsf U_s^*R(\mathsf X_s)
            (\mathsf U_{s,t}^{R\mu})^*\alpha(\mathsf X_t)
            \,dA_s^\mu
        \right)dt
    \right]. \tag{$*$}
\end{align*}
For the second term of $(*)$, we interchange the order of integration and then apply the strong Markov property to the integral over $[s,\tau_F)$. For q.e.\ $x\in X$, this gives
\[
\begin{aligned}
    &\mathbb E_x\!\left[
        \int_0^{\tau_F} e^{-\lambda t}
        \left(
            \int_0^t
            \mathsf U_s^*R(\mathsf X_s)
            (\mathsf U_{s,t}^{R\mu})^*\alpha(\mathsf X_t)
            \,dA_s^\mu
        \right)dt
    \right] \\
    &\quad =
    \mathbb E_x\!\left[
        \int_0^{\tau_F} e^{-\lambda s}
        \mathsf U_s^*R(\mathsf X_s)
        \left(
            \int_s^{\tau_F} e^{-\lambda(t-s)}
            (\mathsf U_{s,t}^{R\mu})^*\alpha(\mathsf X_t)\,dt
        \right)dA_s^\mu
    \right] \\
    &\quad =
    \mathbb E_x\!\left[
        \int_0^{\tau_F} e^{-\lambda s}
        \mathsf U_s^*R(\mathsf X_s)\omega_F(\mathsf X_s)\,dA_s^\mu
    \right].
\end{aligned}
\]
Applying the stopped resolvent formula to the second term of $(*)$, we obtain
\[
    \omega_F
    =
    \widetilde{\breve G_{\lambda,F}(\alpha m)}
    -
    \widetilde{\breve G_{\lambda,F}
        \bigl((\mathbf 1_F R\omega_F)\mu\bigr)}
    \qquad\text{q.e.}
\]
\item We next identify the full resolvent. We show that
\[
    \breve R_\lambda^{\mathsf U^{R\mu}}\alpha
    =
    \widetilde{\breve G_\lambda^{R\mu}\alpha}
    \qquad\text{q.e.}
\]
for every bounded measurable section $\alpha\in\L^2(m;\mathcal H)$ and every sufficiently large $\lambda>0$. Since $\|R\|_{\mathrm{op}}\mu$ is smooth, choose a nest $(F_n)$ such that $\mathbf 1_{F_n}\|R\|_{\mathrm{op}}\mu$ is of finite energy integral for every $n$. For q.e.\ $x\in X$, we have $\tau_{F_n}\uparrow\zeta$ $\mathbb P_x$-a.s. Hence
\[
    \omega_{F_n}(x)
    \longrightarrow
    \mathbb E_x\!\left[
        \int_0^\zeta e^{-\lambda t}
        (\mathsf U_t^{R\mu})^*\alpha(\mathsf X_t)\,dt
    \right]
    =
    \breve R_\lambda^{\mathsf U^{R\mu}}\alpha(x)
    \qquad\text{for q.e.\ }x\in X.
\]
On the other hand, the variational identity from the preceding step gives
\[
    \breve{\mathcal E}^{R\mu}_\lambda
    \bigl(\breve G_\lambda^{R\mu}\alpha-\omega_{F_n},\beta\bigr)
    =0
    \qquad\text{for every }\beta\in\D(\breve{\mathcal E}^{R\mu})_{F_n}.
\]
Thus $\omega_{F_n}$ is the $\breve{\mathcal E}^{R\mu}_\lambda$-orthogonal projection of $\breve G_\lambda^{R\mu}\alpha$ onto $\D(\breve{\mathcal E}^{R\mu})_{F_n}$. Since these subspaces are increasing and their union is dense by Lemma~\ref{lem:perturbed-form-nest-density}, it follows that
\[
    \omega_{F_n}\longrightarrow\breve G_\lambda^{R\mu}\alpha
    \qquad\text{in the perturbed form norm}.
\]
By the closed graph theorem, the inclusion $\D(\breve{\mathcal E}^{R\mu})\hookrightarrow\D(\breve{\mathcal E})$ is continuous for the respective form norms. Thus the convergence also holds in the $\breve{\mathcal E}_1$-norm, and along a subsequence $\omega_{F_n}$ converges q.e.\ to $\widetilde{\breve G_\lambda^{R\mu}\alpha}$. Comparing this with the probabilistic limit gives the asserted identity. In particular, $\mathsf U^{R\mu}$ is associated with $\breve{\mathcal E}^{R\mu}$.
\item To conclude proper association, it remains to show that, for every $t>0$ and every bounded measurable section $\alpha\in\L^2(m;\mathcal H)$,
\[
    \breve p_t^{\mathsf U^{R\mu}}\alpha(x)
    =
    \widetilde{\breve T_t^{R\mu}\alpha}(x)
    \qquad\text{for q.e.\ }x\in X.
\]
Fix $t>0$ and a bounded measurable section $\alpha\in\L^2(m;\mathcal H)$. We apply the full resolvent identity to $\breve p_t^{\mathsf U^{R\mu}}\alpha$, which is bounded by the growth bound and represents $\breve T_t^{R\mu}\alpha$ by association. Since the transport semigroup commutes with its resolvent, we obtain, for every sufficiently large integer $n$,
\[
    \breve p_t^{\mathsf U^{R\mu}}
        \bigl(n\breve R_n^{\mathsf U^{R\mu}}\alpha\bigr)
    =
    n\breve R_n^{\mathsf U^{R\mu}}
        \bigl(\breve p_t^{\mathsf U^{R\mu}}\alpha\bigr)
    =
    \widetilde{n\breve G_n^{R\mu}
        \bigl(\breve T_t^{R\mu}\alpha\bigr)}
    \qquad\text{q.e.}
\]

We first pass to the limit on the left. The sections $n\breve R_n^{\mathsf U^{R\mu}}\alpha$ are uniformly bounded by the growth bound and represent $n\breve G_n^{R\mu}\alpha$, which converge to $\alpha$ in $\L^2$. Choose a subsequence converging to $\alpha$ $m$-a.e. By proper association of the scalar process, for q.e.\ $x$ this convergence holds at $\mathsf X_t$ on $\{t<\zeta\}$, $\mathbb P_x$-a.s. Since $\|\mathsf U_t^{R\mu}\|_{\mathrm{op}}$ is integrable, dominated convergence shows that the left-hand side converges to $\breve p_t^{\mathsf U^{R\mu}}\alpha(x)$ along this subsequence.

On the right, $\breve T_t^{R\mu}\alpha\in\D(\breve{\mathcal E}^{R\mu})$ since $t>0$. Hence the spectral theorem gives
\[
    n\breve G_n^{R\mu}\bigl(\breve T_t^{R\mu}\alpha\bigr)
    \longrightarrow \breve T_t^{R\mu}\alpha
    \qquad\text{in the perturbed form norm}.
\]
As above, this convergence also holds in the $\breve{\mathcal E}_1$-norm, so along a further subsequence the quasi-continuous representatives converge q.e.\ to $\widetilde{\breve T_t^{R\mu}\alpha}$. Comparing the two limits gives
\[
    \breve p_t^{\mathsf U^{R\mu}}\alpha
    =
    \widetilde{\breve T_t^{R\mu}\alpha}
    \qquad\text{q.e.},
\]
which proves proper association.
\end{itemize}
\end{proof}

\subsection{Restriction transport}\label{subsec:restriction-transports}
Let $\widetilde{\mathcal K}$ be a compatible quasi-Hilbert subbundle of $\mathcal H^{\breve{\mathcal E}}$ as in Section~\ref{subsec:restriction-forms}, and let $P$ denote the fiberwise orthogonal projection onto $\widetilde{\mathcal K}$.
\begin{dfn}
A contractive quasi-right transport of $\mathsf X$ on $\widetilde{\mathcal K}$ that is properly associated with $\breve{\mathcal E}^{\widetilde{\mathcal K}}$ is called the restriction transport of $\mathsf U$ to $\widetilde{\mathcal K}$. We denote it by $\mathsf U^{\mathcal K}$.
\end{dfn}
By Proposition~\ref{prp:restriction-kato-form} and Theorem~\ref{thm:multiplicative-transport-existence}, such a transport exists; Theorem~\ref{thm:multiplicative-transport-uniqueness} gives uniqueness. Fix any $\widetilde{\mathcal K}$-detecting measure $\mu$, whose existence follows from Theorem~\ref{thm:detecting-measure-existence}. For every integer $n\geq1$, let $\mathsf U^{(n)}$ be the damped transport of $\mathsf U$ corresponding to the smooth potential $n(I-P)\mu$. In other words (by Theorem~\ref{thm:damped-transport-existence}) it is the (unique) solution of
\[
\mathsf U_t^{(n)}=\mathsf U_t-n\int_0^t\mathsf U_{s,t}(I-P(\mathsf X_s))^*\mathsf U_s^{(n)}\,dA_s^\mu,\qquad 0\leq t<\zeta,
\]
for q.e.\ starting point, almost surely. We set $\mathsf U_t^{(n)}=0$ for $t\geq\zeta$.
\begin{thm}[Weak convergence to restriction transport]\label{thm:restriction-transport-limit}
For every $t>0$ and every $\alpha\in\L^2(m;\mathcal H^{\breve{\mathcal E}})$, we have
\[
    (\mathsf U_t^{(n)})^*\bigl(\alpha(\mathsf X_t)\bigr)
    \rightharpoonup
    (\mathsf U_t^{\mathcal K})^*\bigl((P\alpha)(\mathsf X_t)\bigr)
    \quad\text{in }\L^2\!\left(
        \Omega,\mathcal F_t^{\mathsf X},\mathbb P_m;
        \mathcal H_{\mathsf X_0}^{\breve{\mathcal E}}
    \right)
    \quad\text{as }n\to\infty.
\]
\end{thm}
\begin{proof}\;
    \begin{itemize}
\item We first check that $\mathsf U^{(n)}$ is contractive. Suppose first that $\mu$ belongs to the Kato class. Since $\|n(I-P)\|_{\mathrm{op}}\leq n$, the remark immediately preceding Proposition~\ref{prp:damped-transport-identities} gives the integrability required by Theorem~\ref{thm:damped-transport-association}. Thus $\mathsf U^{(n)}$ is properly associated with $\breve{\mathcal E}^{n(I-P)\mu}$. Since $n(I-P)$ is nonnegative, this form satisfies the Kato inequality. Therefore, by Theorems~\ref{thm:multiplicative-transport-existence} and~\ref{thm:multiplicative-transport-uniqueness}, its associated transport is contractive. For a general smooth measure $\mu$, choose a nest $(F_k)$ such that $\mu_k:=\mathbf1_{F_k}\mu$ belongs to the Kato class for every $k$ \cite[Theorem~2.4]{AlbeverioMa1992}. By the preceding argument, $\mathsf U^{n(I-P)\mu_k}$ is contractive. Let $\tau_k$ be the first exit time from $F_k$. Since $dA_s^{\mu_k}=\mathbf1_{F_k}(\mathsf X_s)\,dA_s^\mu$, we have $A_t^{\mu_k}=A_t^\mu$ for every $0\leq t<\tau_k$. Thus the damped transports $\mathsf U^{(n)}$ and $\mathsf U^{n(I-P)\mu_k}$  satisfy the same integral equation on $[0,\tau_k)$. Pathwise uniqueness in Theorem~\ref{thm:damped-transport-existence} therefore gives
\[
    \mathsf U_t^{(n)}=\mathsf U_t^{n(I-P)\mu_k},
    \qquad 0\leq t<\tau_k.
\]
As $\tau_k\uparrow\zeta$ almost surely for q.e.\ starting point, it follows that $\|\mathsf U_t^{(n)}\|_{\mathrm{op}}\leq1$ for all $t<\zeta$ almost surely for q.e.\ starting point.
\item Fix $t>0$ and $\alpha\in\L^2(m;\mathcal H^{\breve{\mathcal E}})$. Let $F=\prod_{i=1}^r f_i(\mathsf X_{t_i})$, where $0<t_1<\cdots<t_r<t$ and $f_1,\ldots,f_r$ are bounded Borel functions. We next prove that, for every $\xi\in\L^2(m;(\mathcal H^{\breve{\mathcal E}})^*)$,
\[
    \lim_{n\to\infty}
    \mathbb E_m\!\left[
        F\left\langle
            (\mathsf U_t^{(n)})^*\bigl(\alpha(\mathsf X_t)\bigr),
            \xi(\mathsf X_0)
        \right\rangle
    \right]
    =
    \mathbb E_m\!\left[
        F\left\langle
            (\mathsf U_t^{\mathcal K})^*\bigl((P\alpha)(\mathsf X_t)\bigr),
            \xi(\mathsf X_0)
        \right\rangle
    \right].
\]
By the contractivity proved above, Theorem~\ref{thm:damped-transport-association} applies to $\mathsf U^{(n)}$. Using multiplicativity and conditioning successively at $t_r,\ldots,t_1$, we obtain
\[
\begin{aligned}
    &\mathbb E_m\!\left[
        F\left\langle
            (\mathsf U_t^{(n)})^*\bigl(\alpha(\mathsf X_t)\bigr),
            \xi(\mathsf X_0)
        \right\rangle
    \right] \\
    &\quad =
    \left(
        \breve T_{t_1}^{n(I-P)\mu}
        \Bigl(f_1\breve T_{t_2-t_1}^{n(I-P)\mu}
        \bigl(\cdots f_r\breve T_{t-t_r}^{n(I-P)\mu}\alpha\bigr)\Bigr),
        \xi
    \right)_2 \\
    &\quad \longrightarrow
    \left(
        \breve T_{t_1}^{\widetilde{\mathcal K}}
        \Bigl(f_1\breve T_{t_2-t_1}^{\widetilde{\mathcal K}}
        \bigl(\cdots f_r\breve T_{t-t_r}^{\widetilde{\mathcal K}}P\alpha\bigr)\Bigr),
        \xi
    \right)_2 \\
    &\quad =
    \mathbb E_m\!\left[
        F\left\langle
            (\mathsf U_t^{\mathcal K})^*\bigl((P\alpha)(\mathsf X_t)\bigr),
            \xi(\mathsf X_0)
        \right\rangle
    \right].
\end{aligned}
\]
Here we used the monotone convergence of $\breve{\mathcal E}^{n(I-P)\mu}$ to $\breve{\mathcal E}^{\widetilde{\mathcal K}}$ from Theorem~\ref{thm:restriction-as-supremum}, which implies strong convergence of the associated semigroups \cite[Theorems~4.1--4.2]{Simon1978}.
    \end{itemize}
Note that the same calculation applies to cylinder functions $F$ with factors at times $0$ and $t$, since these factors can be absorbed into $\xi$ and $\alpha$, respectively. By contractivity of the transports and $m$-symmetry of $\mathsf X$, we have
\[
    \mathbb E_m\!\left[
        \left\|(\mathsf U_t^{(n)})^*\bigl(\alpha(\mathsf X_t)\bigr)\right\|^2
    \right]
    \leq \mathbb E_m\!\left[|\alpha(\mathsf X_t)|^2\right]
    \leq \|\alpha\|_2^2,
\]
and the same bound holds for $(\mathsf U_t^{\mathcal K})^*((P\alpha)(\mathsf X_t))$. Since the finite linear combinations of $F\xi(\mathsf X_0)$ with $F$ a bounded cylinder function are dense in $\L^2(\Omega,\mathcal F_t^{\mathsf X},\mathbb P_m;(\mathcal H_{\mathsf X_0}^{\breve{\mathcal E}})^*)$, the uniform bound extends the convergence of the pairings to every element of this space. This proves the asserted weak convergence.
\end{proof}
\begin{ex}
The penalized transports $\mathsf U_t^{(n)}$ need not converge pathwise to the restriction transport $\mathsf U_t^{\mathcal K}$ simultaneously for all $0<t<\zeta$. Indeed, consider the reflecting Brownian motion on $[0,\infty)$ started at $x>0$ and the quasi-subbundle $\widetilde{\mathcal K}$ from Example~\ref{ex:transport-no-left-limit}.
\end{ex}

\subsection{Jump transport}\label{subsec:jump-transports}
In this subsection, we fix a symmetric $\sigma$-finite Borel
measure $j$ on $X\times X$ satisfying
\[
    j(\operatorname{diag}_X)=0,
    \qquad
    j(N\times X)=0
\]
for every exceptional Borel set $N\subset X$,
together with a connection rule $T$ with respect to $j$,
as in Definition~\ref{dfn:connection-rule}.
Thus $T$ is a $j$-measurable family of contractions
\[
    T(x,y)\colon
    (\mathcal H_x^{\breve{\mathcal E}})^*
    \longrightarrow
    (\mathcal H_y^{\breve{\mathcal E}})^*
\]
satisfying $T(y,x)=T(x,y)^\dagger$
for $j$-a.e.\ $(x,y)$.
We assume that $\mathcal E^j$ and
$\breve{\mathcal E}^{j,T}$ are densely defined on
$\L^2(m)$ and $\L^2(m;\mathcal H)$, respectively,
and that $\mathcal E^j$ is quasi-regular.
By Proposition~\ref{prp:jump-perturbation-kato},
$\mathcal E^j$ is a quasi-regular Dirichlet form and
$\breve{\mathcal E}^{j,T}$ is a Kato form with respect
to $\mathcal E^j$.
\begin{dfn}
Let $(\mathsf X_t^j)_{t\geq0}$ be a symmetric right
process properly associated with $\mathcal E^j$.
A contractive quasi-right transport
$(\mathsf U_t^{j,T})_{t\geq0}$ of
$(\mathsf X_t^j)_{t\geq0}$ on
$\mathcal H^{\breve{\mathcal E}^{j,T}}$
that is properly associated with $\breve{\mathcal E}^{j,T}$
is called a jump transport with respect to the
connection rule $T$.
\end{dfn}
Theorems~\ref{thm:multiplicative-transport-existence} and~\ref{thm:multiplicative-transport-uniqueness}, applied with scalar form $\mathcal E^j$, give existence and uniqueness of this transport.
\begin{ex}\label{ex:finite-activity-jump-transport}
    To see how a jump perturbation affects the transport, we give an explicit description of the resulting transport. We assume that
\[
    j(dx,dy)=m(dx)\,k(x,dy),
    \qquad
    M:=\operatorname*{ess\,sup}_{x\in X}k(x,X)<\infty,
\]
where $k$ is a positive kernel on $X$ and the
essential supremum is taken with respect to $m$.
We also assume that $\mathsf X$ has continuous paths before its lifetime.
Set $\tau_0:=0$ and define the successive jump times of $\mathsf X^j$ by
\[
    \tau_{n+1}
    :=
    \inf\{s>\tau_n\mid s<\zeta^j,\ 
        \mathsf X^j_{s-}\ne\mathsf X^j_s\},
    \qquad n\geq0.
\]
For q.e.\ $x\in X$, $\mathbb P_x$-a.s., only finitely many $\tau_n$, $n\geq1$, lie in $(0,t]$ for every $t<\infty$ (\cite[Section~4.5]{FukushimaOshimaTakeda2011}).
The symmetry of $j$ and the contractivity of $T$ give
\[
    \mathcal J(f)\leq2M\|f\|_2^2,
    \qquad
    \breve{\mathcal J}^{T}(\alpha)\leq2M\|\alpha\|_2^2.
\]
Thus Example~\ref{ex:jump-perturbation-equivalent-norms}
applies with $a=0$ and $b=2M$.
The perturbed domains coincide with the original domains,
and the corresponding form norms are equivalent.
In particular, $\mathcal E^j$ and $\mathcal E$ have the
same nests, and $(\mathcal H^{\breve{\mathcal E}},J)$
is also a quasi-realization of $\breve{\mathcal E}^{j,T}$
with respect to $\mathcal E^j$.
Henceforth, we identify
$\mathcal H^{\breve{\mathcal E}^{j,T}}$ with
$\mathcal H^{\breve{\mathcal E}}$ via the unique isometric quasi-Hilbert bundle isomorphism compatible with the realization maps, and regard $(\mathsf U_t^{j,T})_{t\geq0}$ as a transport on $\mathcal H^{\breve{\mathcal E}}$. 
Write $q(x):=k(x,X)$ (the jump rate) and define the operator $K$ on
$\L^2(m;\mathcal H)$ by
\[
    (K\alpha)(x):=\int_X T(x,y)^*\alpha(y)\,k(x,dy).
\]
Then the generator $\breve L^{j,T}$ with respect to the jump-perturbed form $\breve{\mathcal E}^{j,T}$ is obtained by
\[\breve L^{j,T}=\breve L-q+K,\quad\D(\breve L^{j,T})=\D(\breve L).\]
With these preparations, we can write the jump transport in the following form: for q.e.\ $x\in X$, $\mathbb P_x$-a.s., for all $0\leq t<\zeta^j$:
\[
\begin{aligned}
    \mathsf U_t^{j,T}
    ={}&
    \mathsf U_{\tau_n,t}
    T(\mathsf X^j_{\tau_n-},\mathsf X^j_{\tau_n})
    \mathsf U_{\tau_{n-1},\tau_n-}
    \cdots
    T(\mathsf X^j_{\tau_1-},\mathsf X^j_{\tau_1})
    \mathsf U_{0,\tau_1-},
    \qquad \tau_n\leq t<\tau_{n+1}.
\end{aligned}
\]
Here $\mathsf U_{\tau_{i-1},\tau_i-}$ denotes the original
transport along the continuous path from
$\mathsf X^j_{\tau_{i-1}}$ to $\mathsf X^j_{\tau_i-}$. In fact, the probability of making no jump contributes the term $-q$ to the generator, while averaging the transport across a jump against the jump kernel gives $K$.
\end{ex}

\clearpage
\section{On Riemannian manifolds}\label{sec:riemannian-transports}
Let $M$ be a smooth Riemannian manifold with smooth boundary $\partial M$ and let $E\to M$ be a smooth metric vector bundle of rank $r$ equipped with a metric connection $D$. We denote the induced dual connection on $E^*$ again by $D$. We use the notation and conventions of Section~\ref{sec:riemannian-bundle-forms}. Let $(\mathsf B_t)_{t\geq0}$ be a (reflecting) Brownian motion properly associated with $\mathcal E$, with lifetime $\zeta$.
\subsection{Stochastic parallel transport with Neumann boundary conditions}\label{subsec:riemannian-parallel-transport}
\begin{dfn}
A contractive quasi-right transport $(\mathsf U_t)_{t\geq0}$ of $(\mathsf B_t)_{t\geq0}$ on $E$ that is properly associated with the connection energy $\mathcal E^C$ is called a stochastic parallel transport with Neumann boundary conditions.
\end{dfn}
As shown in Section~\ref{subsec:neumann-connection-energy}, $\mathcal E^C$ is a Kato form and $(E,\operatorname{Id}_E)$ is its quasi-realization. Theorems~\ref{thm:multiplicative-transport-existence} and~\ref{thm:multiplicative-transport-uniqueness} therefore guarantee the existence of such a transport and its uniqueness up to $\mathbb P_x$-indistinguishability for q.e.\ $x\in M$. We fix such a transport and write
\[
    \mathsf U_t\colon E_{\mathsf B_0}^*\longrightarrow E_{\mathsf B_t}^*,
    \qquad 0\leq t<\zeta.
\]

This transport agrees, up to $\mathbb P_x$-indistinguishability for q.e.\ $x\in M$, with the classical stochastic parallel transport defined by the horizontal lift; see \cite{Hsu2002} and, for reflecting Brownian motion, \cite[Section~2]{Hsu2002Boundary} and \cite[Section~4]{ArnaudonLi2017}. Indeed, the classical construction is a quasi-right transport associated with the same Neumann connection energy, so Theorem~\ref{thm:multiplicative-transport-uniqueness} applies. In particular, for q.e.\ $x\in M$, $\mathbb P_x$-a.s., $\mathsf U_t$ is an isometry for every $0\leq t<\zeta$.

\subsection{Stochastic parallel transport with mixed boundary conditions}\label{subsec:stochastic-tangent-transport}
Fix an integer $k\geq0$ and take $E=\otimes^kT^*M$ with the connection induced by the Levi--Civita connection, so that $\mathsf U_t$ acts on $\otimes^kTM$. We use the metric to identify tensor powers of the tangent and cotangent bundles.
\begin{dfn}
The restriction transport of $\mathsf U$ to $\otimes^kT_{\mathrm{tan}}M$ is called the stochastic parallel transport with mixed boundary conditions. We denote it by $\mathsf U^{\mathrm{tan}}$, with the degree understood from the bundle, so that
\[
    \mathsf U_t^{\mathrm{tan}}\colon \otimes^kT_{\mathrm{tan},\mathsf B_0}M\longrightarrow \otimes^kT_{\mathrm{tan},\mathsf B_t}M,
    \qquad 0\leq t<\zeta.
\]
Equivalently, on the section side it is the restriction transport to $\otimes^kT_{\mathrm{tan}}^*M$. For $0\leq k\leq\dim M$, the same definition applies with $\otimes^k$ replaced by $\Lambda^k$, and we use the same notation $\mathsf U^{\mathrm{tan}}$.
\end{dfn}
By Section~\ref{subsec:mixed-connection-energy}, this transport is properly associated with the connection energy $\mathcal E^C$ with mixed boundary conditions on the corresponding tensor or exterior power.
Let $P$ be the orthogonal projection onto the corresponding tangential quasi-subbundle, and choose any detecting measure $\mu$. For the damped transports $\mathsf U^{(n)}$ corresponding to $n(I-P)\mu$, Theorem~\ref{thm:restriction-transport-limit} gives weak $\L^2(\mathbb P_m)$ convergence of $(\mathsf U_t^{(n)})^*\alpha(\mathsf B_t)$ to $(\mathsf U_t^{\mathrm{tan}})^*(P\alpha)(\mathsf B_t)$ for each $t>0$ and every $\L^2$ section $\alpha$. The limit is independent of the choice of $\mu$. In particular, one may take $\mu=\sigma$, since the boundary measure $\sigma$ is a detecting measure by Section~\ref{subsec:mixed-connection-energy}.
In this case, the construction can be summarized schematically as
\[
    \mathsf U_t
    \xrightarrow[n(I-P)\sigma]{\text{damping}}
    \mathsf U_t^{(n)}
    \xrightarrow[n\to\infty]{\text{weak limit}}
    \mathsf U_t^{\mathrm{tan}}.
\]
For $k=1$, under the hypotheses of \cite{Hsu2002Boundary}, this transport agrees with the transport obtained from the construction in Sections~3 and~6 of that paper by setting both the interior and boundary damping endomorphisms to zero. Both transports are quasi-right and associated with the same mixed connection energy, so Theorem~\ref{thm:multiplicative-transport-uniqueness} gives $\mathbb P_x$-indistinguishability for q.e.\ $x$, with initial vectors in $T_{\mathrm{tan},x}M$.

\subsection{Ricci-damped transport with mixed boundary conditions}\label{subsec:weitzenbock-ricci-damped-transport}
Fix an integer $0\leq k\leq\dim M$ and let $\mathsf U^{\mathrm{tan}}$ be the stochastic parallel transport with mixed boundary conditions on $\Lambda^kT_{\mathrm{tan}}M$. Throughout this subsection, we assume the hypotheses of the final example in Section~\ref{subsec:hodge-energy-mixed-boundary} (either $M$ is compact, or $k=1$ and $M$ is complete with Ricci curvature bounded from below and convex boundary.)
\begin{dfn}
A quasi-right transport of $\mathsf B$ on $\Lambda^kT_{\mathrm{tan}}^*M$ that is properly associated with the Hodge energy $\mathcal E^H$ on $k$-forms with mixed (absolute) boundary conditions is called the Ricci-damped transport with mixed boundary conditions. We denote it simply by $\mathsf Q$, with the degree understood from the bundle, so that
\[
    \mathsf Q_t\colon \Lambda^kT_{\mathrm{tan},\mathsf B_0}M\longrightarrow \Lambda^kT_{\mathrm{tan},\mathsf B_t}M,
    \qquad 0\leq t<\zeta.
\]
\end{dfn}
By Theorem~\ref{thm:damped-transport-association}, this is the damped transport of $\mathsf U^{\mathrm{tan}}$ corresponding to the smooth potential
\[
    \frac12\mathscr R_k\,\mathrm{vol}+\frac12S_k\,\sigma.
\]
Here $\mathscr R_k$ and $S_k$ are the curvature and boundary endomorphisms defined in Section~\ref{subsec:hodge-energy-mixed-boundary}; in particular, $\mathscr R_1$ is the Ricci endomorphism. This gives existence, and Theorem~\ref{thm:multiplicative-transport-uniqueness} gives uniqueness. Thus the construction can be summarized schematically as
\[
    \mathsf U_t
    \xrightarrow{\text{tangential restriction}}
    \mathsf U_t^{\mathrm{tan}}
    \xrightarrow{\text{damping}}
    \mathsf Q_t.
\]
For $k=1$, under the respective hypotheses of \cite[Theorem~4.2]{Hsu2002Boundary} and \cite[Theorem~5.2]{ArnaudonLi2017}, their classical damped transports are quasi-right and associated with the same Hodge energy. Hence they agree with $\mathsf Q$ up to $\mathbb P_x$-indistinguishability for q.e.\ $x\in M$, with initial vectors in $T_{\mathrm{tan},x}M$.

\clearpage
\renewcommand{\bibname}{References}

\end{document}